\documentclass[10pt,hidelinks]{article}
\providecommand{\keywords}[1]{\textbf{\textit{Key words}} #1}

\usepackage{cancel}

\usepackage[normalem]{ulem}

\usepackage{pstricks-add}
\usepackage{yfonts}
\usepackage{pstricks-add}

\usepackage{pst-func}

\usepackage{bbold}
\usepackage[normalem]{ulem}
\usepackage{etex}
\usepackage[left=2cm,top=2cm,right=2cm,bottom=2cm]{geometry}
\usepackage{amsmath, amsthm, amssymb}
\allowdisplaybreaks
\usepackage{pst-all}
\usepackage[textwidth=0.5in]{todonotes}
\usepackage[rightcaption]{sidecap}

\usepackage{relsize}
\usepackage{float}
\usepackage[caption = false]{subfig}

\usepackage[hypertexnames=false]{hyperref}

\usepackage[scaled=1]{helvet}
\usepackage{titlesec}
\usepackage{multido}
\usepackage{graphicx}
\usepackage{multirow}
\usepackage{blindtext}
\usepackage{color}
\usepackage{lipsum}
\usepackage{mathtools}
\usepackage[vcentermath]{youngtab}
\usepackage{young}
\usepackage[all,ps,dvips,graph]{xy}
\usepackage[misc,geometry]{ifsym}

\usepackage{titlesec}
\usepackage{lineno,hyperref}

\let\OLDthebibliography\thebibliography
\renewcommand\thebibliography[1]{
  \OLDthebibliography{#1}
  \setlength{\parskip}{0pt}
  \setlength{\itemsep}{0pt plus 0.3ex}
}

\titleformat{\section} {\normalfont\scshape \large \centering}{ \thesection}{1em}{}
\definecolor{morado}{rgb}{0.5,0,0.5}

\newcommand{\QQ}{\mathbb Q}

\newcommand{\shape}{\text{shape}}

\newcommand{\FF}{ { \mathbb F}_{\! p}}

\newcommand{\N}{ { \mathbb N}}

\newcommand{\NN}{ { \color{black}{n}}}
\newcommand{\NNone}{ { \color{black}{n_1}}}
\newcommand{\NNtwo}{ { \color{black}{n_2}}}
\newcommand{\RR}{ { \color{black}{r}}}

\newcommand{\JWk}{ {\mathbf{JW}_{\! k}}}
\newcommand{\JWn}{ {\mathbf{JW}_{\! n}}}
\newcommand{\JWntwo}{ {\mathbf{JW}_{\! n_2}}}

\newcommand{\JWtres}{ {\mathbf{JW}_{\! 3}}}

\newcommand{\JWdos}{ {\mathbf{JW}_{2}}}
\newcommand{\JWdoce}{ {\mathbf{JW}_{\! 12}}}

\newcommand{\TLkQ}{ {\mathbb{TL}_k^{\! \mathbb Q}}}
\newcommand{\TLnQ}{ {\mathbb{TL}_n^{\! \mathbb Q}}}
\newcommand{\TLnk}{ {\mathbb{TL}_n^{\! \Bbbk}}}

\newcommand{\TLZpn}{  {\mathbb {TL}}_{n}^{\! {\mathbb Z}_{(p)} }}
\newcommand{\TLZpntwo}{  {\mathbb {TL}}_{n_2}^{\! {\mathbb Z}_{(p)} }}
\newcommand{\TLFpntwo}{  {\mathbb {TL}}_{n_2}^{\! {\mathbb F}_{(p)} }}
\newcommand{\TLQpntwo}{  {\mathbb {TL}}_{n_2}^{\! \mathbb  Q  }}

\newcommand{\TLn}{ {\mathbb{TL}_n}}
\newcommand{\TLndiag}{ {\mathbb{TL}_n^{\! diag}}}
\newcommand{\TLnF}{ {\mathbb{TL}_n^{\! {\mathbb F}_p}}}

\newcommand{\TLntwoF}{ {\mathbb{TL}_{n_2}^{\! {\mathbb F}_p}}}

\newcommand{\TLnZp}{ {\mathbb{TL}_n^{\! {\mathbb Z}_{(p)}}}}

\newcommand{\TLtrestres}{ {\mathbb{TL}_3^{\! {\mathbb F}_3}}}

\newcommand{\A}{\mathcal{A}}
\newcommand{\Z}{\mathbb{Z}}

\newcommand{\Par}{{\rm Par}_n   }
\newcommand{\ParTwo}{{\rm Par}^{\le 2}_n}
\newcommand{\ParTwoNuno}{{\rm Par}^{\le 2}_{n_2}}

\newcommand{\RKLR}{{ \mathcal R}_n }
\newcommand{\RKLRZ}{ {\mathcal R}_{n}^{ \Z_{(p)}}}

\newcommand{\s}{\mathfrak{s}}

\newcommand{\V}{\mathfrak{v}}
\newcommand{\T}{  \mathfrak{t}}
\newcommand{\U}{  \mathfrak{u}}

\newcommand{\aaa}{\mathfrak{a}}

\newcommand{\Si}{\mathfrak{S}}
\newcommand{\std}{{\rm Std}}

\newcommand{\R}{ \mathcal{R}_n}

\newcommand{\UU}{\mathbb{u}}

\newcommand{\UUU}{\mathbb{U}}

\newcommand{\LL}{\mathbf{L}}
\newcommand{\LLL}{ \textswab{L}}  
\newcommand{\EE}{\mathbb{E}}
\newcommand{\EEE}{{\color{black}{ \mathbb{E}}}}
\newcommand{\ee}{\mathbf{e}}

\newcommand{\JM}{ \mathbf{JM} }

\newcommand{\one}{\mathbb{1}}
\newcommand{\sss}{\mathbb{s}}

\newcommand{\greekrho}{\varrho}

\newcommand{\Std}{{\rm Std}}

\newcommand{\ii}{\mathbf{i}}
\newcommand{\jj}{\mathbf{j}}

\modulolinenumbers[5]

\newtheorem{theorem}{Theorem}[section]
\newtheorem{lemma}[theorem]{Lemma}

\newtheorem{definition}[theorem]{Definition}
\newtheorem{corollary}[theorem]{Corollary}

\newtheorem{remark}[theorem]{Remark}

\newenvironment{dem}{\noindent \textit{Proof:} }{\quad \hfill $\square$}

\numberwithin{equation}{section}

\newgray{plomoclaro}{.90}
\newgray{plomooscuro}{.70}
\newgray{blanco}{1}

\begin{document}
\Yvcentermath1
\sidecaptionvpos{figure}{lc}


\title{Seminormal forms for the Temperley-Lieb algebra. }

\author{
  \, Katherine Orme\~no Bast\'ias{\thanks{Supported in part by ANID beca
doctorado nacional 21202166 and EPSRC grant EP/W007509/1
  }} \, and  Steen Ryom-Hansen{\thanks{Supported in part by FONDECYT grant 1221112 
and EPSRC grant EP/W007509/1}}}

\date{\vspace{-5ex}}
\maketitle
\begin{abstract}
Let $\TLnQ$ be the rational Temperley-Lieb algebra, with loop parameter $ 2 $. 
In the first part of the paper we study the seminormal idempotents $ \EEE_{\T} $ for $\TLnQ$ for 
$ \T $ running over two-column standard tableaux. 
Our main result is here a concrete combinatorial construction of $ \EEE_{\T} $ 
using Jones-Wenzl idempotents $ \JWk $ for $\TLkQ$ where $ k \le n $. 

In the second part of the paper
we consider 
the Temperley-Lieb algebra $\TLnF$ over the finite field $ \FF$, where $ p>2$.
The KLR-approach to $\TLnF$ gives rise to an action 
of a symmetric group $ \Si_m $ on $\TLnF$, for some $ m < n $. 
We show that the $ \EEE_{\T} $'s from the first part of the paper are 
simultaneous eigenvectors for the associated Jucys-Murphy elements for $ \Si_m $. This leads to 
a KLR-interpretation of the $p$-Jones-Wenzl idempotent 
$ ^{p}\!\JWn $ for $\TLnF$, that was introduced recently by Burull, Libedinsky and Sentinelli. 
\end{abstract}


\section{Introduction}
In the present paper we introduce a new perspective on both the semisimple and the non-semisimple
representation theory of the Temperley-Lieb algebra $ \TLn$. Some of the ingredients of
this new perspective are very classical and well-known, whereas other ingredients
are based on ideas developed within the last few
years. The unifying element of these ingredients
are {\it seminormal forms}. 

\medskip
In the representation theory of the symmetric group $ \Si_n$, 
the seminormal form has a long history 
going back to
Young's papers a century ago. A milestone in this history was 
Murphy's discovery in the eighties, 
realizing the {\it seminormal form} as common eigenvectors for the Jucys-Murphy elements
of the symmetric group. 
Later on, Mathas showed that Murphy's results hold in the 
general framework of a {\it cellular algebra} $ \A $ endowed with a {\it family of $ \JM$-elements}.
In the present paper we
{\color{black}{use}}
that $ \TLn$ fits into this framework, with
cellular structure given by the diagram basis and $ \JM $-elements induced
from the $ \JM $-element from the symmetric group. 

\medskip
In Mathas' framework there is a dichotomy between the {\color{black}{\lq separated\rq\ }}and
the {\color{black}{\lq unseparated\rq\ }}cases 
and our treatment of $ \TLn$ follows this dichotomy,
with $ \TLnQ$ corresponding to the separated case and $  \TLnF$ to the unseparated case. 
In the separated case, $ \A $ is semisimple and 
there exists a {\it complete} set of {\it orthogonal idempotents}
$ \{ \EE_{\T} \, | \, \T \in \std(\Lambda) \} $ that are common eigenvectors for the $ \JM $-elements. 
In the unseparated case however, $ \A $ is not semisimple and the $ \EE_{\T} $'s are undefined, but
summing over a {\it $p$-class $ [\T]$ of standard tableaux} we have well defined
{\it class idempotents}
in $ \A $, given by 
\begin{equation}
\EE_{[\T]} := \sum_{ \s \in [\T]  } \EE_{\s}  
\end{equation}

\medskip
In general one is especially interested in the {\it primitive idempotents }in $ \A $, since
they correspond to the {\it indecomposable projective modules} for $ \A$, that carry
essential relevant information about the representation theory of $\A$, such as decomposition
numbers, etc. The primitive idempotent for the projective cover of the trivial $ \TLnF$-module,
the {\it $p$-Jones-Wenzl idempotent} $\,  ^{p}\!\JWn $, 
was determined recently
by Burull, Libedinsky and Sentinelli in \cite{BLS}, via a recursive construction involving the
expansion of $ n+1 $ in base $ p$. 
The class idempotents $ \EE_{[\T]} $
are not always primitive, 
but in our final Theorem \ref{finalTheo}
we show that $ ^{p}\!\JWn $ can in fact be built from
the class idempotents, via a  
recursive construction in KLR-theory.

\medskip

In the separated case, our main results are  
Theorem \ref{mentionedabove} and Corollary \ref{finalcorsection2} that together give a 
realization of the idempotents $ \EE_{\T} $ in terms 
of a concrete diagrammatic construction involving 
{\it Jones-Wenzl} idempotents
$ \JWk $ for $\TLkQ$ where $ k \le n $. A key ingredient for the proof of these results is
contained in our Theorem \ref{YSFfirst}, 
giving a
close connection between the Jones-Wenzl idempotents
and the seminormal form, in which 
the well-known recursive formula for $ \JWn$
\begin{equation}
  \raisebox{-.35\height}{\includegraphics[scale=0.7]{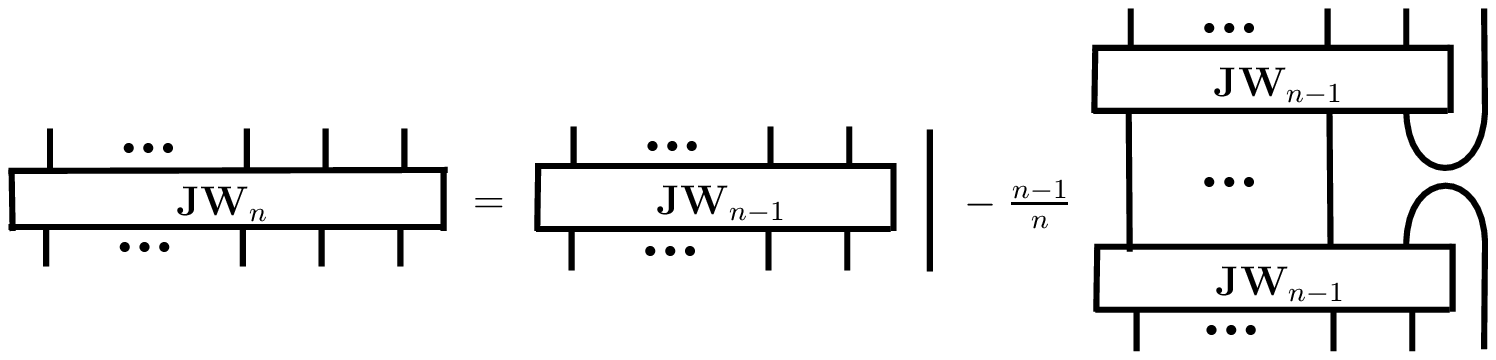}}
\end{equation}
becomes
{\color{black}{a special case of the classical formula for the action of the symmetric group on the seminormal form,
see Remark \ref{YSFremark}}}.

\medskip
Our treatment of the unseparated case $ \TLnF$ relies heavily on the KLR-algebra (Khovanov-Lauda-Rouquier)
approach to $ \TLnF$. The KLR-algebra was introduced independently
by Khovanov-Lauda and Rouquier in
\cite{KhovanovLauda} and \cite{Rouquier}, 
in order to categorify the positive part of the type
$ A$ quantum
group. In seminal papers by Brundan-Kleshchev and Rouquier, see \cite{brundan-klesc}
and \cite{Rouquier}, an isomorphism $  \R \cong  \FF \Si_n  $ was established (in fact 
in the greater generality of cyclotomic Hecke algebras).  
Hu and Mathas gave in \cite{hu-mathas2} a new
simpler proof of this isomorphism
{\it using seminormal forms}, and via this they were able to lift it to an
isomorphism $\RKLRZ \cong   \Z_{(p)} \Si_n   $, where $ \RKLRZ $
is an integral version of $ \RKLR $, defined over the localization $ \Z_{(p)} $ of $ \Z $ at the prime $p$.
This isomorphism and its proof
{\color{black}{are}} important to us. 
It firstly induces an isomorphism $ \RKLRZ \!  / {\cal I }_n \cong  \TLnZp  $ where $ {\cal I }_n $
is an {\color{black}{explicitly given graded}}
ideal of $ \RKLRZ $, generalizing a result by D. Plaza and the second author, see
Theorem \ref{steendavid} and
\cite{PlazaRyom}.

\medskip
Under the isomorphism $  \RKLRZ \! /  {\cal I }_n   \cong   \TLnZp $, the KLR-generator
$ e(\ii) $ for $ \ii =(0,-1,-2, \ldots, -n+1) $ corresponds to the class idempotent $ \ee= \EE_{[ \T_n] } $
where $ \T_n $ is the unique one-column standard tableau of length $ n $, and so we consider the
idempotent truncated subalgebra $ \ee \TLnZp \ee $ of $ \TLnZp $. It contains block intertwining
elements $ \UUU_i $, that are represented diagrammatically by
{\color{black}{\lq diamonds\rq\ }}as follows, for $ p=3$. 
\begin{align}
& \UUU_1 = \raisebox{-.5\height}{\includegraphics[scale=0.75]{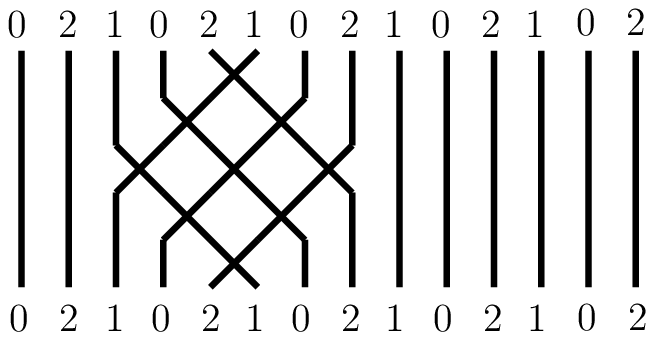}} & 
& \UUU_2 = \raisebox{-.5\height}{\includegraphics[scale=0.75]{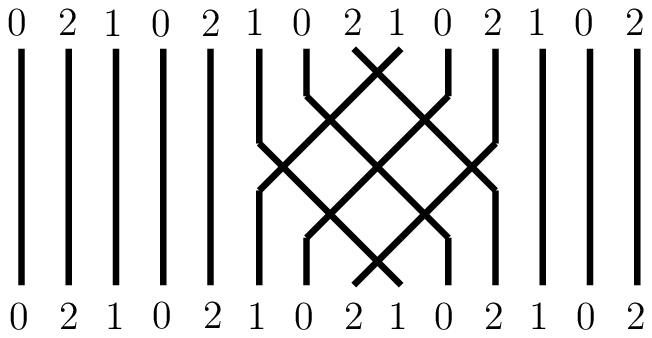}} 
\end{align} 
Similar elements have already been been considered before for example
in \cite{KMR}, \cite{LiPl}, \cite{Lo} and \cite{LPR},
but only in the context of the original KLR-algebra $ \R $ defined over a field.

\medskip
Our main results in the 
unseparated case 
revolve around the
action of the $ \UUU_i $'s on 
the seminormal form for $  \TLnQ$, given by Hu-Mathas' isomorphism. The key result 
is our Theorem \ref{YSFthird},
which establishes a formula reminiscent of the classical formula for the
action of the symmetric group on 
Young's seminormal form.
The fact that such a formula is possible
in the representation theory of $ \TLnZp $, and hence also 
of 
$ \TLnF$, 
could surely not have been conceived before
the introduction of KLR-algebra in representation theory, in particular before Hu-Mathas'
proof of the isomorphism $\RKLRZ \cong   \Z_{(p)} \Si_n   $.

\medskip
As a first consequence of Theorem \ref{YSFthird} we obtain in
Corollary \ref{thereisaninjection}
an injection
$\iota_{KLR}:   {\mathbb {TL}}_{n_2}^{\! {\mathbb Z}_{(p)} } \rightarrow \TLZpn $ for
$ n_2 $ a concretely given number $ n_2 < n$. 
The small Temperley-Lieb algebra $ {\mathbb {TL}}_{n_2}^{\! {\mathbb Z}_{(p)} } $
is of course endowed with its own
$ \JM$-elements $ \LLL_i $, with associated seminormal form idempotents $ \{ \EEE_{s} \} $
that apriori are unrelated to the original seminormal form idempotents $  \{ \EE_{\T} \} $
for $\TLnQ$. But nevertheless we establish in Theorem \ref{indutionseminormalKLR}, 
Corollary \ref{idempotentJMcor} and
Corollary \ref{idempotentJMcordos}
another series of
{\color{black}{interesting}}
facts, showing that the
$   \EE_{\T}  $'s are eigenvectors for the $ \LLL_i $'s, and that the product
$  \EEE_{\s} \EE_{\T} $ is given by a simple formula. 

\medskip

These last results are the essence of our recursive construction of $   ^{p}\!\JWn $.
We obtain a chain 
\begin{equation} 
0 \subseteq   {\mathbb {TL}}_{n^{k-1}_2}^{\! {\mathbb Z}_{(p)} } \subseteq
  {\mathbb {TL}}_{n^{k-2}_2}^{\! {\mathbb Z}_{(p)} } \subseteq  \cdots \subseteq  
  {\mathbb {TL}}_{n_2^0}^{\! {\mathbb Z}_{(p)} }  \subseteq
  \mathbb {TL}\strut_{n}^{\! {\mathbb Z}_{(p)} }
\end{equation}  
of subalgebras, and let $ \EEE_{[\T_{n^{k-1}_2}]} $ be the class idempotent for 
$ {\mathbb {TL}}_{n_2^{k-1}}^{\! {\mathbb Z}_{(p)} }  $, considered as an element of 
$   {\mathbb {TL}}_{n}^{\! {\mathbb Z}_{(p)} } $,
{\color{black}{via this chain.}}
In our final Theorem \ref{finalTheo}
we show that $\,  ^{p}\!\JWn $ is {\color{black}{equal to
$ {\mathbb {TL}}_{n_2^{k-1}}^{\! {\mathbb Z}_{(p)} }  $.}}

\medskip
The organization of our paper is as follows. In section \ref{basic notions}
we fix notation and recall some basic results concerning the objects that are studied in the paper, that is
the Temperley-Lieb algebra $\TLn$, the group algebra of the symmetric group, and so on.
In section \ref{separated case}, we construct for each two-column standard tableau $ \T $ an
element $ \EE_{\T}^{\prime} \in \TLnQ$ and show that $  \EE_{\T}^{\prime} = \EE_{\T}$.
In section \ref{The unseparated case}, we initiate our study of the unseparated case. We recall
the construction of $\,  ^{p}\!\JWn $ and give various examples that motivate and illustrate the main
results of the last section. In section \ref{The integral KLR-algebra}, we recall the KLR algebra
$ \R $ and in particular Hu and Mathas' integral version $\RKLRZ $ of $ \R$. We also recall the basic
ingredients in Hu and Mathas' proof of the isomorphism $\RKLRZ \cong   \Z_{(p)} \Si_n   $. They
involve a series of formulas for the action of the KLR-generators on the seminormal basis.

In section \ref{Seminormal form for}, we prove our main results concerning the unseparated
case,
starting with Theorem \ref{YSFthird}, which describes the action of the $ \UUU_i $ on the seminormal
basis. We found the results of this section with 
considerable aid from the
Sage computer system. The proof of Theorem \ref{YSFthird} is a lengthy calculation involving Hu-Mathas' formulas. 
We firmly believe that it is possible to generalize all our results of the present paper to hold 
for the 
Tempeley-Lieb algebra 
at loop parameter $ q + q^{-1} $ and have in fact already made some progress in this direction.
On the other hand,
we also realized that the generalization of Theorem \ref{YSFthird} to loop parameter
$ q+ q^{-1} $ involves substantially more
calculations
and we therefore decided to develop
our results firstly in the loop parameter 2 setting.

\medskip
We thank N. Libedinsky, D. Plaza and S. Griffeth for useful comments on a previous version
of this work. We specially thank G. Burull for explaining us the construction of $ ^{p}\!\JWn $.
{\color{black}{Finally, but not least, we thank the two anonymous referees for their interested and detailed reports that helped us
improve and clarify the article. }}

\section{Basic notions}\label{basic notions}
In this section we fix notation and recall the basic notions related to the
symmetric group and 
Temperley-Lieb algebra and their representation theories.

\medskip
Let $ n $ be a positive integer. A partition of $ n $ of {\it length} $ k $ is a weakly
decreasing sequence $ \lambda = (\lambda_1, \lambda_2,
\ldots, \lambda_k ) $ of positive integers with total sum $  n$. 
The set of partitions of $ n $ is denoted
$ \Par$. A partition $ \lambda  = (\lambda_1, \lambda_2,
\ldots, \lambda_k )\in \Par $ is represented diagrammatically by its {\it Young diagram}, which consists 
of $ k $, left aligned, rows of nodes with $ \lambda_1 $ nodes in the first row, $ \lambda_2 $
nodes in the second row, and so on. For example, 
$ \lambda = (7,4,2^2) = (7,4,2,2) \in {\rm Par}_{15} $ is represented by 
the following Young diagram
\begin{equation}\label{dib47}
\lambda = \raisebox{-.5\height}{\includegraphics[scale=0.7]{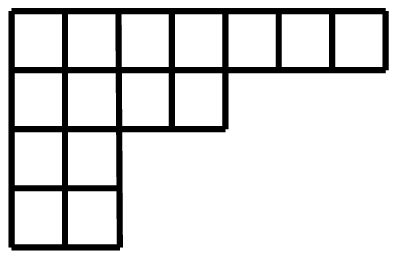}}  
\end{equation}
As in \eqref{dib47}, we shall identify $ \lambda \in \Par$ with its Young diagram.
The set of partitions of $n$ whose Young diagram has less than {\color{black}{$k$}} columns
is denoted $ {\rm Par}_n^{\le k } $. For $ \lambda, \mu \in \Par $ we write 
$ \lambda \trianglelefteq \mu $
if $ \lambda_1 \le \mu_1 $, $ \lambda_1 + \lambda_2 \le \mu_1 + \mu_2 $,
$ \lambda_1 + \lambda_2 + \lambda_3 \le \mu_1 + \mu_2 +\mu_3 $,
and so on. This is the
{\it dominance order} on $ \Par$. For $ \lambda \in \Par $, a $ \lambda$-tableau $ \T $ is a filling of
the nodes of the Young diagram of $ \lambda$ using the numbers
from the set $ \{ 1,2,\ldots, n \}$, with each number occurring 
once. A $ \lambda$-tableau $ \T $ is said to be {\it standard} if the numbers from
$ \{ 1,2,\ldots, n \}$ appear increasingly from left to right along the rows of $ \lambda $, and increasing from top to
bottom along the columns of $ \lambda $. The set of standard $ \lambda$-tableaux is denoted $ \std(\lambda)$. 
Here is an example of a standard $ \lambda$-tableau $ \T $, for $ \lambda $ as in \eqref{dib47}
 \begin{equation}\label{dib48}
\T = \raisebox{-.5\height}{\includegraphics[scale=0.7]{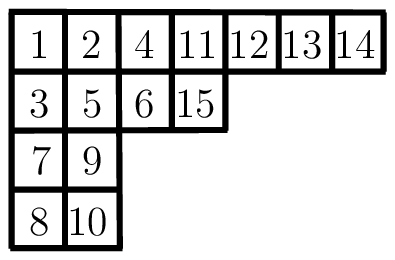}}  
\end{equation}
 We define $ \shape(\T) = \lambda $ if $ \T \in \Std(\lambda) $. For $ \T $ any standard tableau we define
 $ \T |_{\le m } $ to be the tableau obtained from $ \T $ by deleting all the nodes that contain numbers
 strictly larger than $ m $. 
 We then extend the dominance order to standard tableaux
 via $ \s \trianglelefteq \T $ if $ \shape( \s |_{\le m } ) \trianglelefteq \shape( \T |_{\le m } ) $ for
 all $ m=1,2,\ldots, n$. 

\medskip
For $\lambda \in \Par$ we define $ \T^{\lambda} \in  \Std(\lambda) $ as the {\it row-reading} tableau,
and similarly we define $ \T_{\lambda} \in  \Std(\lambda) $ as the {\it column-reading} tableau.
 When restricted to $ \std(\lambda) $, we have that $ \T^{\lambda} $ is the unique maximal tableau
 and $ \T_{\lambda} $ is the unique minimal tableau, with respect to $\trianglelefteq $.
For example, for $ \lambda $ as in \eqref{dib47} we have 
 \begin{equation}\label{dib48second}
   \T^{\lambda} = \raisebox{-.5\height}{\includegraphics[scale=0.7]{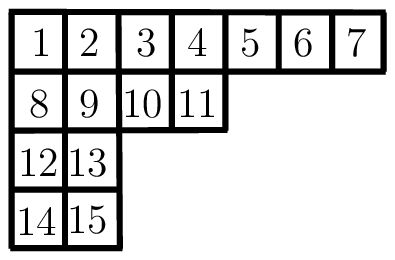}}, \, \, \, \, \, \, \, \, 
\T_{\lambda} = \raisebox{-.5\height}{\includegraphics[scale=0.7]{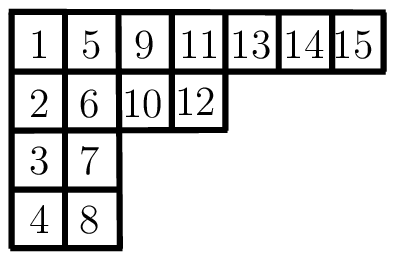}}     
\end{equation}

 \medskip
{\color{black}{ For $ n= 0 $ we use the conventions that $  {{\rm Par}_0   } = \{ \emptyset  \} $ and
 $ \std(\emptyset) = \emptyset $, and so $|  {\rm Par}_0  |  = | \std(\emptyset)| =1  $.}}

\medskip
The Temperley-Lieb algebra was introduced in the seventies from motivations in statistical mechanics.
{\color{black}{Since then it 
has been generalized in several interesting ways and 
has been shown to be related to many
areas of mathematics as well, including knot theory, categorification theory and Soergel bimodules, 
see for example \cite{GL2}, 
\cite{Lo2}, \cite{LPR}, \cite{Lobos-Ryom-Hansen}, \cite{PMartin}, \cite{PlazaRyom},
\cite{Stro}}}. 
In this paper we shall use the variation of the Temperley-Lieb algebra that has loop parameter
equal to $2$. It is defined as follows.  
\begin{definition}\label{defTL}
  The Temperley-Lieb algebra $ \TLn$ is the associative unitary $ \Z $-algebra on 
  generators $ \UU_1, \UU_2, \ldots, \UU_{n-1} $ subject to the relations
	\begin{align}
\label{eq oneTL}	\UU_i^2	& =   2\UU_i,  & &  \mbox{if }  1\leq i < n \\
\label{eq twoTL}	\UU_i\UU_j\UU_i & =\UU_i,      &  & \mbox{if } |i-j|=1   \\
\label{eq threeTL}	\UU_i\UU_j& = \UU_j\UU_i,      & &   \mbox{if } |i-j|>1
	\end{align}
{\color{black}{For $ n=0 $ or $ n=1 $ we define $  {\mathbb {TL}}_{n} := \Z$. }}
\end{definition}
For $ \Bbbk $ a commutative ring 
we shall also consider
the {\it specialized version} $ \TLnk$ of $ \TLn $, defined as 
\begin{equation}
 \TLnk:=  \TLn \otimes_{\Z} \Bbbk 
\end{equation}
Here we are mostly interested in the cases 
where 
$  \Bbbk $ is the rational field $ \QQ $, 
the finite field with
$ p $ elements $ \FF $ or the localization $ \Z_{(p)} $ of $ \Z $ at the prime $ p$.
The corresponding Temperley-Lieb algebras are  
$ \TLnQ $, $  \TLnF $ and $ \TLZpn$.

\medskip
A well-known and important feature of $ \TLn $ is the fact that it is a diagram algebra. Concretely, 
$ \TLn $ is isomorphic to the diagrammatically defined algebra $ \TLndiag $ with basis given by {\it non-crossing
planar matchings} of $ n $ northern points of a(n invisible) rectangle with $ n $ southern points of the rectangle. 
Here are three examples for $ n= 5$. 
\begin{equation}\label{dib1}
\raisebox{-.5\height}{\includegraphics[scale=0.75]{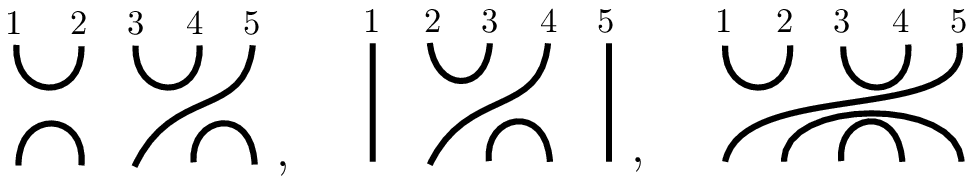}}  
\end{equation}
We refer to such matchings as {\it Temperley-Lieb diagrams}. 
For two Temperley-Lieb diagrams $ D_1 $ and $ D_2 $ the product $ D_1 D_2 $ in $ \TLndiag $
is given by concatenation
with $ D_1 $ on top of $ D_2 $. For example, choosing $ D_1 $ and $ D_2 $ as the first two diagrams in
\eqref{dib1} we have that
\begin{equation}\label{dib2}
D_1 D_2 \, \, =\, \,  \raisebox{-.4\height}{\includegraphics[scale=0.75]{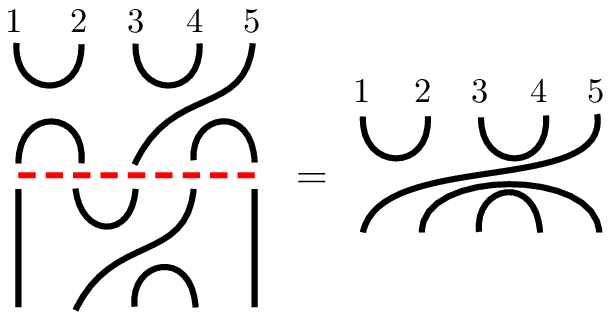}}   \, \, =\, \,D_3
\end{equation}
where $ D_3 $ is the third diagram of \eqref{dib1}. This concatenation product may give rise to
diagrams with internal loops. Each internal loop is removed from
the diagram, and the resulting diagram is multiplied by the scalar $ 2 \in \Z $.
For example, if $ D_1 $ and $ D_3 $ are as above, we have that 
\begin{equation}\label{dib3}
D_1 D_3 \, \, =\, \,  \raisebox{-.4\height}{\includegraphics[scale=0.75]{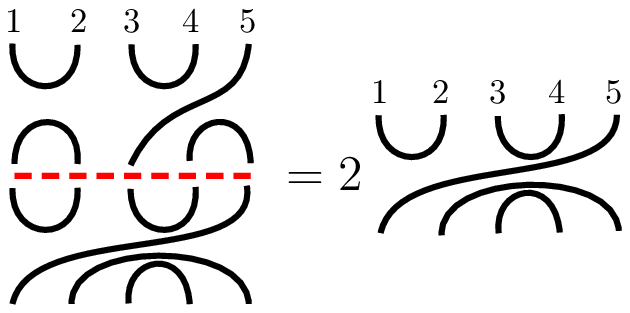}}   \, \, =\, \, 2D_3
\end{equation}
The isomorphism between $ \TLn$ and $ \TLndiag $ is given by 
\begin{equation}\label{dib4} 
  \one \, \, \, \,  \mapsto\, \, \,  \, \raisebox{-.35\height}{\includegraphics[scale=0.75]{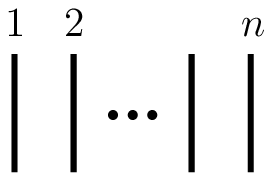}}, \, \, \, 
    \UU_i  \, \, \mapsto\, \,  \raisebox{-.35\height}{\includegraphics[scale=0.75]{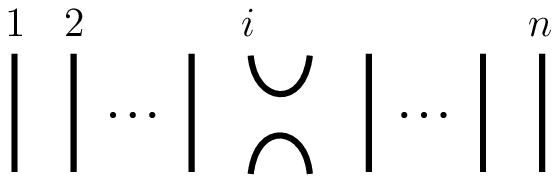}} 
\end{equation}
where $ \one $ is the {\color{black}{unit}}-element of $ \TLn$.  
From now on we shall identify $ \TLn$ with $ \TLndiag$ 
via this isomorphism. There is {\color{black}{a}} similar
isomorphism for the specialized Temperley-Lieb algebra $ \TLnk $
and here we shall also identify $ \TLnk $ with the corresponding diagrammatic algebra, defined over $ \Bbbk $. 

\medskip
Throughout the paper we shall be interested in the {\it Jones-Wenzl idempotent}
$\JWn $ of $ \TLnQ $,
{\color{black}{see \cite{Jo} and \cite{Wenzl}. }}
It is the unique nonzero idempotent of 
$ \TLnQ $ satisfying
\begin{equation}\label{definingpropertyJW}
 \UU_i  \JWn = \JWn \UU_i  = 0  \mbox{ for all } i 
\end{equation}
We use the following standard diagrammatic notation
for $\JWn $ 
\begin{equation}
\JWn = 
  \raisebox{-.35\height}{\includegraphics[scale=0.7]{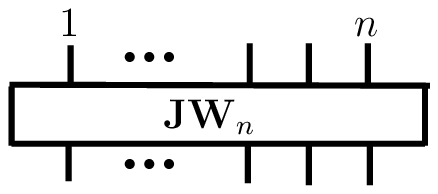}} \in \TLnQ
\end{equation}
For example we have that 
\begin{equation}\label{dib7}
   \raisebox{-.35\height}{\includegraphics[scale=0.7]{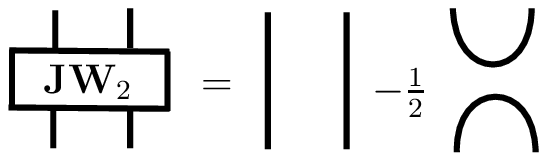}} 
 \end{equation}
\begin{equation}\label{dib8}
   \raisebox{-.35\height}{\includegraphics[scale=0.7]{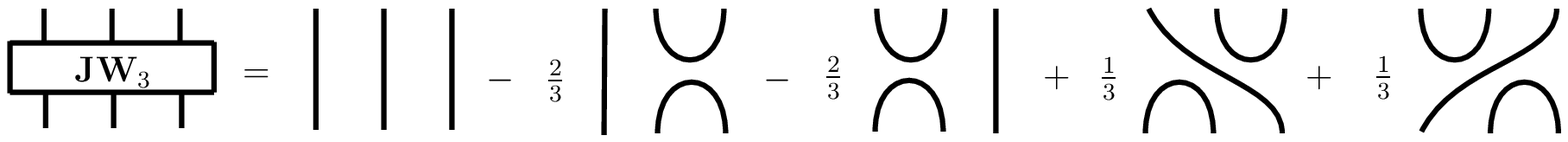}} 
 \end{equation}
In general, as one already observes in \eqref{dib7} and \eqref{dib8}, 
when expanding $ \JWn  $ in terms of the diagram basis for
$  \TLnQ $, 
the coefficient of $ \one $ is 1, whereas the other coefficients in general
are rational numbers $ \dfrac{a}{b} $ with 
non-trivial denominator. These denominators $ b $ prohibit the
specialization of $ \JWn  $ to fields $ \Bbbk $ whose characteristic $ p $ divides $ b $. 

\medskip
On the other hand, 
we always have $ \mathbf{JW}_n^{\ast} = \JWn $ where $ \ast $ is the antiautomorphism of $ \TLnQ $
given by reflection
through a horizontal axis, and similarly $ \JWn $ is symmetric with respect to reflection through
a vertical axis.
These properties can be observed in \eqref{dib7} and \eqref{dib8}. 

\medskip
For general $ n $ there is no known closed formula for calculating 
the coefficients of $ \JWn  $ in terms of the diagram basis for $ \TLnQ$; all known formulas are recursive.
We shall need the following recursive formula that goes back to
Jones and Wenzl, see \cite{Jo} and \cite{Wenzl}.

\begin{equation}\label{goesbackto}
  \raisebox{-.35\height}{\includegraphics[scale=0.7]{dib9.eps}}
\end{equation}
Combining it with \eqref{eq oneTL} 
we obtain the following well{\color{black}{-}}known formula 
\begin{equation}\label{closing}
  \raisebox{-.35\height}{\includegraphics[scale=0.7]{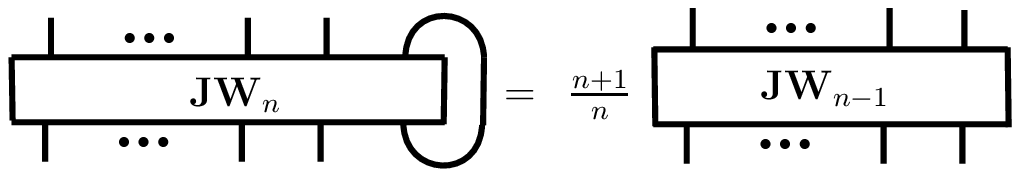}}
\end{equation}
which can be repeated to arrive at 
\begin{equation}\label{basicelements}
  \raisebox{-.35\height}{\includegraphics[scale=0.7]{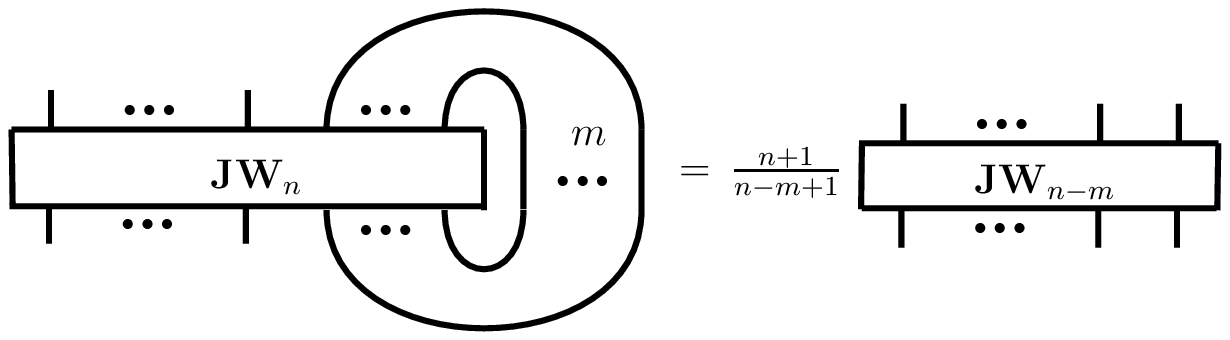}}
\end{equation}

Using \eqref{goesbackto} one proves that
for $ m<n $ we have 
$ \mathbf{JW}_{m} \mathbf{JW}_{n} =
\mathbf{JW}_{n} $, or diagrammatically 
\begin{equation}\label{absorbtion}
  \raisebox{-.35\height}{\includegraphics[scale=0.7]{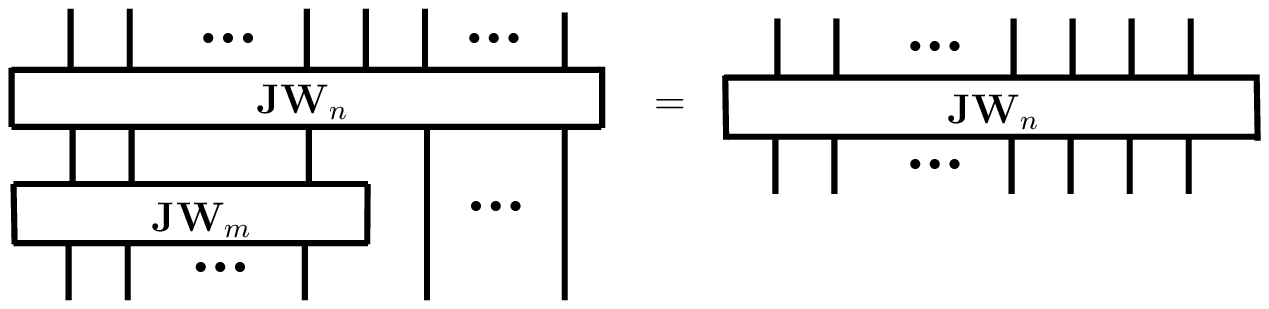}}
\end{equation}

We next recall the basic elements of the representation theory of $ \TLn$, using the language
of {\it cellular algebras}.
The notion of cellular algebras was introduced by Graham and Lehrer in
\cite{GL} and in fact $ \TLn $ was one of their motivating examples. 

\medskip
{\color{black}{Throughout the paper, when referring to \lq representations\rq\ and \lq actions\rq\ we
shall in general mean \lq right representations\rq\ and \lq right actions\rq\ \!\!.}}

\begin{definition}\label{cellular}
  Let $ A $ be an associative $ R$-algebra with unit, where $ R $ is a commutative ring. A cell datum for
  $ A $ is a triple $ (\Lambda, T, C) $ where $ \Lambda = (\Lambda, >) $ is a finite poset, 
  $ T$ is a function from $ \Lambda $ to finite sets and $ C $ is an injective function 
\begin{equation}  
  C: \prod_{\lambda \in \Lambda} T(\lambda) \times T(\lambda) \rightarrow A,  (s,t) \in T(\lambda ) \times
T(\lambda )  \mapsto C_{s t}^{\lambda} 
\end{equation}    
These data should satisfy
{\color{black}{the following conditions.}}
\begin{enumerate}
\item  The image of $ C $, that is $im \, C= \{ C_{st}^{\lambda}  \, | \, (s,t) \in
 T(\lambda) \times T(\lambda), \lambda \in \Lambda \} $, is an $ R $-basis for $ A $.
\item For any $ a \in A $ we have
\begin{equation}\label{mult}  C_{st}^{\lambda}a 
  = \sum_{v \in T(\lambda)} r_{tva}C_{sv}^{\lambda}  \, \mbox{ mod }  A^{\lambda}
\end{equation}  
where $ A^{\lambda} $ is the $ R$-submodule of $ A $ spanned by $ \{ C_{st}^{\mu}  \, | \, \mu > \lambda \}$.
\item The $R$-linear map $ \ast$ of $ A$ given by $ C_{st}^{\lambda} \mapsto C_{ts}^{\lambda} $ is an algebra 
anti-isomorphism of $A$. 
\end{enumerate}
\end{definition}
An algebra endowed with a cell datum $ (\Lambda, T, C) $ is called a cellular algebra, with
cellular basis $ im \, C$.

\medskip
$\TLn $ is an example of a cellular algebra. To see this one lets $ \Lambda $ be the set of two-column
integer partitions $ \ParTwo $,  
endowed with the usual dominance order.
{\color{black}{Thus, for $ \lambda=(2^{l_2}, 1^{l_1-l_2} ), \mu=(2^{m_2}, 1^{m_1-m_2} ) \in  \ParTwo$
one has $ \lambda \unlhd \mu $ if and only if $ l_2 \le m_2 $.}}
For $ \lambda \in \ParTwo $ one lets $ T(\lambda) $ be the set of standard $\lambda$-tableaux
$ \std(\lambda) $. 
To explain $ C $, one first constructs for $ \lambda \in \ParTwo $ and $ \T \in \std(\lambda) $ 
a {\it Temperley-Lieb half-diagram} $ C_{\T}^{\lambda}$ for $ \TLn$ as follows. Going through the numbers
$ \{1, 2\, \ldots, n \} $ in increasing order, one raises for any $ i $ occurring in the first column of $ \T $
a vertical line from the $ i $'th lower position of the rectangle 
and for any $ i $ occurring in the second column of $ \T $, one joins the $i$'th lower position with the 
top end of the first vacant line to the left, always avoiding line crossings. Here is an example for
$ \lambda = (2^4,1^3) \in {\rm Par}_{11} $. 
\begin{equation}
  \T :=  \raisebox{-.5\height}{\includegraphics[scale=0.7]{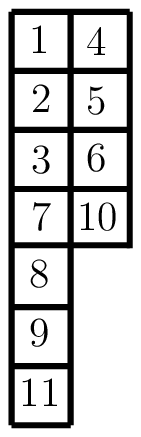}}
  \mapsto  \, \, \, \,   C_{\T}^{\lambda} = \,    \raisebox{-.3\height}{\includegraphics[scale=0.7]{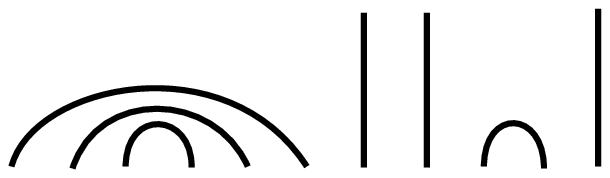}}
\end{equation}
For a pair of standard $\lambda$-tableaux $ (\s, \T) $, one then defines $ C_{\s \T}^{\lambda} $ as the
diagram obtained from $  C_{\s}^{\lambda} $ and $ C_{\T}^{\lambda}  $ by reflecting
$ C_{\s}^{\lambda} $ horizontally and concatenating below with $C_{\T}^{\lambda} $. 
Here is an example. 
\begin{equation}\label{diagrambasisTL}
 (\s,  \T) :=  \raisebox{-.5\height}{\includegraphics[scale=0.7]{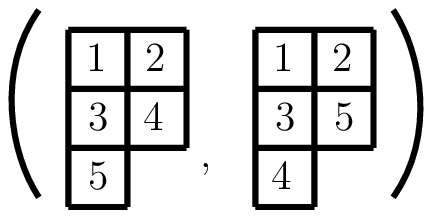}}
  \mapsto  \, \, \, \,   C_{\s\T}^{\lambda} = \,    \raisebox{-.5\height}{\includegraphics[scale=0.7]{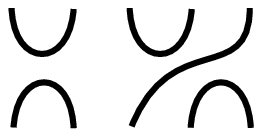}}
\end{equation}
Using the multiplication rules explained in \eqref{dib2} and \eqref{dib3}, 
one now checks that $ \TLn $ indeed is a cellular algebra over $ \Z$, with the ingredients just introduced,
and similarly $ \TLnk $ is a cellular algebra over $ \Bbbk $.

\medskip
In general, for a cellular algebra $ A $ there is a family of {\it cell modules}
$ \{ \Delta(\lambda) \, | \, \lambda \in \Lambda \} $ that play a key role
when studying the representation theory of $ A $. To define $ \Delta(\lambda) $ one chooses an
arbitrary $ s_0 \in T(\lambda) $ and sets
\begin{equation}
\Delta(\lambda) := {\rm span}_R \{ C_{ s_0 t}^{\lambda} \, | \, t \in T(\lambda) \}
\end{equation}
The action of $ A $ on $ \Delta(\lambda) $ is given by
$ C_{ s_0 t}^{\lambda} a  =  \sum_{v \in T(\lambda)} r_{tva}C_{s_0v}^{\lambda} $ where 
$ r_{tva} \in R$ is the scalar that appears in \eqref{mult}.
We shall sometimes write $  \Delta^R(\lambda) $ for $  \Delta(\lambda) $ to indicate the dependence on
the ground ring $ R $. 

\medskip
In the case of $ \TLn$, we identify for $ \s_0, \T \in \std(\lambda) $ the $ \Delta(\lambda) $ basis element 
$ C_{\s_0 \T}^{\lambda} $ with the half-diagram $ C_{ \T}^{\lambda} $. Under this identification, 
for a Temperley-Lieb diagram $ D $ we have that $ C_{ \T}^{\lambda} D $ 
is the concatenation with $ C_{ \T}^{\lambda}$ on top of $ D $, where internal loops are removed
by multiplying by $ 2$, and where half-diagrams that do not belong to $ \{ C_{ \T}^{\lambda} \,| \, \T
\in \std(\lambda) \} $ are set equal to zero.

\medskip
In the present paper, Jucys-Murphy elements play an important role. Their key properties
were developed in Murphy's papers in the eighties,
see \cite{Murphy2}, \cite{Murphy}, \cite{Murphy1}. 
These properties were formalized by Mathas as follows, see \cite{Mat-So}. 
\begin{definition}\label{JM}
  Let $ A $ be a cellular algebra over $ R $ with triple $ (\Lambda, T, C) $, and suppose that for 
  each $ \lambda \in \Lambda $ there is a poset structure on $ T(\lambda) $ with order relation $ <$.
  Then a family of elements $ \{L_1, L_2, \ldots, L_M\} $
  is called a set of $\JM$-elements for $ A $ if it satisfies the following conditions. 
\begin{enumerate}
\item The $L_i  $'s commute and satisfy $ L_i^{\ast} = L_i$.
\item For each $ t \in T(\Lambda) := \coprod_{\lambda \in \Lambda} T(\lambda) $ there
  is a function $ c_{t}: \{1, 2, \ldots, M\} \rightarrow R $ such that the following
  triangularity formula holds in $ \Delta(\lambda) $
\begin{equation}  
C^{\lambda}_t L_i = c_t(i) C^{\lambda}_t + \sum_{s\in T(\lambda), s>t } a_s C_s^{\lambda}
\end{equation}  
for some $ a_s \in R$. 
\end{enumerate}
\end{definition}
If $ \{ L_1, L_2, \ldots, L_M \} $ is a family of $\JM$-elements for $ A $, then the
corresponding functions $ c_{t}$ are called {\it content functions}. 

\medskip
$\JM$-elements were first constructed for the group algebra of the symmetric group,
and
from these one obtains $\JM$-elements for $ \TLn$,
as we shall shortly see.

\medskip
Let $ \Si_n $ be the symmetric group of bijections
{\color{black}{on}}
$ \{1,2, \ldots, n \} $. For $ \sigma \in \Si_n $ and
$ i \in \{1,2, \ldots, n \} $ we write $ (i)\sigma $ for the image of $i $ under $ \sigma$.
We shall use standard cycle notation
for elements of $ \Si_n $, that is $ \sigma = (i_1, i_2, \ldots, i_k ) $ is the element of $ \Si_n $
defined by $ (i_1)\sigma = i_2, (i_2) \sigma = i_3, \ldots, (i_k)\sigma = i_1 $.
For elements $ \sigma_1, \sigma_2 \in \Si_n $ the product $ \sigma_1 \sigma_2 \in \Si_n $ is given by
$ (i) (\sigma_1 \sigma_2) = ((i) \sigma_1 ) \sigma_2 $.

\medskip
Let 
$ \{ L_1, L_2, \ldots, L_n \}  \subseteq  \Z \Si_n \,  {(\mbox{or }} \Bbbk \Si_n) $ be defined by 
\begin{equation}\label{defJM}
L_1 := 0, \mbox{ and } L_i := (1,i) + (2,i) + \ldots + (i-1, i ) \mbox{ for } i = 2,3, \ldots, n 
\end{equation}  
Define moreover for $\T \in \std(\lambda) $ the function 
$ c_{\T}:  \{1,2,\ldots, n\} \rightarrow \Z \, {(\mbox{or }} \Bbbk) $ by
\begin{equation}\label{contentdef}
 c_{\T}(i) := c-r  \mbox{ for } \T[r,c] = i
\end{equation}
where $ \T[r,c] $ is the number that appears in the $ r$'th row of $ \T $, counted from top to bottom,
and in the $ c$'th column of $ \T$, counted from the left to right.

\medskip
$ \Si_n $ is a Coxeter group on the {\it simple transpositions} $ s_i := (i,i+1) $. 
We need the following well{\color{black}{-}}known fundamental Lemma, which is easily verified. 
\begin{lemma}\label{wellknownfunda}
There is a 
surjection $ \Phi: \Z \Si_n \rightarrow \TLn $, given by $ s_i \mapsto \UU_i -\one$.
The kernel of $ \Phi $ is the ideal in $  \Z \Si_n  $ generated by
$ s_1 s_2 s_1 +s_1 s_2 + s_2 s_1  +s_1 + s_2 + 1 $. 
A similar statement holds over $ \Bbbk $. 
 \end{lemma}  

Let $ \LL_i := \Phi(L_i) $. We represent  $   \LL_i $ diagrammatically as follows
\begin{equation}
\LL_i =  \raisebox{-.35\height}{\includegraphics[scale=0.7]{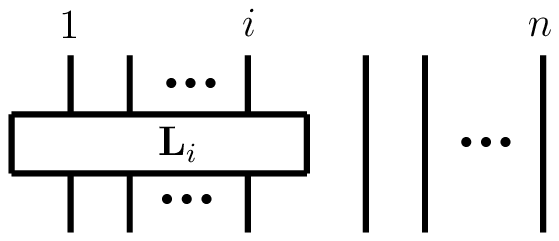}}
\end{equation}

We now have the following key result. 
\begin{theorem}\label{followingkeyresult}
$ \{ \LL_1,    \LL_2,   \ldots,  \LL_n \} $ 
is a family of $\JM$-elements for $ \TLn$ with respect to the content functions
$ c_{\T}$, defined in \eqref{contentdef}. 
\end{theorem}
\begin{dem}
  $ \{ L_1, L_2, \ldots, L_n \} $ is known to be a family of $\JM$-elements
  for the cellular structure on $ \Z \Si_n $ given by the specialization $ q = 1 $ of
  Murphy's standard basis for the Hecke algebra, see \cite{Mat} and \cite{Murphy1}. 
  On the other hand, $ \Phi: \Z \Si_n  \rightarrow \TLn $ maps the standard basis cellular structure
  on $ \Z \Si_n $ to the diagram basis cellular structure on $ \TLn $ and therefore
  $ \{ \LL_1,    \LL_2,   \ldots,  \LL_n \} $ is a family of $\JM$-elements
  for $ \TLn$, as claimed. For more details one should consult \cite{Orm}. 
\end{dem}  

{\color{black}{
 \begin{remark}\label{JMremark}
  \normalfont
Jucys-Murphy elements for the Temperley-Lieb algebra have been considered before 
in \cite{Enyang} and in \cite{HaMazzRam}. The Jucys-Murphy elements in
\cite{Enyang} are different from ours. The Jucys-Murphy elements in \cite{HaMazzRam} are also different from ours
since they are \lq multiplicative\rq\ and do not specialize to the $ q=1 $ setting of the
present paper.
 \end{remark}
}}
\section{The separated case}\label{separated case}
In this section we consider the rational Temperley-Lieb algebra $\TLnQ $.
The ground ring for $ \TLnQ $ is $ \QQ $ which implies that for
two-column partitions $ \lambda $ and $ \mu $ and for 
standard tableaux
$ \s \in \std(\lambda) $ and $ \T \in \std(\mu) $ we have that 
\begin{equation}\label{separationcondition}
c_{\s}(i) = c_{\T}(i) \mbox{ for } i= 1, 2, \ldots, n \Longrightarrow \s= \T
\end{equation}
In other words, the {\it separation condition} in \cite{Mat-So} is fulfilled and so
$ \TLnQ$ is semisimple. The separation condition also implies that
the simultaneous action of the $ \LL_i$'s
on $  \TLnQ $ via right multiplication is diagonalizable
with eigenvalues given by the $ c_{\T}(i)$'s, and similarly for the left action.
Moreover, under the separation condition we have the
following expression for the idempotent projector $ \EE_{\T } $ for the common eigenvector for
all the $ \LL_i$'s with eigenvalues $ c_{\T}(i)$
\begin{equation}\label{IdempotentHecke1}
  \EE_{\T} = \prod_{c \in {\cal C}} \, \, \prod_{\substack{i=1, \ldots, n\\ c \neq c_{\T}(i)} } \dfrac{\LL_i-c}{c_{\T}(i)-c}
  \in \TLnQ
\end{equation}
where $ {\cal C } $ is the set of contents for standard tableaux of two-column partitions
{\color{black}{of $n$}}, that is 
\begin{equation}\label{IdempotentHecke2}
{\cal C }:= \{ c_{\T}(i) \, | \,  i=1,2, \ldots, n \mbox{ and } \T \in \std(\ParTwo) \} \, \, \mbox{where} \, \, \, 
\std(\ParTwo) := \bigcup_{\lambda \in \ParTwo} \std(\lambda)
\end{equation}
With $ \T $ running over $\std(\ParTwo)$, the $ \EE_{\T}$'s form a 
complete
set of orthogonal
primitive
idempotents for $ \TLnQ $, that is we have 
\begin{equation}\label{IdempotentHecke3}
  \one = \sum_{\T \in \std(\ParTwo)} \EE_{\T}, \, \, \, \, \,  \LL_i \EE_{\T} = \EE_{\T} \LL_i=  c_{\T}(i) \EE_{\T} , \,\, \, \, \,  \EE_{\s} \EE_{\T} =
  \delta_{\s \T} \EE_{\s}
\end{equation}
where $  \delta_{\s \T}  $ is the Kronecker delta. 
  
\medskip
The formulas in \eqref{IdempotentHecke1} and \eqref{IdempotentHecke3} are
consequences of the general theory for $ \JM$-elements developed in \cite{Mat-So}. 
For $ \QQ \Si_n $ the analogues of \eqref{IdempotentHecke1} and 
\eqref{IdempotentHecke3} were first found by Murphy in \cite{Murphy3}. We find it worthwhile to mention
that the corresponding 
properties do not hold for the Young symmetrizer idempotents for $ \QQ \Si_n $, since these are
not orthogonal.

\medskip
The expression for $ \EE_{\T} $ given in \eqref{IdempotentHecke1} contains many redundant factors and
is in general intractable, in the symmetric group case as well as in the Temperley-Lieb case.

\medskip
The purpose of this section is to give a new expression for $ \EE_{\T} $ in the Temperley-Lieb case,
using Jones-Wenzl idempotents. In view of this, one may now consider seminormal forms and Jones-Wenzl
idempotents as two aspects of the same theory.

\medskip
Let $ \T \in \std(\lambda) $ be a two-column standard tableau. Then $ \T $ can be written in the form 
\begin{equation}\label{intheform}
\T = \raisebox{-.5\height}{\includegraphics[scale=0.75]{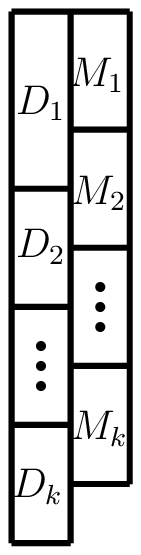}}
\end{equation}
where each $ D_i $ and $ M_i $ is a non-empty block
of consecutive natural
numbers, except that $ M_k $ is allowed to be empty, satisfying
that
{\color{black} the numbers of $ D_i $ are less than the numbers of $ M_i$ and that
the numbers of $ M_i $ are less than the numbers of
$ D_{i+1}$ for all $i$.
Let $  d_i := | D_i| $  and $ m_i  = | M_i| $ be the cardinalities of $ D_i $ and $ M_i $, respectively. 
Then $ d_1+d_2 + \ldots +d_i \ge m_1+m_2 + \ldots +m_i $ for all $ i $. Moreover, 
each sequence of blocks $ D_1, M_1, D_2, \ldots, M_k $ 
satisfying all these conditions gives rise to a two-column standard tableau}
{\color{black}and in this way we obtain a bijective correspondence between
such sequences of blocks and two-column standard tableaux.
}

\medskip
For $ i=1, 2, \ldots, k $ define $ n_i $ via 
\begin{equation}\label{wenowassociate}
n_1 := d_1 \mbox{ and }  n_i = (d_1+d_2+ \ldots + d_i) - ( m_1+m_2+ \ldots + m_{i-1}) \mbox{ for } i >1 
\end{equation}
We now associate with $ \T $ an element $ f_{\T} \in \Delta^{\QQ}(\lambda) $ in the following recursive way. 
Suppose first that $  M_k \neq \emptyset  $. If 
$ k= 1 $ we set
\begin{equation}\label{exA}
{\color{black}{f_\T = \raisebox{-.35\height}{\includegraphics[scale=0.7]{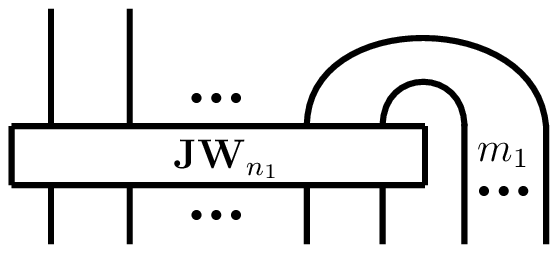}}}}
\end{equation}
and if $ k= 2 $ we set
\begin{equation}\label{exB}
{\color{black}{f_\T \, \, \, = \, \, \,  \raisebox{-.5\height}{\includegraphics[scale=0.7]{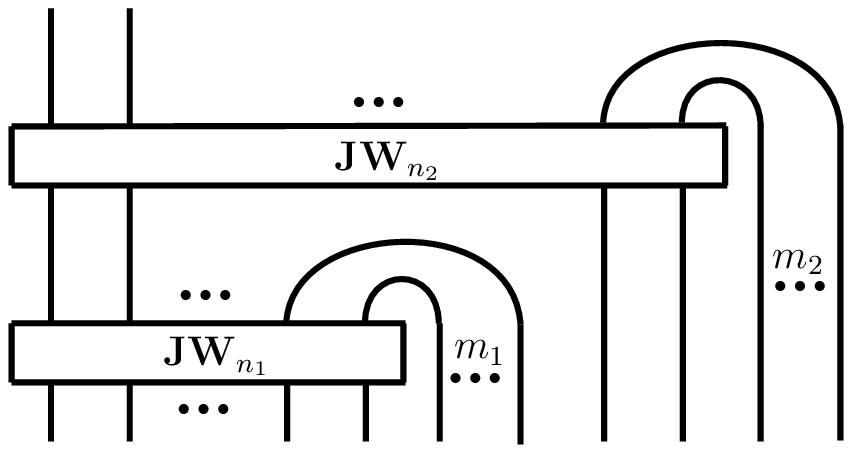}}}}
\end{equation}
We repeat this recursively, that is in the $i$'th step we first concatenate on top with $ \mathbf{JW}_{\! n_i } $
and then bend down the $ m_i$ top and rightmost lines 
to the bottom. If $ M_k = \emptyset $ the construction is the same as for
$ M_k \neq \emptyset  $, except that in the last step the bending down of the $ m_k$ top and rightmost 
lines is omitted. 

\medskip
For example,
\begin{align}\label{exC}
& \text{if }\T \,  := \,   \raisebox{-.5\height}{\includegraphics[scale=0.7]{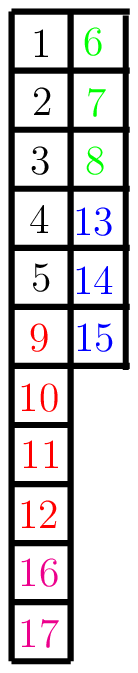}} & 
& \text{we have that} & 
& f_\T \, \, \, = \, \, \,  \raisebox{-.5\height}{\includegraphics[scale=0.6]{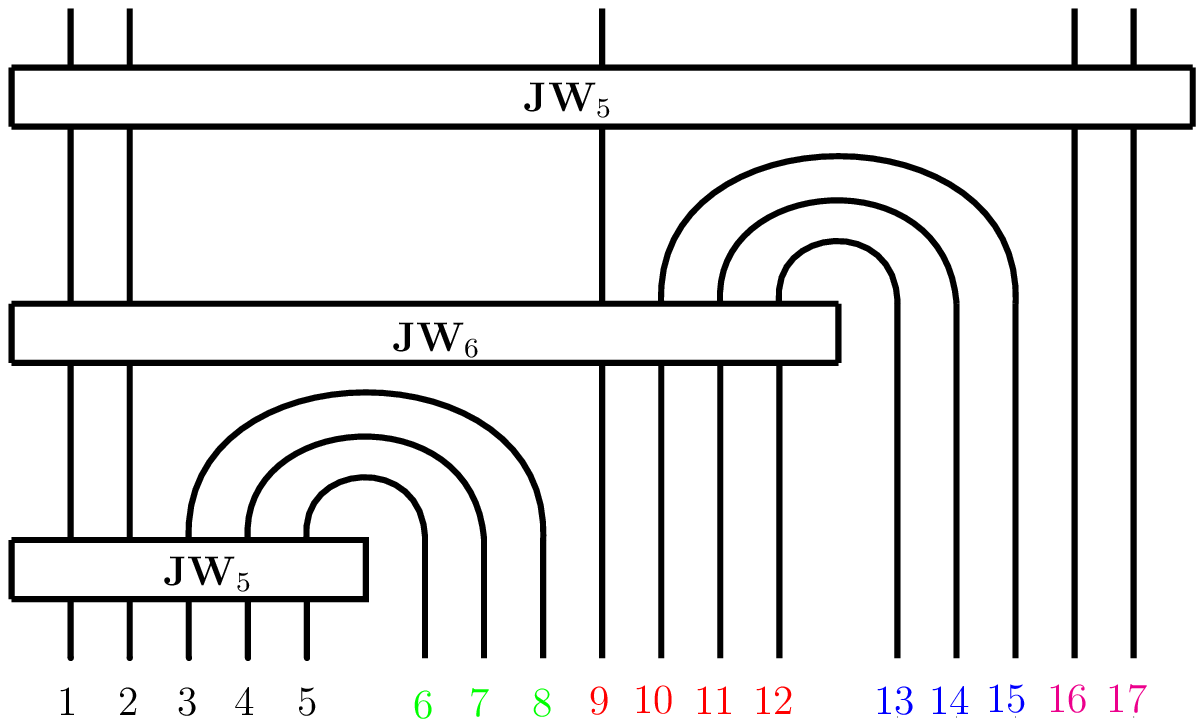}}
\end{align}
In general, if $ i $ appears in the first column of $ \T $ then $ i $ is connected
in $ f_\T $ to the southern border of a Jones-Wenzl element, and 
if $ i $ appears in the second column of $ \T $ then $ i $ is connected
in $ f_\T $ to the northern border of a Jones-Wenzl element. We have 
indicated this in \eqref{exC}, using colors. 

\medskip
In general, for $ \T $ as in \eqref{intheform}
we shall sometimes represent $ f_\T $ in the following schematic way
\begin{equation}\label{shorthand}
f_\T = \raisebox{-.5\height}{\includegraphics[scale=0.7]{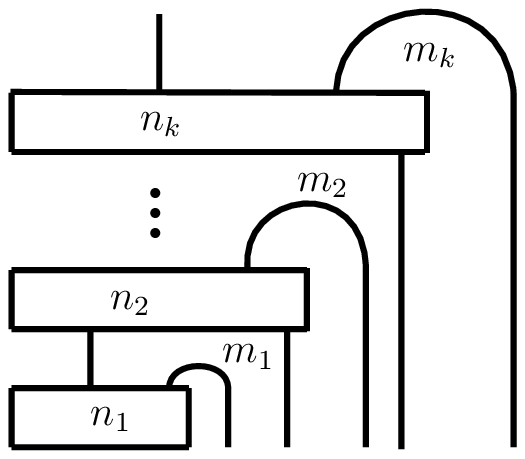}}
\end{equation}
where $ n $ is a shorthand for $ \JWn $ and 
where $m_i $ indicates the number of lines being bent down, which may be zero for $ m_k $.

\medskip
For $ \T $ a two-column standard tableau we set 
\begin{equation}
\gamma_{\T} := \prod_{j=1}^{k }\dfrac{n_i +1}{n_i-m_i +1} 
\end{equation}
We define 
$  f_{\T \T} $ as the concatenation of $ f_{\T}^{\ast} $ with $ f_{\T} $ with
$ f_{\T}^{\ast} $ on top of $ f_{\T} $ and 
finally we define $ \EE_{\T}^{\prime} \in \TLnQ$ as 
$ \EE_{\T}^{\prime} := \frac{1}{\gamma_{\T}} f_{\T \T} $, 
or diagrammatically
\begin{equation}\label{mainconstruction} 
  \EE_{\T}^{\prime} := \dfrac{1}{\gamma_{\T}}
\raisebox{-.5\height}{\includegraphics[scale=0.7]{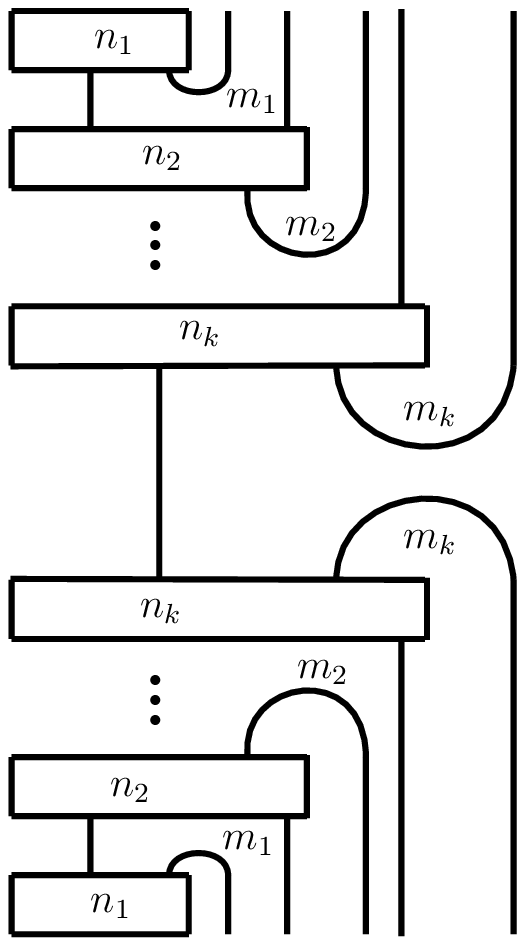}}
\end{equation}

The elements $ \EE_{\T}^{\prime} $ have already appeared in the literature,
see \cite{CoHo}, \cite{GW} and \cite{PMartin}, 
{\color{black}{with our diagrammatic approach essentially being the one of \cite{CoHo}.}}
{\color{black}{The purpose of this section is to show $ \EE_{\T}  = \EE_{\T}^{\prime} $.
This identity is mentioned in \cite{ElH}}}.

\medskip
{\color{black}{The following Theorem has already appeared in \cite{CoHo}, see also
\cite{GW} and \cite{PMartin}, but we still include it for completeness.}}
\begin{theorem}\label{startby}
$ \{ \EE_{\T}^{\prime} \, | \, \T \in \std(\ParTwo) \} $ is a
set of
orthogonal idempotents
  for $ \TLnQ$. 
\end{theorem}
\begin{dem}
  We first observe that \eqref{basicelements} implies 
$ f_{\T \T}^2 = \gamma_\T f_{\T \T} $ and so $ \EE_{\T}^{\prime} $ is indeed an idempotent.
  Similarly, we observe that $ f_{\T}^{ \ast} C_{\T }^{} = \gamma_{\T} { \rm JW}_{n_k -m_k } $
  from which it follows that $ f_{\T}  \neq 0 $, and hence also
  $ \EE_{\T}^{\prime}  \neq 0 $.
  
\medskip
We next assume that $ \T \neq \overline{\T} $ and must show 
that $ \EE_{\T}^{\prime} \EE_{\overline{\T}}{}^{\! \! \prime} = 0 $ which can be done by showing
that $ f_{\T}^{} f_{\overline{\T}}{}^{\! \! \ast} = 0 $.  
Letting $ \{\overline{D}_i \, | \, i=1,2, \ldots, \overline{k} \} $ and 
$ \{\overline{M}_i \, | \, i=1,2, \ldots, \overline{k} \} $
be the blocks for $ \overline{\T} $, as in 
\eqref{intheform},
and defining $ \overline{n}_i $ and $ \overline{m}_i $ as in
\eqref{wenowassociate}, 
we must show that the following diagram is zero
\begin{equation}\label{mustshow}
f_{\T}^{} f_{\overline{\T}}{}^{\! \! \ast} = \raisebox{-.5\height}{\includegraphics[scale=0.65]{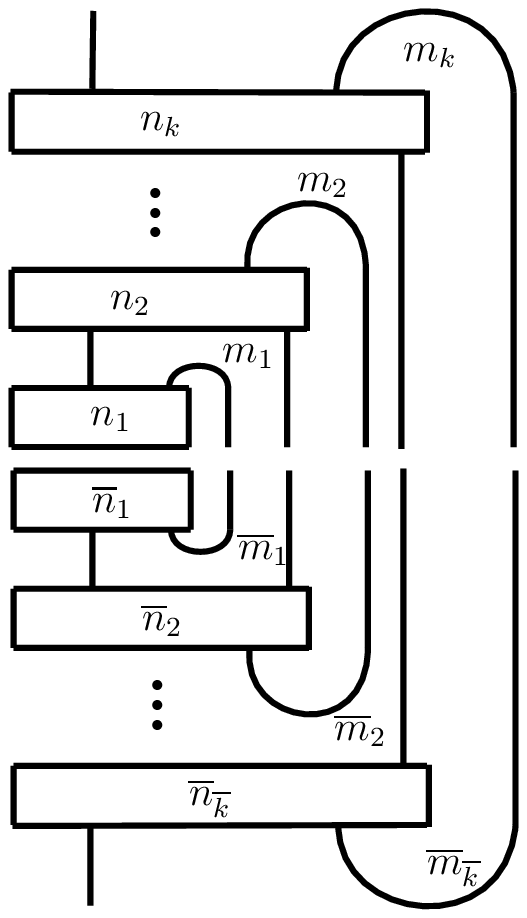}}
\end{equation}
If $ n_1 = \overline{n}_1 $ and $ m_1 = \overline{m}_1 $ then \eqref{mustshow} is equal
to $ \dfrac{n_1 +1}{n_1-m_1 +1}  $ times 
$ f_{\T_1}^{} f_{\overline{\T}_1}{}^{\! \! \ast}  $ 
where $ \T_1 $ and $ \overline{\T}_1 $ are the standard tableaux obtained from
$ \T $ and $ \overline{\T} $  
by removing the blocks $ M_1, D_1 $ and $ \overline{M_1}, \overline{D_1} $
and so 
we may assume that $ n_1 \neq \overline{n}_1  $ or $ m_1 \neq \overline{m}_1 $. 
If $ n_1 < \overline{n}_1 $ then at least one line from $ \mathbf{JW}_{n_1} $ is bent down to
$ \mathbf{JW}_{n_1^{\prime}} $, and so it follows from
\eqref{absorbtion} that 
the resulting diagram is zero: to illustrate this we take 
$ n_1 = 3 $ and $ n_1^{\prime} = 4 $ where the relevant part of
$ f_{\T_1}^{} f_{\overline{\T}_1}{}^{\! \! \! \! \ast}  $ is 
\begin{equation}\label{mustshow1}
\raisebox{-.55\height}{\includegraphics[scale=0.65]{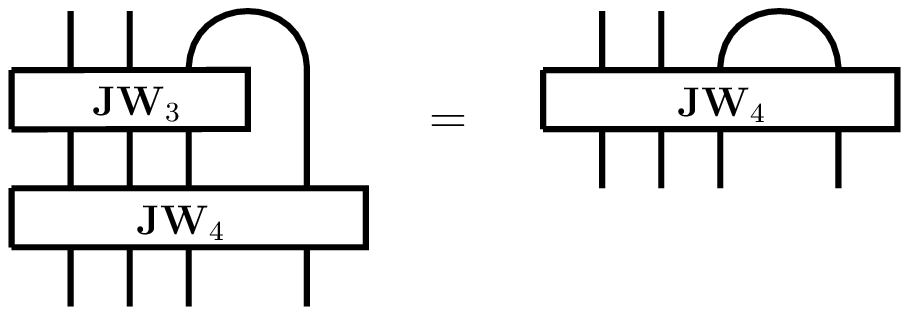}} \, \, = \, \, 0
\end{equation}  
If $ n_1 > \overline{n}_1 $ one applies $ \ast $ to \eqref{mustshow} and is then reduced to
the previous case $ n_1 < \overline{n}_1$.
If $ n_1 = \overline{n}_1 $ and $ m_1>  \overline{m}_1 $ then
a line from $ \mathbf{JW}_{n_1} $ is bent down to $ \mathbf{JW}_{n_2^{\prime}} $ 
and so the resulting diagram is also zero in this case. Let us illustrate this using $ n_1 = \overline{n}_1 = 4 $ and
$ m_1 =3, \overline{m}_1 =2 $ and $ n_2  = 3$ where the relevant part of 
\eqref{mustshow} is as follows
\begin{equation}
\raisebox{-.5\height}{\includegraphics[scale=0.7]{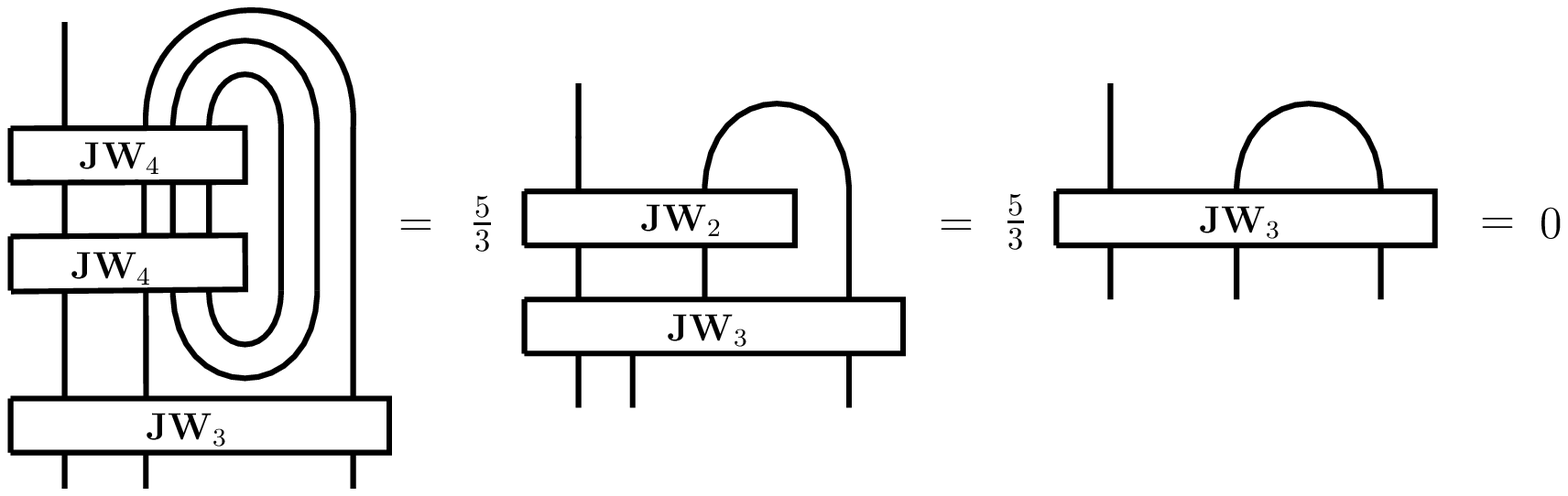}} 
\end{equation}  
Finally, if $ n_1 = \overline{n}_1 $ and $ m_1<  \overline{m}_1 $ we once again first apply $ \ast$ and
are then reduced to the previous case. This proves that $\{ \EE_{\T}^{\prime}\, | \, \T \in \std(\ParTwo) \} $
is a set of orthogonal idempotents.

\end{dem}

\medskip
\begin{corollary}
  Let $ \lambda $ be a two-column partition. Then 
  $ \{ f_\T \, | \, \T \in \std(\lambda) \} $ is a $\QQ$-basis for $ \Delta^{\QQ}(\lambda)$. 
\end{corollary}
\begin{dem}
  We have that $  f_{\T} \EE^{\prime}_{\s} = \delta_{ \s \T} f_{\T} $ and so it
  follows from Theorem \ref{startby} that 
    $ \{ f_\T \, | \, \T \in \std(\lambda) \} $ is a $ \QQ$-linearly independent subset of $ \Delta^{\QQ}(\lambda)$. 
Since $ \dim \Delta^{\QQ}(\lambda) = | \std(\lambda) | $ it is also a basis for $ \Delta^{\QQ}(\lambda)$. 
\end{dem}

\medskip
The action of $ \Si_n $ on $ \{ 1,2, \ldots, n \} $ extends naturally to an action
of $ \Si_n $ on $ \lambda$-tableaux, by permuting the entries. For $ \T $ a $ \lambda$-tableau
and $ \sigma \in \Si_n $, we denote by $ \T \sigma $ the action of $ \sigma $ on $ \T$.

\medskip
We now set out to prove 
$ \EE_{\T}^{\prime} = \EE_{\T} $. Our proof will be an induction over the dominance order on standard
tableaux and for this the following Theorem is a key ingredient.

\begin{theorem}\phantomsection\label{YSFfirst}
  Suppose that $ \T \in \std(\lambda) $ where $ \lambda \in \ParTwo$. Suppose first that
  for a simple transposition $ s_i \in \Si_n$ we have that $ \T s_i \in \std(\lambda) $ and that 
  $ \T  \unlhd \T s_i $. Then, setting $ \T_d := \T, \T_u := \T s_i $
and $ r := c_{\T_u}(i) - c_{\T_d}(i)  $, the following formulas hold  
  \begin{description}
  \item[a)]
$  f_{\T_d} \UU_i=
\dfrac{ r +1 }{ r }
f_{\T_d} + \dfrac{ r^2-1}{ r^2} f_{\T_u}$
  \item[b)]
$  f_{\T_u} \UU_i=
\dfrac{ r -1 }{ r }
f_{\T_u} +  f_{\T_d}$
    \end{description}
Suppose next that $ \T s_i \notin \std(\lambda) $. Then 
  \begin{description}
  \item[c)]
    $  f_{\T} \UU_i= 0  \,\, \, \,\, \, \,   \mbox{ if } i, i+1 \mbox{ are in the same column of } \T $ 
 \item[d)]
    $  f_{\T} \UU_i= 2 f_{\T}  \, \mbox{ if } i, i+1 \mbox{ are in the same row of } \T $ 
    \end{description}
\end{theorem}
\begin{dem}
  We first show {\bf a)}.
We have blocks $ D_j $ and $ M_j $ for $ \T $, as in \eqref{intheform}.
  By the assumptions, $ i $ lies in the first column of $ \T$,
  as the biggest number of a block $ D_j $, whereas $ i+1 $ lies in the second column of $ \T$, as 
the smallest number of the block $ M_j $, as indicated in the example below. 
\begin{equation}\label{2.17}
\T= \raisebox{-.5\height}{\includegraphics[scale=0.7]{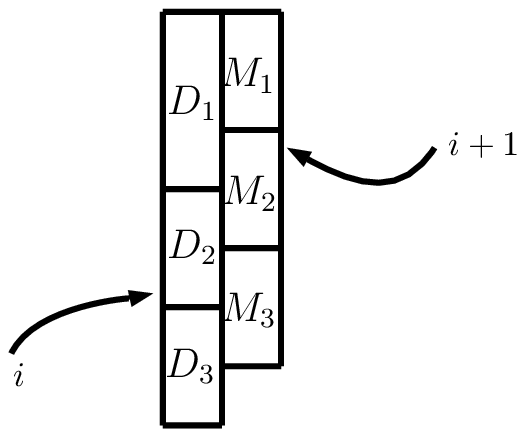}}, \, \, \,  \, \, \,  \, \, \, 
f_\T= \raisebox{-.5\height}{\includegraphics[scale=0.7]{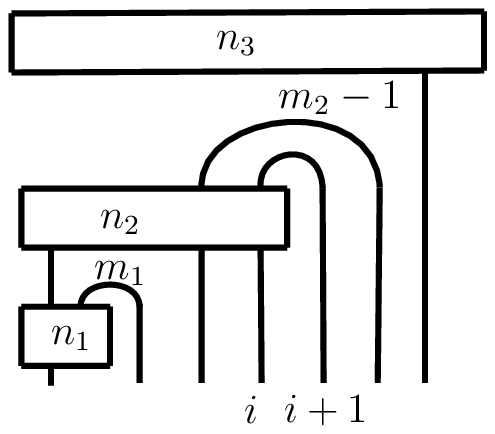}} 
\end{equation}
In \eqref{2.17}, we have indicated the corresponding $ f_{\T} $ and have singled out the lines for $ f_{\T} $ 
that are connected to $ i $ and $ i+1 $.
We now get, using \eqref{closing}
\begin{equation}\label{inserting}
  f_\T \UU_i= \raisebox{-.5\height}{\includegraphics[scale=0.7]{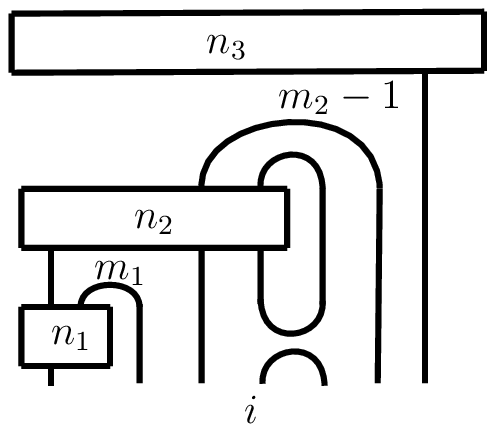}} = 
\raisebox{-.5\height}{\includegraphics[scale=0.7]{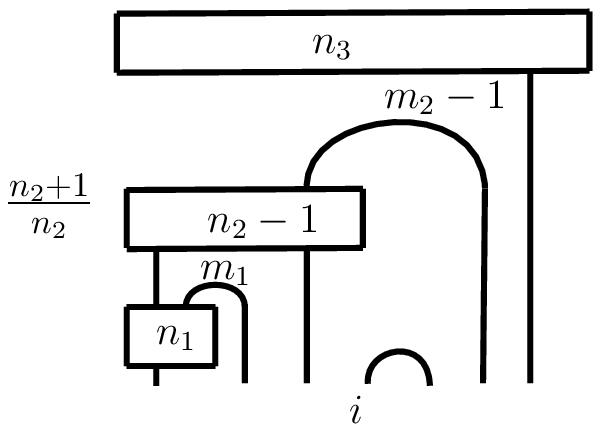}}   
\end{equation}  
On the other hand, bending down the last top line of the recursive formula \eqref{goesbackto}
for $ \JWn$ we have
\begin{equation}\label{bendingdown}
\raisebox{-.5\height}{\includegraphics[scale=0.7]{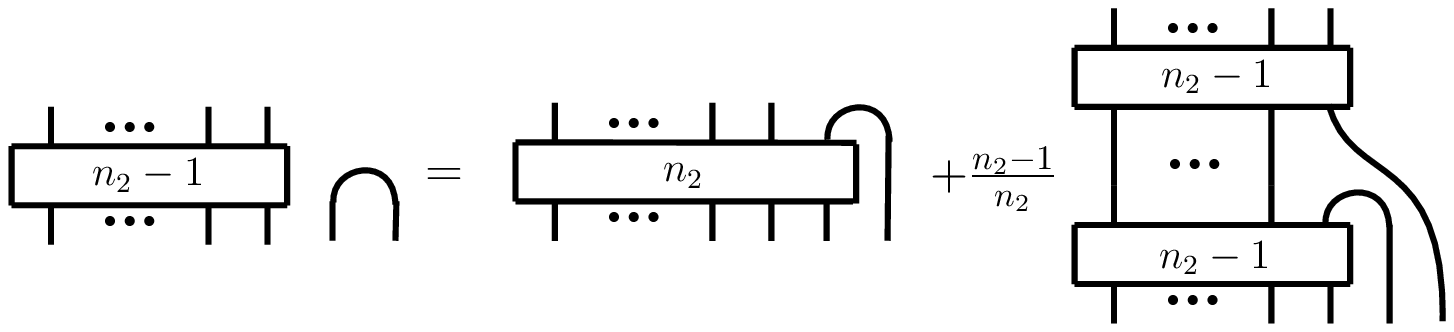}}   
\end{equation}
and inserting this in the right hand side of \eqref{inserting} we obtain
\begin{equation}\label{followsfrom}
\begin{array}{ll}  {\color{black}{  f_\T \UU_i }} &= 
   \raisebox{-.35\height}{\includegraphics[scale=0.75]{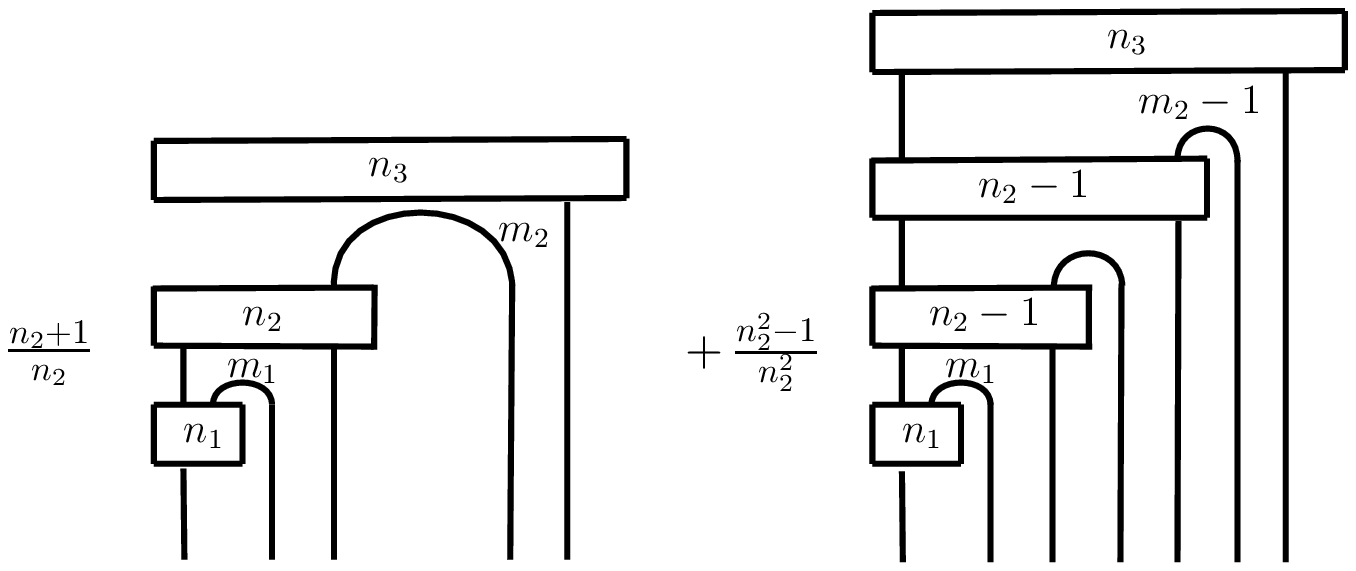}} \\& \\
   &= \dfrac{ n_2+1}{ n_2} f_{\T_u}+ \dfrac{ n_2^2 -1}{ n_2^2} f_{\T_d}
\end{array}
\end{equation}
One finally checks that $ n_2 =  c_{\T_u}(i) - c_{\T_d}(i)  = r $ and so
{\bf a)} follows 
from \eqref{followsfrom}, at least for $ \T $ as in \eqref{2.17}.
For general $ \T$ the proof of {\bf a)} is carried out the same way. 
From this {\bf b)} follows by applying $ \UU_i $ to both sides of {\bf a)}.

\medskip
Finally, 
{\bf c)}
{\color{black}{and {\bf d)}}}
are direct consequences of the definitions, 
{\color{black}{with {\bf c)} corresponding to $\UU_i $ annihilating a Jones-Wenzl element,
and 
{\bf d)} to $\UU_i $ acting on a cap.}}
\end{dem}

\medskip
Theorem \ref{YSFfirst} is an analogue of Young's seminormal form known from the representation
theory of $ \QQ \Si_n $. 
To make this explicit we set 
\begin{equation}
\sss_i := \Phi(s_i ) = \UU_i - \one
\end{equation}
Then we have the following Corollary to Theorem \ref{YSFfirst}. 
\begin{corollary}\phantomsection\label{YSFsecond}
  {\rm (Young's seminormal form YSF for $\TLnQ$)}.
  Let $ \T, s_i, \T_u, \T_d $ and $ r $ be as in Theorem \ref{YSFfirst}. Then we have 
  \begin{description}
  \item[a)]
$  f_{\T_d} \sss_i=
\dfrac{ 1 }{ r }
f_{\T_d} + \dfrac{ r^2 -1 }{ r^2 } f_{\T_u}$ 
  \item[b)]
    $  f_{\T_u} \sss_i= -\dfrac{ 1 }{ r }
f_{\T_u} +  f_{\T_d}
$
  \end{description}
Suppose next that $ \T s_i \notin \std(\lambda) $. Then 
  \begin{description}
  \item[c)]
    $  f_{\T} \sss_i= -f_{\T}  \,\, \, \, \, \,   \mbox{ if } i, i+1 \mbox{ are in the same column of } \T $ 
 \item[d)]
    $  f_{\T} \sss_i= f_{\T}    \,\, \, \,\, \, \,\, \, \, \mbox{ if } i, i+1 \mbox{ are in the same row of } \T $ 
    \end{description}
\end{corollary}  

\begin{dem}
This follows immediately from Theorem \ref{YSFfirst}. 
\end{dem}

\medskip

 \begin{remark}\label{YSFremark}
  \normalfont
{\color{black}{
Note that the main ingredient for proving Theorem \ref{YSFfirst}, and hence also 
Corollary \ref{YSFsecond}, was the 
recursive formula \eqref{goesbackto}
for $ \JWn$. In fact, \eqref{goesbackto} may be viewed as a special case of Theorem \ref{YSFfirst}. 
Indeed, setting $ \lambda=(2, 1^{n-2}) $ and letting $ \T= \T^{\lambda}s_{n-1}$
we have 
\begin{equation}f_{\T}= \raisebox{-.42\height}{\includegraphics[scale=0.75]{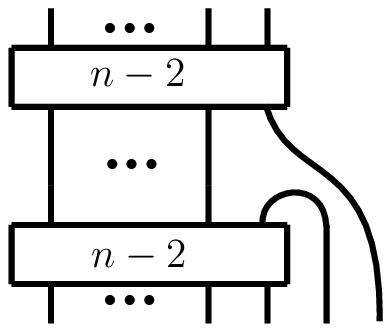}}
\end{equation}
Moreover $      \T s_{n-1}  \unlhd \T  $ and so {\bf b)} of Theorem \ref{YSFfirst}
is the formula $ f_{\T} \UU_{n-1} = \frac{n-2}{n-1} f_{\T} + f_{\T s_{n-1}} $, that is 
\begin{equation} \raisebox{-.42\height}{\includegraphics[scale=0.75]{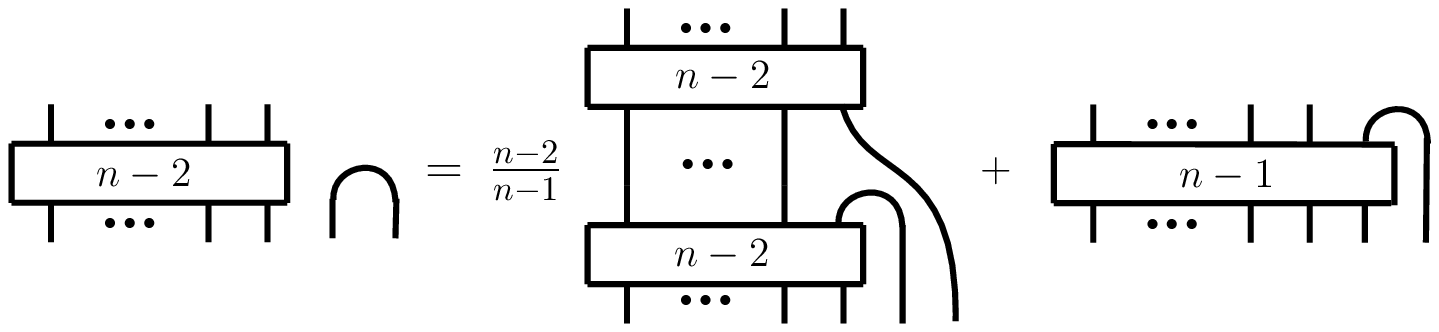}}
\end{equation}
After bending up the last line, this becomes \eqref{goesbackto} for $n-1$. 
In view of this one may consider \eqref{goesbackto} and YSF,
that is Corollary \ref{YSFsecond}, as two sides of the same coin. 
}}
\end{remark}

\medskip
We next aim at proving that $ f_\T$'s is an eigenvector for $ \LL_i$ with eigenvalue $ c_\T(i) $.
The argument for this will be an induction on $ \std(\lambda)$ over $\trianglelefteq $.
We may either carry out this induction from top to bottom, using
$ \T^{\lambda} $ as inductive basis, or from bottom to top, using $ \T_{\lambda} $ as inductive basis.
In either case it turns
out that the inductive step, using Theorem \ref{YSFfirst}, is relatively straightforward and
similar to the inductive step for the $ \QQ \Si_n $-case, whereas the
inductive basis is the most complicated part of the proof.
The $ \T_{\lambda} $-case is slightly simpler than the $ \T^{\lambda} $-case
and so we choose to carry out the induction from bottom to top. In other words, to prove the inductive basis
we should take $ \T= \T_{\lambda} $ where $ \lambda = \ParTwo $ and must show that
$ f_{\T_{\lambda}} \LL_i = c_{\T_{\lambda}}(i) f_{\T_{\lambda}}$ for all $ i =1, 2, \ldots, n$. This is the content(!)
of the next Lemma. 
\begin{lemma}\label{2.4}
  Let the situation be as just described, that is $ \T= \T_{\lambda} $
  where $ \lambda = (2^{l_2} , 1^{l_1-l_2} )  \in \ParTwo$ and $ l_1 $ and
$ l_2 $ are the lengths of the two columns of $ \lambda$. 
  Then
  we have that $   f_{\T} \LL_{i} = (1-i)  f_{\T} $  for $ i=1,\ldots, l_1 $ 
  and $   f_{\T} \LL_{i} = (2 -i + l_1 )  f_{\T} $
  for $ i = l_1+1, \ldots, n$, that is 
  \begin{equation}\label{thatis}
  f_{\T} \LL_{i} = c_{\T}(i)  f_{\T}
  \end{equation}
\end{lemma}  
\begin{dem}
{\color{black}{We have 
\begin{equation}\label{2.23A}
  f_{\T } =   \raisebox{-.47\height}{\includegraphics[scale=0.7]{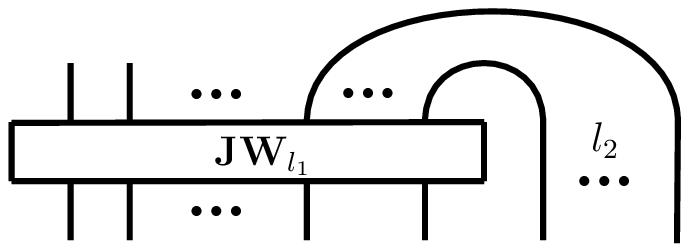}}
\end{equation}
On the other hand, from the definition of the $\JM$-elements in 
\eqref{defJM} we have the following formula, valid for $ i=1, 2, \ldots, n-1 $
\begin{equation}\label{225first}
  \LL_{i+1} = \sss_i \LL_i \sss_i + \sss_i = (\UU_i-\one) \LL_i (\UU_i-\one) +  \UU_i-\one 
\end{equation}
Since $ f_{\T }  \UU_{i } = 0 $ for $ i = 1, 2, \ldots, l_1-1 $ we get from this 
that
\begin{equation}\label{wededucefrom225}
  f_{\T }  \LL_{i } = (1-i ) f_{\T }  \mbox{ for } i= 1,2\ldots, l_1
\end{equation}  
which shows \eqref{thatis} for these value of $ i $.

\medskip
Now, if we assume that \eqref{thatis} also holds for $ i = l_1+1 $,
we would deduce from \eqref{225first} that \eqref{thatis} holds for $ i = l_1+2,\ldots, n $ as well, 
since $ f_{\T }  \UU_{i } = 0 $ for $ i = l_1+1, \ldots, n-1$, and so \eqref{thatis} would have been proved
for all $i$.

\medskip
We are therefore reduced to showing \eqref{thatis} for $ i = l_1+1 $,
which is equivalent to showing that $   f_{\T }  \LL_{l_1+1 } =  f_{\T }$.
But in the diagrammatic expression for $ f_{\T }  \LL_{l_1+1 } $, that is  
\begin{equation}\label{2.23}
f_{\T } {\mathbf{L}}_{l_1+1} = 
   \raisebox{-.4\height}{\includegraphics[scale=0.7]{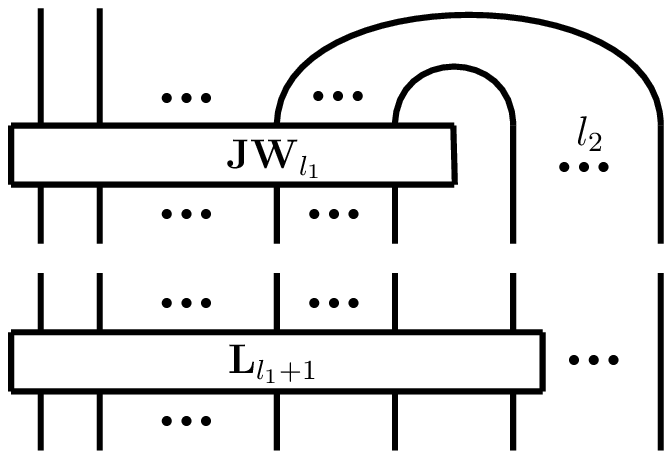}} 
\end{equation}
the multiplication with $  \LL_{l_1+1 } $ 
only involves the leftmost $ l_1 +1 $ bottom lines of $ f_{\T }  $ and so we 
may assume that $ l_2 =1$ when proving
$   f_{\T }  \LL_{l_1+1 } =  f_{\T }$.
We therefore proceed to prove by induction on $ l_1 $ that 
$   f_{\T} \LL_{l_1+1} =  f_{\T} $ where $ \lambda= (2, 1^{l_1 -1}) $. }}

\medskip
For this the basis case $ l_1=1 $ is the claim that  
\begin{equation}
 \raisebox{-.5\height}{\includegraphics[scale=0.75]{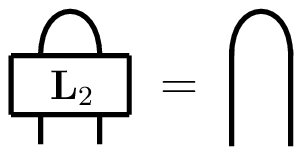}} 
\end{equation}
or equivalently that $ \UU_1 \LL_{2} =  \UU_1 ( \UU_1 - \one )$ 
which is immediate from the definitions.


\medskip
{\color{black}{We then treat the inductive step from $ l_1 $ to $l_1 +1$}}. 
In view of \eqref{225first}, 
we first calculate an expression for $ f_{\T }  (\UU_{l_1 }-\one ) $. 
We find 
\begin{equation} \begin{array}{l}
f_{\T }  (\UU_{l_1 }-\one ) =  \raisebox{-.6\height}{\includegraphics[scale=0.7]{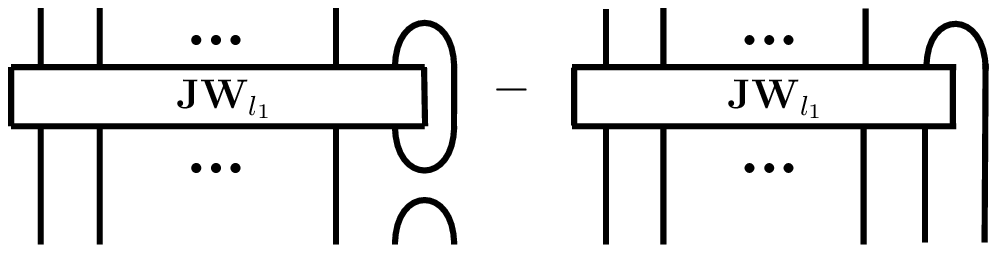}} 
\end{array}
  \end{equation}
\begin{equation} \hspace{-5.8cm}  
=  \raisebox{-.55\height}{\includegraphics[scale=0.7]{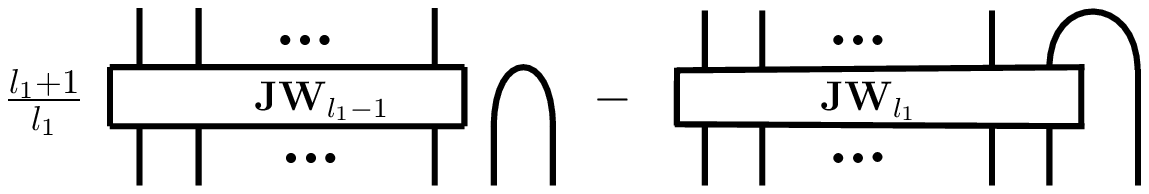}} 
  \end{equation}
\begin{equation}\label{2.27}  
=  \raisebox{-.38\height}{\includegraphics[scale=0.7]{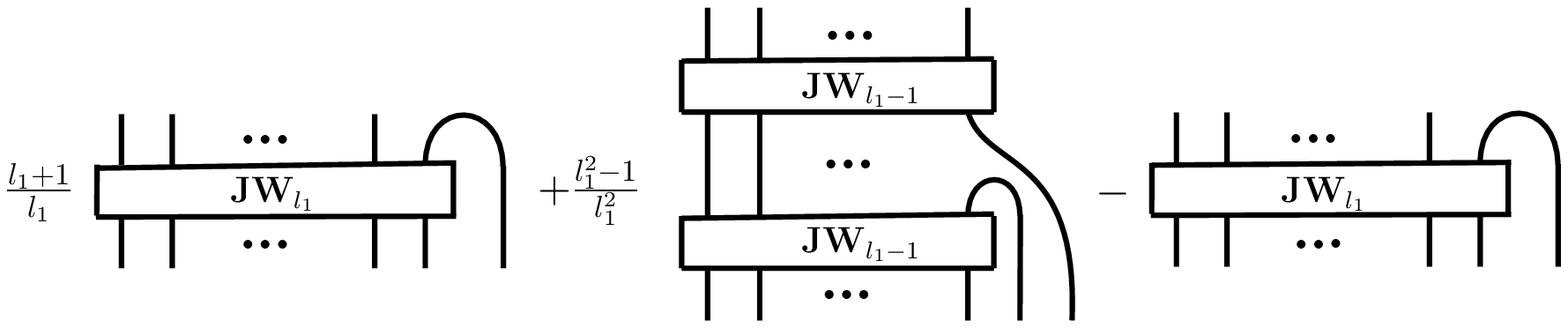}} 
  \end{equation}
\begin{equation}\label{2.28}
\hspace{-5.2cm}
=  \raisebox{-.38\height}{\includegraphics[scale=0.7]{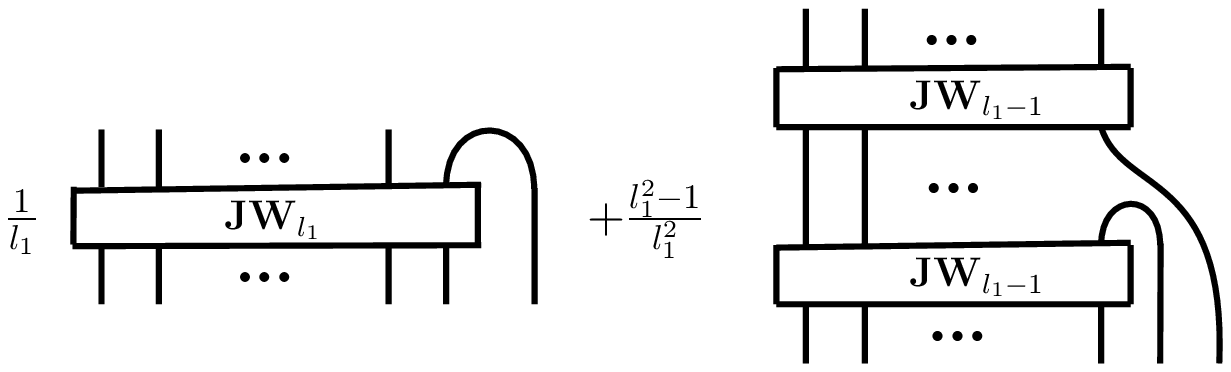}}   
\end{equation}
where we used the \eqref{bendingdown} variation of 
\eqref{goesbackto} for \eqref{2.27}. We next apply $ \LL_{l_1} $ to 
\eqref{2.28} in order to arrive at an expression for $  f_{\T}( \UU_{l_1}-\one) \LL_{l_1} $. 
Using \eqref{wededucefrom225} we see that 
$ \LL_{l_1} $ acts on the first term of \eqref{2.28} by multiplication with $ 1-l_1$ 
and, by inductive hypothesis, it acts on the second term of
\eqref{2.28} by multiplication with $ 1$. Combining, we get that 
\begin{equation}\label{2.29}
\hspace{-3cm}
f_{\T} (\UU_{l_1}-\one) \LL_{l_1} =  \raisebox{-.38\height}{\includegraphics[scale=0.7]{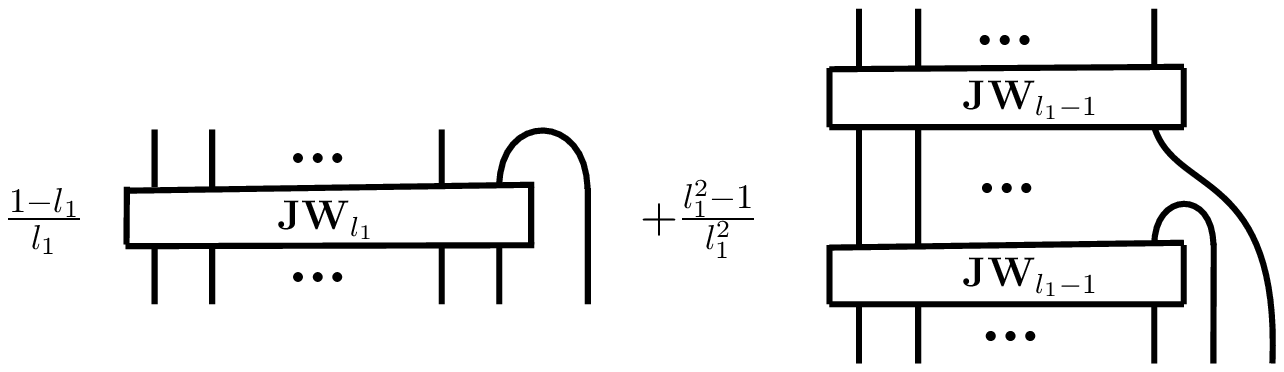}}   
\end{equation}
We now get
\begin{equation}\label{2.30}
\hspace{-3cm}
f_{\T} (\UU_{l_1}-\one) \LL_{l_1}  (\UU_{l_1}-\one)=
\raisebox{-.73\height}{\includegraphics[scale=0.7]{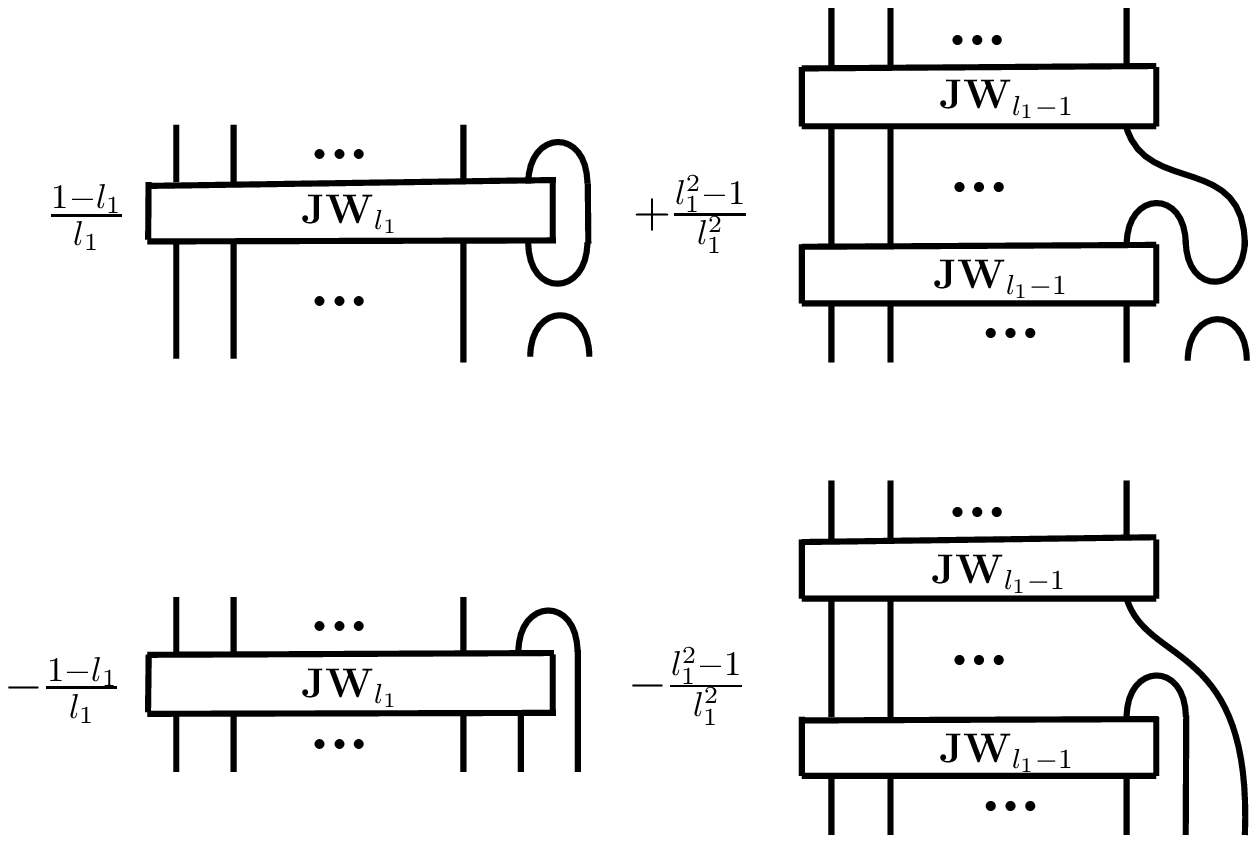}}   
\end{equation}
\begin{equation}\label{2.31}
 =  \raisebox{-.4\height}{\includegraphics[scale=0.7]{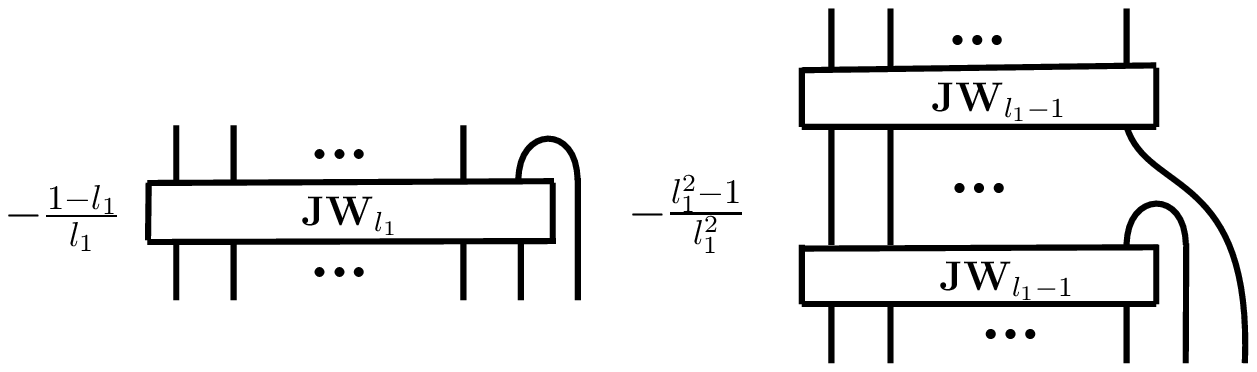}}   
\end{equation}
Finally, adding \eqref{2.28} and \eqref{2.31} we get, using \eqref{225first}
\begin{equation}\label{2.32} 
  f_{\T}  \LL_{l_1+1}
  = \raisebox{-.5\height}{\includegraphics[scale=0.7]{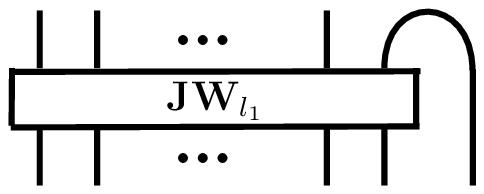}}   = f_{\T }
\end{equation}
This proves the induction step
and then also the Lemma.

\end{dem}  

\medskip

\begin{lemma}\phantomsection\label{commutation}
  We have the following commutation relations between $ \LL_k $ and $ \UU_i $. 
  \begin{description}
  \item[a)] If $ k \neq i, i+1 $ then $   \UU_i \LL_k = \LL_k \UU_i$
  \item[b)] We have $(\UU_i -\one)  \LL_i  =  \LL_{i+1} ( \UU_i -\one) -1 $
     \item[c)] We have $(\UU_i -\one)  \LL_{i+1}  =  \LL_{i} ( \UU_i -\one) +1 $
  \end{description}    
\end{lemma}  
\begin{dem}
  This follows immediately from $ \LL_k = \Phi(L_i ) $ and the definition of $ L_i $ in
  \eqref{defJM}. 

\end{dem}  

\medskip
We can now show the Theorem that was mentioned above. 
\begin{theorem}\label{mentionedabove}
  Let $ \T \in \std(\lambda) $ where $ \lambda \in \ParTwo$. Then
for all $ i =1,2, \ldots, n $ 
  we have 
  \begin{equation}\label{2.34}
  f_{\T} \LL_{i} = c_{\T}(i)  f_{\T}
  \end{equation}  
\end{theorem}  
\begin{dem}  
As already mentioned, we show the formula \eqref{2.34} by upwards induction on $ \std(\lambda)$. The basis
  case $ \T = \T_{\lambda} $ is given by Lemma \ref{2.4}, so let us assume that
 $ \T \neq \T_{\lambda} $ and that 
  \eqref{2.34} holds
  for all $ \s $ such that $ \s \lhd \T $. We must then check it for $ \T$. Since
  $ \T \neq \T_{\lambda} $ there is an $ i $ appearing in the second column of $ \T $, but with $ i+1$
  appearing in the first column of $ \T$, in a lower position, 
  and so $ \T s_i \lhd \T $. Setting $ f_d = f_{\T  s_i} $, $ f_u = f_{\T} $
  and $ r :=    c_u(i ) -  c_d(i ) $ where $  c_u(k ) := c_\T(k ) $ and $  c_d(k ) := c_{\T s_i}   (k ) $, 
  we have from {\bf a}) of Theorem \ref{YSFfirst} that
  \begin{equation}\label{YSFinduction}
 f_{d} \UU_i=\dfrac{ r +1 }{ r }
f_{d} + \dfrac{ r^2-1}{ r^2} f_{u}
  \end{equation}

By  
induction hypothesis we have that $ f_d \LL_k = c_d (k) f_{\color{black}{d}}$ for all $ k $. 
  Suppose first that $ k \neq i, i+1 $. Then we get from Lemma \ref{commutation} that
  $ \UU_i \LL_k =  \LL_k \UU_i   $. Acting upon $ f_d $, this equation becomes via
\eqref{YSFinduction}
 \begin{equation}\label{compare}
\dfrac{ r +1 }{ r }
c_d(k)   f_{d} + \dfrac{ r^2-1}{ r^2}  f_{u} \LL_k = \dfrac{ r +1 }{ r }
c_d(k) f_{d} + \dfrac{ r^2-1}{ r^2} c_d(k) f_{u}
 \end{equation}
from which we deduce 
that $ f_u \LL_k = c_d(k ) f_u $. But in this case $ c_u(k) = c_d(k ) $ and so $ f_u \LL_k = c_u(k ) f_u $,
as claimed.

Suppose now that $ k = i $.
We have from Lemma \ref{commutation}
that $ (\UU_i - \one ) \LL_i = \LL_{i+1} (\UU_i - \one ) -1 $. Acting upon $ f_d $, this becomes 
$ f_d (\UU_i - \one ) \LL_i  = f_d (\LL_{i+1} (\UU_i - \one ) -1 )  $. 
Using \eqref{YSFinduction}, the left hand side of this is 
\begin{equation}\label{compareA}
{\rm LHS}= \dfrac{ 1 }{ r } c_d(i) 
f_{d} + \dfrac{ r^2-1}{ r^2} f_{u} \LL_i^{} 
\end{equation}
whereas the right hand side is 
\begin{equation}\label{compareB}
 \begin{split}
{\rm RHS}=&  \dfrac{ 1 }{ r } c_d(i+1) 
f_{d} + \dfrac{ r^2-1}{ r^2} c_d(i+1)  f_{u} - f_d
= \dfrac{- r+ c_d(i+1) }{ r } 
f_{d} + \dfrac{ r^2-1}{ r^2} c_d(i+1)  f_{u}  \\
=&  \dfrac{ c_d(i) }{ r } 
f_{d} + \dfrac{ r^2-1}{ r^2} c_d(i+1)  f_{u} 
 \end{split}
 \end{equation}
where we used $ c_d(i+1) = c_u(i) $ and $ r =    c_u(i ) -  c_d(i ) $ for the last equality. 
Comparing \eqref{compareA} and \eqref{compareB}
we conclude that $ f_u \LL_i  = c_d(i+1 ) f_d =  c_u(i ) f_d $, proving the Theorem in this case
as well.

\medskip
Finally, the case $ k = i+1 $ is proved the same way. The Theorem is proved.
\end{dem}

\medskip

\begin{corollary}\label{finalcorsection2}
  For $ \lambda $ a two-column partition and $ \T \in \std(\lambda) $ we have that
  $ \EE^{\prime}_{\T} = \EE_{\T}  $. In particular, the $ \{  \EE_{\T}^{\prime} \} $ form
  a complete set of primitive idempotents for $ \TLnQ$. 
\end{corollary}
\begin{dem}
  It follows from Theorem \ref{mentionedabove} and
  the formula $ \EE_{\T}^{\prime} := \frac{1}{\gamma_{\T}} f_{\T \T} $ that
  $ \EE^{\prime}_{\T} \LL_k = \LL_k  \EE^{\prime}_{\T} = c_\T(k)  \EE^{\prime}_{\T} $ for all $ k$. 
  But this property characterizes the idempotent $ \EE_\T $ and so $ \EE^{\prime}_{\T} = \EE_{\T} $,
  as claimed. 
\end{dem}

\medskip
 \begin{remark}\label{Okounkovremark}
  \normalfont
      {\color{black}{In the Okounkov-Vershik theory for the representation theory of $ \QQ \Si_n $ one
          derives 
Young's seminormal form via the {\it Gelfand-Zetlin} subalgebra of  $ \QQ \Si_n $, see \cite{OkoVer}.
It should be possible to establish an analogue of this theory for $ \TLnQ $, using
our $ \LL_k $'s. It should also be possible to show that $ \EEE_{\T}^{\prime} = P_{\T}$ where
$ P_{\T} $ is a product of central idempotents as in \cite{OkoVer}. This would
      give an alternative way of proving Corollary \ref{finalcorsection2}. }}
 \end{remark}         

\section{The unseparated case}\label{The unseparated case}
We shall from now on focus on the Temperley-Lieb algebra $\TLnF $ defined over the finite
field $ \FF$, where $ p>2$. We are interested in idempotents in $\TLnF $.

\medskip
If $ p > n $ the condition \eqref{separationcondition} still holds and 
so $\TLnF $ is a semisimple algebra and in fact all the results from the
previous section remain valid. Let us therefore assume that $ p \le n $.
Under that assumption \eqref{separationcondition} does not hold,
and so we are in the {\it unseparated case}
in the
{\color{black}{terminology}} of \cite{Mat-So}.
Moreover, the coefficients of $ \JWn $ and of $ \EE_{\T} $ 
cannot be reduced from $ \QQ $ to $ \FF$, and hence these idempotents {\it do not exist} in 
$\TLnF $. In fact, if $ p \le n $
there are in general no nonzero idempotents in $\TLnF $ satisfying \eqref{definingpropertyJW}. 

\medskip
On the other hand, we can still apply the general theory
of $\JM$-elements to construct idempotents for $\TLnF $. Let us briefly explain this. 

\medskip
For $ \T \in \std(\lambda) $ where $ \lambda \in \ParTwo$ we define the {\it $p$-class} $ [ \T]  $ of $ \T$
via
\begin{equation}\label{firstclass}
[\T]  = \{ \s \in \std(\ParTwo) \, | \, c_{\s}(i) \equiv c_{\T}(i)  \mbox{ mod }  p \mbox{ for all } i=1,2,\ldots, n  \} 
\end{equation}
We now set 
\begin{equation}\label{mainprinciple}
\EE_{[\T]} := \sum_{ \s \in [\T]  } \EE_{\s}  
\end{equation}
By definition $ \EE_{[\T]} \in \TLnQ $, but 
it follows from the general theory developed in \cite{Mat-So} that $ \EE_{[\T]} $ in fact belongs to 
$ \TLZpn $ where $ \Z_{(p)} :=\{ \dfrac{a}{b} \in \QQ \, | \,   p
\mbox{ does not divide } b \} $. 
We have that $  \Z_{(p)} $ is a local ring with maximal ideal $ \pi := p \Z_{(p)} $ and
$   \Z_{(p)}/\pi \cong \FF$.
and hence $  \EE_{[\T]}  $ can be reduced to an element of $ \TLnF $, that we shall also denote 
$   \EE_{[\T]}  $. 

\medskip
The $   \EE_{[\T]}  $'s clearly are idempotents in $ \TLnF $, called {\it class idempotents}, 
but they are 
not primitive idempotents in general, as we shall shortly see. 

\medskip
Let $ M_{triv} := \Delta^{\FF}( 1^{n} ) $
be the trivial $ \TLnF$-module, in other words $ M_{triv} $ is the one-dimensional $ \TLnF$-module on which
$ \UU_k $ acts as zero for all $ k $. 
Let $ P_{triv} $ be the projective cover for $ M_{triv} $. 
By general principles there exists 
a primitive idempotent $ ^{p} \EE_{triv} \in \TLnF$ such that $   ^{p} \EE_{triv} \TLnF = P_{triv} $.
Recently, it was  observed in \cite{StuSpe} that the idempotent $ ^{p} \EE_{triv} $ coincides with the 
{\it $p$-Jones-Wenzl idempotent} $ ^{p}\!\JWn  $ that was introduced by Burull, Libedinsky and Sentinelli,
see  
\cite{BLS}. We need this fact in the following, 
and shall therefore recall the definition of $ ^{p}\!\JWn  $.

\medskip
For $ n \in \N$ we define non-negative integers $ a_i $ satisfying $ 0 \le a_i < p, a_k \neq  0 $
and 
\begin{equation}\label{padic}
 n+1 = a_k p^k + a_{k-1} p^{k-1} + \ldots +a_1 p + a_0 
\end{equation}
In other words, $ (a_k, a_{k-1}, \ldots, a_1, a_0 )$ are the coefficients
of $ n+1 $ when written in base $p$. We then define $ {\mathcal I}_n  \subseteq \N $ via 
\begin{equation}\label{mathcal I}
{\mathcal I}_n := \{ a_k p^k \pm  a_{k-1} p^{k-1} \pm  \ldots  \pm a_1 p  \pm a_0 \} -1
\end{equation}
{\color{black}{where for $ A \subseteq \N$ we define $ A - 1 := \{ a-1 \, | \, a \in A\}.$}}
One checks that each $ m \in {\mathcal I}_n $ is given uniquely by the corresponding
sequence of signs for the nonzero
$ a_k$'s. Using this, 
for $  m \in {\mathcal I}_n $ 
we now define a tableau $ \T_m \in \std(\ParTwo) $
in terms of a block decomposition for standard tableaux as in \eqref{intheform}, 
using 
blocks $ D_1, M_1, D_2, M_2, \ldots, D_k, M_k $ of consecutive numbers, as follows.

\medskip
Suppose first that $ i_1 \ge  0 $ is maximal such that $ (a_k, a_{k-1}, \ldots, a_{k-i_1}) $ all
appear in $ m $ with non-negative sign. Then $ D_1 $ is defined by the condition that it be of cardinality 
$ | D_1|  = a_k p^k  + \ldots + a_{k-i_1} p^{k-i_1}  -1 $. 
Suppose next that $ i_2 > i_1 $ is maximal such that
$ (a_{k-i_1-1}, a_{k-i_1-2}, \ldots, a_{k-i_2 }) $ all appear in $ m $ with non-positive sign. Then we define
$ M_1 $ by the condition that it be of cardinality 
$ | M_1|  = a_{k-i_1-1} p^{k-i_1-1} +  \ldots + a_{k-i_2} p^{k-i_2}   $. We then continue the same way, defining
$ D_2, M_2, \ldots $ 
except that the $ -1 $ term should only appear for $ D_1$. 

\medskip
The $p$-Jones-Wenzl idempotent $ ^{p}\!\JWn  $ is now defined as follows
\begin{equation}\label{pJW}
 ^{p}\!\JWn  := \sum_{ m \in {\cal I}_n } \EE_{ \T_m}^{\prime} 
\end{equation}

Note that, unlike the definition in \eqref{pJW},
the original definition of $  ^{p}\!\JWn  $
in \cite{BLS} was formulated recursively. {\color{black}{The definition in \eqref{pJW} 
is the left-right mirror of Definition 2.22 in \cite{TuWe1}, although we have here formulated it 
in terms of 
standard tableaux.}}
Note also that the original definition of $  ^{p}\!\JWn  $, 
{\color{black}{and the definition in \cite{TuWe1},}}
was carried out for the Temperley-Lieb
algebra with loop parameter $ -2 $, as opposed to loop parameter $ 2 $ as in the present paper. 
To switch between the two settings one
should apply the isomorphism $ \UU_i \mapsto - \UU_i $. 

\medskip
Let us give a couple of examples. If $ n= 3 $ and $ p= 3 $ we have $ {\cal I }_3 = \{ 3 \pm 1 \} -1 = 
\{ 3,1 \} $. The tableaux corresponding to the elements of $  {\cal I }_3 $ are as follows
\begin{equation}\label{dib51}
  \T_3 =  \raisebox{-.4\height}{\includegraphics[scale=0.7]{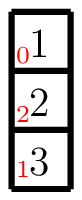}}, \, \, \, \, 
    \T_1 = \raisebox{-.45\height}{\includegraphics[scale=0.7]{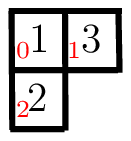}}
\end{equation}
and so we get 
\begin{equation}\label{dib53}
  ^{3}\!\JWtres     =\EE_{\T_3}^{\prime} + \EE_{\T_1}^{\prime}   = 
  \raisebox{-.42\height}{\includegraphics[scale=0.7]{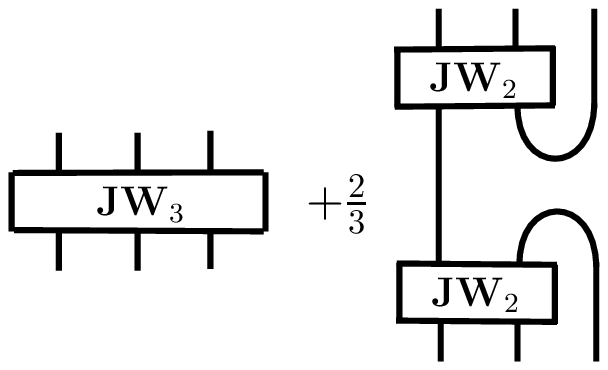}}
\end{equation}
To verify that  $ ^{3}\!\JWtres $ belongs to $ \TLtrestres$, one 
uses \eqref{dib7} and \eqref{dib8} to expand $ \JWdos $ and $ \JWtres $ and finds
\begin{equation}\label{dib54}
 ^{3}\!\JWtres     = \raisebox{-.42\height}{\includegraphics[scale=0.7]{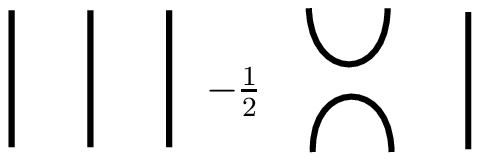}}
\end{equation}
which indeed belongs to $ \TLtrestres$.

\medskip

In the tableaux in \eqref{dib51} we have indicated with color red, for each $ i =1,2,3 $, 
the {\it residue} $ c_{\T}(i) \, \,  {\rm   mod  }\, \, p $
of the content $  c_{\T}(i) $. Using this
we get that the $3$-class of $ \T_3 $ is $ [ \T_3] = \{ \T_3, \T_1 \} $.
We now 
use Corollary
\ref{finalcorsection2} and get that 
\begin{equation}
\EE_{[\T_3]}  = \,  ^{3} \! \JWtres     
\end{equation}
Thus in this case
the class idempotent $ \EE_{[\T_3]} $ is in fact primitive. 

\medskip
To give an example where the class idempotent is not primitive we
choose $ p= 3 $ and $ n=12$. 
We then have $ n+1 = 9+3+1 $ and so $ {\cal I}_{n} = \{ 9 \pm 3 \pm 1 \} -1 = \{ 12, 10, 6,4 \} $
and so we have that $    ^{3} \! \JWdoce =
\EE_{\T_{12}}^{\prime} + \EE_{\T_{10}}^{\prime} +
\EE_{\T_{6}}^{\prime} + \EE_{\T_{4}}^{\prime}  $ 

\medskip
The corresponding standard tableaux, with $3$-residues indicated with color red as before, are
as follows
\begin{align}\label{redcolortableaux1}
&    \T_{12}=  \raisebox{-.5\height}{\includegraphics[scale=0.7]{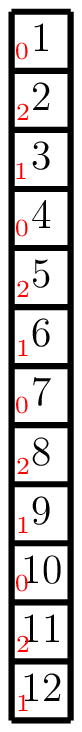}}  & 
&   \T_{10}=  \raisebox{-.455\height}{\includegraphics[scale=0.7]{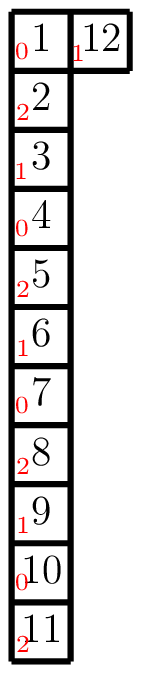}} & 
&      \T_{6}=  \raisebox{-.34\height}{\includegraphics[scale=0.7]{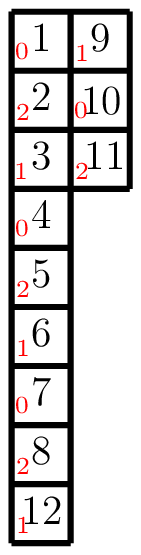}} & 
&   \T_{4}=  \raisebox{-.26\height}{\includegraphics[scale=0.7]{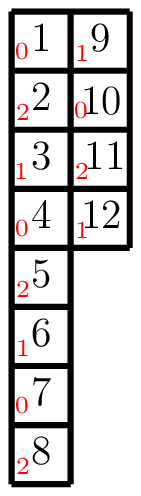}}
 \end{align}
Note that
$ \T_{12}, \T_{10}, \T_{16}$ and $ \T_{4} $ all belong to the same $3$-class, as can be seen
by comparing the residues modulo $3$. But the class $ [\T_{12}] $ 
contains two more tableaux, namely 
\begin{align}\label{redcolortableaux2}
&     \s=  \raisebox{-.50\height}{\includegraphics[scale=0.7]{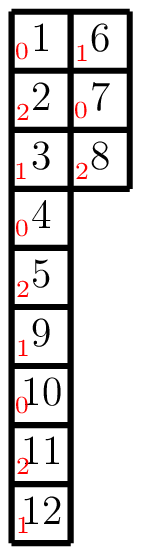}} & 
&   \T=  \raisebox{-.44\height}{\includegraphics[scale=0.7]{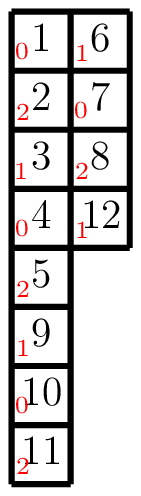}}
\end{align}
obtained by interchanging $ \{6,7,8\} $ and $ \{9,10,11\} $ in $ \T_6 $ and $ \T_4$. From this we get that
\begin{equation}\label{notprimitive}
 \EE_{[\T_{12}]} = \,  ^{3} \!  \JWdoce  + \EE_{\s}^{\prime} +
 \EE_{\T}^{\prime}
\end{equation}
which shows that
$ \EE_{\s}^{\prime} +
\EE_{\T}^{\prime} \in \TLnF $. By expanding in terms of the diagram basis
for $ \TLnF $, one gets $  \EE_{\s}^{\prime} +
 \EE_{\T}^{\prime} \neq 0 $ in $ \TLnF $ and clearly
${\color{black}{  ^{3} \!}}  \JWdoce $ and $  \EE_{\s}^{\prime} +
 \EE_{\T}^{\prime}  $ are orthogonal. 
Hence  
$ \EE_{[\T_{12}]} $ is not a primitive idempotent in $ \TLnF $.  

\medskip
The purpose of the rest of the paper is to show that a variation of the principle for constructing
idempotents given in \eqref{mainprinciple}, this time using KLR-theory, 
can be applied recursively to derive the $p$-Jones-Wenzl idempotents
for $ \TLnF$, that is the primitive idempotents. 

\medskip
Let us start out by proving the following Lemma, which is a generalization of
\eqref{notprimitive}.
\begin{lemma}
  Let $ \EE_{[\T_n]} \in \TLnF  $ be the class idempotent for the $p$-class $ [ \T_n] $, given by
  the one-column tableau $ \T_{n} = \T_{ 1^n} = \raisebox{-.5\height}{\includegraphics[scale=0.5]{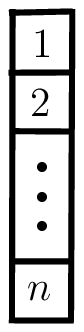}}$.
  Then $ \EE_{[\T_n ]} = \,   ^{p}\!\JWn + \EE $ for some idempotent $ \EE $
  in $ \TLnF $, orthogonal to $  ^{p}\!\JWn $. 
\end{lemma}
\begin{dem}
  We must show that $ \T_m \in [ \T_{n}] $ for all $ m \in{ \cal I }_n $ as in \eqref{mathcal I}.
  Let $ D_1, M_1, \ldots, D_k, M_k $ be the sequence of blocks defining $ \T_m $,
  as in the paragraph preceding \eqref{pJW}. Then clearly $ c_{\T_{n}}(i) \equiv  c_{\T_m}(i)  \mbox{ mod } p$ for
  $ i \in D_1 $, since in fact $ c_{\T_{n}}(i) =  c_{\T_m}(i) $ for these $ i $.
  Suppose now that $ M_1 \neq \emptyset $ and that $ m_{1, min} $ is the first number in $M_1 $. Then
  by the cardinality of $ D_1 $ we have that
  $ c_{\T_{n}}( m_{1, min})  \equiv  c_{\T_m}( m_{1, min}) \equiv 1  \mbox{ mod } p$ and
  then $ c_{\T_n}(m)  \equiv  c_{\T_m}(m)  \mbox{ mod } p$ for all $ m \in M_1 $. 
  This patterns repeats itself. If
  $ D_2 \neq \emptyset $ we let $ d_{2,min} $ be the first number of $ D_2$. Then by the cardinality of
  $  D_1 \cup M_1 $
  we have that $ c_{\T_n}(d_{2,min})  \equiv  c_{\T_m}(d_{2,min}) \equiv 1  \mbox{ mod } p$
  and then $ c_{\T_n}(d)  \equiv  c_{\T_m}(d)  \mbox{ mod } p$ for all $ d\in D_2 $, and so
  on recursively. This
  proves the Lemma. 
\end{dem}

\medskip

For the rest of the article we fix {\color{black}{integers}}
$ \NNone, \NNtwo, \RR $ using integer division as follows 
\begin{equation}\label{fixnotation}
\NN= (p-1) + \NNone, \,  \NNone = p \NNtwo + \RR \, \mbox{ where }  0 \le \RR < p
\end{equation}
Recall that $ n \ge p $
{\color{black}{and so $ n_2 $ is non-negative.}}

\medskip
The next Lemma gives us a kind of recursive description of the class $ [ \T_{\NN}] $.

\begin{lemma}\label{kindofrecursive}
If $ \RR = 0 $ there is a bijection
  \begin{equation}\label{313}
f_1: [ \T_{\NN}] \rightarrow \std(\ParTwoNuno) 
  \end{equation}
Otherwise, if $ \RR > 0 $, there is a bijection 
  \begin{equation}\label{314} 
f_2: [ \T_{\NN}] \rightarrow \std(\ParTwoNuno) \times \{1,2\}
  \end{equation}
\end{lemma}
\begin{dem}
  Suppose first that $ \RR=0$ and let $ \T \in [ \T_{\NN}] $. We must define $ f_1(\T) $
  and must show that it is a bijection.
Since $  \T \in [ \T_{\NN}] $, 
the numbers
$ (1,2,\ldots, p-1) $, whose content residues in $ \T$ are $ (0,p-1,\ldots, 2) $, 
all appear in the first column 
of $ \T$. We now consider consecutive blocks of consecutive numbers
$ B_1, B_2, \ldots, B_{\NNtwo} $ in $ \T$, 
all of length $ p $, starting with the block
$ B_1 := (p,p+1, \ldots, 2p-1) $. 
For each $B_i$, the 
content residues are $ (1,0,p-1,p-2,\ldots, 3,2)$. The numbers of each $ B_i $ may appear in either column
of $ \T$, 
but they all 
appear in the same column of $ \T$, since $ \T \in [ \T_{\NN}] $.
Using this observation, 
we can define $ f_1(\T) $ as the two-column standard
tableau {\color{black}{of $\NNtwo$}} that has $ i $ in the
first column iff the numbers of $ B_i $ are in the first column of $ \T$.

Here are two examples 
of $ f_1(\T) $, using $ p=3 $, in which we have indicated the blocks $ B_1, B_2, B_3 $ and $ B_4 $ with colors.

\begin{equation}
f_1:  \raisebox{-.5\height}{\includegraphics[scale=0.7]{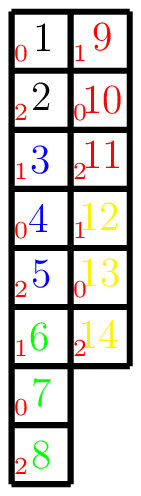}} \mapsto
      \raisebox{-.4\height}{\includegraphics[scale=0.7]{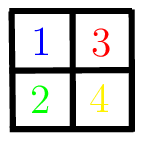}}, \, \, \,  \, \, \,  \, \, \,  \, \, \, 
f_1:  \raisebox{-.5\height}{\includegraphics[scale=0.7]{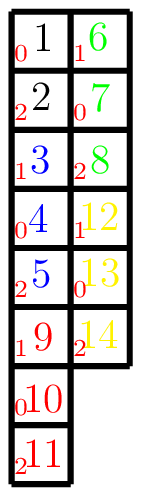}} \mapsto
      \raisebox{-.4\height}{\includegraphics[scale=0.7]{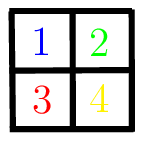}} 
\end{equation}
One readily checks that $ f_1 $, defined this way, is a bijection, proving \eqref{313}.

\medskip
In order to show \eqref{314}, we choose $ \T \in  [\T_{\NN}] $ and proceed as before, defining blocks
$ B_1, B_2, \ldots, B_{\NNtwo} $ of consecutive
numbers of length $ p $. But since $ \RR > 0 $ there will this
time be an
{\color{black}{\lq extra\rq\ }}block $ B_{\color{black}{\NNtwo+1}} $ of length $ \RR $.
The numbers of $ B_{\color{black}{\NNtwo+1}} $ may appear in
either column of $ \T$, but they all appear in the same column. 
Let $ \T_1 :=  \T |_{\le \NN-\RR } $. We now define $ f_2(\T) := (f_1(\T_1), 1) $ if the 
numbers of $  B_{\color{black}{\NNtwo+1}} $ are all in the
first column of $\T$, and otherwise we define $ f_2(\T) := (f_1(\T_1), 2) $. Here are two examples, using
$ p= 3$ and $ \RR=2 $.

\begin{equation}
f_2:  \raisebox{-.5\height}{\includegraphics[scale=0.7]{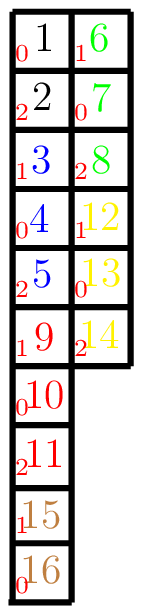}} \mapsto
     \left( \raisebox{-.4\height}{\includegraphics[scale=0.7]{dib64.eps}},1\right), \, \, \,  \, \, \,  \, \, \,  \, \, \, 
f_2:  \raisebox{-.5\height}{\includegraphics[scale=0.7]{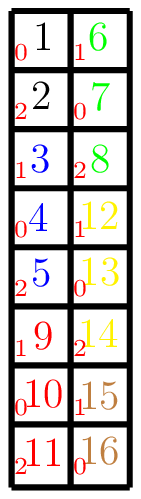}} \mapsto
    \left(  \raisebox{-.4\height}{\includegraphics[scale=0.7]{dib64.eps}}, 2 \right)
\end{equation}
    {\color{black}{Note that if $ n_2 = 0 $, corresponding to $ n+1 = p +r $, one has $ f_2(\T_1 ) = \emptyset
    \in \std({\rm Par}^{\le 2}_0)$}}.

\medskip
Once again, one checks that $ f_2 $ is a bijection, which proves \eqref{314}, and hence the Lemma. 

\end{dem}

\medskip
Returning to the examples \eqref{redcolortableaux1} and \eqref{redcolortableaux2},
where $ \NN=12 $ and $ p = 3$, 
we have that $ [ \T_{12}] = \{ \T_{12}, \T_{10}, \T_{6}, \T_{4}, \s, \T \} $ and writing $ f= f_2 $ we get
\begin{equation}\label{tableauxclassA}
  f(\T_{12}) =     \left( \raisebox{-.45\height}{\includegraphics[scale=0.7]{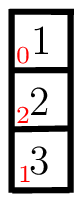}},1\right), \, 
  f(\T_{10}) =     \left( \raisebox{-.45\height}{\includegraphics[scale=0.7]{dib68.eps}},2\right), \, 
  f(\T_{6}) =     \left( \raisebox{-.45\height}{\includegraphics[scale=0.7]{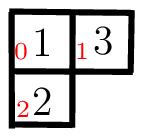}},1\right), \, 
    f(\T_{4}) =     \left( \raisebox{-.45\height}{\includegraphics[scale=0.7]{dib69.eps}},2\right)
\end{equation}
whereas 
\begin{equation}\label{tableauxclassAB}
  f(\s) =     \left( \raisebox{-.45\height}{\includegraphics[scale=0.7]{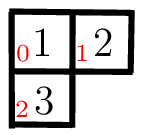}},1\right), \, 
  f(\T) =     \left( \raisebox{-.45\height}{\includegraphics[scale=0.7]{dib70.eps}},2\right) \, 
\end{equation}
Note now that $ [ \T_{3} ] =\left\{ 
\raisebox{-.45\height}{\includegraphics[scale=0.7]{dib68.eps}},
\raisebox{-.45\height}{\includegraphics[scale=0.7]{dib69.eps}} 
\right\} $, which are the two tableaux that appear in \eqref{tableauxclassA}, but that 
$ \raisebox{-.45\height}{\includegraphics[scale=0.7]{dib70.eps}} $ does not belong
to $ [ \T_{3} ]  $. Our second goal is to explain, in general,
that this is the reason why the tableaux $ \s $ and $ \T $ should
not be taken into account when giving the primitive idempotent.

\section{The integral KLR-algebra}\label{The integral KLR-algebra}
Brundan-Kleshchev and independently Rouquier found a new presentation
for the group algebra $ \FF \Si_n $, proving that it is isomorphic to the {\it KLR-algebra} $ \RKLR $
(in fact they worked in the greater generality of cyclotomic Hecke algebras). 
If $ n\ge p $ it follows from
their work that $  \FF \Si_n $ is endowed with a non-trivial $ \Z$-grading, since
$ \RKLR $ is endowed with a non-trivial $ \Z$-grading in that case.
The isomorphism $ \FF \Si_n \cong \RKLR $ is important to us since it induces, via
Lemma \ref{wellknownfunda}, an isomorphism $ \TLnF \cong \RKLR/ {\mathcal I}_{n} $
where $ {\mathcal I}_{n} $ is a graded ideal in $ \RKLR $, and hence in particular 
$ \TLnF  $ inherits a $\Z$-grading from $ \FF \Si_n  $, see \cite{PlazaRyom} for more details on this.

\medskip
Hu and Mathas gave in \cite{hu-mathas2} a new simpler proof of the Brundan-Kleshchev and Rouquier isomorphism
using seminormal forms, and via this they were able to lift it to an
isomorphism $  \Z_{(p)} \Si_n \cong \RKLRZ $, where $ \RKLRZ $
is an integral version of $ \RKLR $ (once again the result was proved in the greater generality
of cyclotomic Hecke algebras).
We shall need this isomorphism and its proof so let us recall
the precise definition of $ \RKLRZ $ from \cite{hu-mathas2}.

\medskip

We first arrange the elements of $ \FF= \{0,1,2,\ldots, p-1 \} $ in a cyclic quiver as follows
\begin{equation}
  \raisebox{-.45\height}{\includegraphics[scale=0.7]{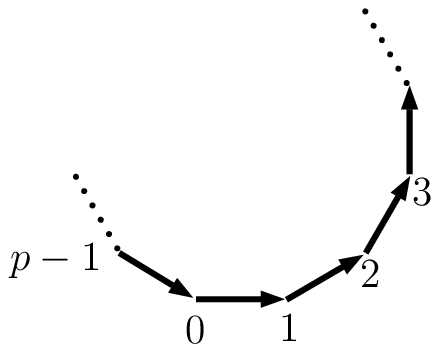}}
 \end{equation}
and for $ i, j \in \FF $ we 
write $ i \rightarrow j $ if $ i $ and $ j $ are adjacent in the quiver in the way that the arrows
indicate. We shall refer to the elements $ \ii
=(i_1, i_2, \ldots, i_n ) $ of $    \FF^n $ as residue sequences.

\begin{definition}\label{RnZ}
  The integral KLR-algebra $ \RKLRZ $ is the $  \Z_{(p)} $-algebra generated
  by the elements
\begin{equation}  
 \{ e( \ii ) \,  | \,  \ii \in  \FF^n \} \cup \{ \psi_k \,  | \,  1 \le k < n \} \cup
  \{  y_l \,  | \,  1 \le l \le n  \} 
\end{equation}    
with identity 
$   1  = \sum_{\ii \in  \FF^n } e(\ii) $,
subject to the relations
\begin{align}
\label{firstR} e(\ii) e(\jj) & = \delta_{\ii, \jj} e(\ii)  &   y_l e(\ii) & = e(\ii) y_l  \\
  \psi_k e(\ii) & =  e(\ii \cdot s_k) \psi_k    &   y_l y_m &= y_m y_l \\
  \psi_k  y_{k+1} e(\ii) & = (y_k \psi_k + \delta_{i_k, i_{k+1}}) e(\ii)     &
  y_{k+1} \psi_k e(\ii) & = (\psi_k y_k  + \delta_{i_k, i_{k+1}} )e(\ii)       \\
  \psi_k y_l & =  y_l \psi_k     &    \mbox{ if }&  l  \neq k, k+1 \\
   \psi_k \psi_m   & =   \psi_m \psi_k      &    \mbox{ if }&      | {\color{black}{k}}-m |  > 1 \\ 
   e(\ii)  & = 0      &    \mbox{ if } &  i_1  \neq 0
   \\
      y_1 e(\ii) &= 0  & &  
\end{align}
\begin{align}
      \big(\psi_k\psi_{k+1}\psi_k-\psi_{k+1}\psi_k\psi_{k+1}\big)e(\ii)=
        \begin{cases}
          -e(\ii)&\text{if }i_{k+2}=i_k\rightarrow i_{k+1}\\
          e(\ii)&\text{if }i_{k+2}=i_r\leftarrow i_{k+1}\\
          0&\text{otherwise}
        \end{cases}
\end{align}

\begin{align}\label{finalrelations}
  \psi_k^2 e(\ii) =
\begin{cases}
        (y_k-y_{k+1}) e(\ii) & \mbox{if }i_k\rightarrow i_{k+1}\neq0\\
        (y_{k}+p-y_{k+1})e (\ii)& \mbox{if }i_k\rightarrow i_{k+1}=0\\
        (y_{k+1}-y_k)e(\ii)& \mbox{if }0\neq i_k\leftarrow i_{k+1}\\
        (y_{k+1}+p-y_k)e(\ii) & \mbox{if } 0=i_k\leftarrow i_{k+1}\\
        0&\mbox{if }i_k=i_{k+1}\\
        e(\ii)&\mbox{otherwise}
\end{cases}        
\end{align}
\end{definition}
\noindent
where $ \ii \cdot s_k = (i_1, i_2, \ldots, i_k, i_{k+1}, \ldots, i_n ) \cdot s_k :=
(i_1, i_2, \ldots, i_{k+1}, i_{k}, \ldots, i_n ) $. 

\medskip
It is easy to check that $ \RKLRZ \otimes_{\Z_{  { \color{black}{(}}  p  { \color{black}{)}}  }} \FF \cong \RKLR $ where $  \RKLR $ is the original
cyclotomic KLR-algebra. Note however, that the degree function for
$  \RKLR $, {\color{black}{given for example in \cite{brundan-klesc}, \cite{KhovanovLauda}
    and \cite{Rouquier},}}
does {\it not} induce a $ \Z$-grading on $\RKLRZ $, since the relations in 
\eqref{finalrelations} are not homogeneous.

\medskip
We have already alluded to the following Theorem, that was proved by Hu and Mathas in
\cite{hu-mathas2}. 
\begin{theorem}\label{humathas}
There is an isomorphism of $ \Z_{(p)}$-algebras $F: \RKLRZ   \cong  \Z_{(p)}\Si_n    $. 
\end{theorem}

We next recall the diagrammatics for $ \RKLRZ$, as given in \cite{hu-mathas2}.
It is an 
extension of the diagrammatics for $ \RKLR$. 
A {\it KLR-diagram} $ D $ for $ \RKLR$ consists of $ n $ strands connecting $ n $ northern
points with $n$ southern 
points of a(n invisible) rectangle. Crossings are allowed in $ D $, but only crossings involving two strands.
Isotopic diagrams are considered to be equal. 
The strands of $ D $ are decorated with elements of $ \FF$, and the segments
of a strand are decorated with a nonnegative number of dots.
The product $ D_1 D_2 $ of KLR-diagrams $ D_1 $ and $ D_2 $ is realized by vertical concatenation
with $ D_1 $ on top of $D_2 $ where $ D_1 D_2 $ is set to zero if the bottom residue sequence
for $ D_1 $ does not coincide with the top
residue sequence for $ D_2 $. 
Here is an example of a KLR-diagram, using $ n= 6$
and $  p = 3$. 
\begin{equation}
  \raisebox{-.45\height}{\includegraphics[scale=0.7]{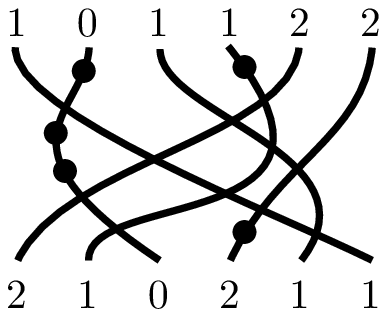}}
\end{equation}
The diagrammatics for $ \RKLRZ$ is given by

\begin{equation}
 \begin{split}
    &e(\ii) \mapsto  \raisebox{-.5\height}{\includegraphics[scale=0.8]{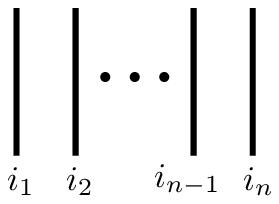}}, \, \, \, \, \, \, 
    y_l e(\ii) \mapsto  \raisebox{-.5\height}{\includegraphics[scale=0.8]{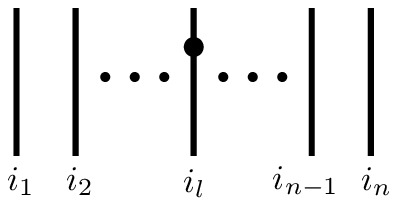}} \\
    &  \psi_k e(\ii) \mapsto  \raisebox{-.5\height}{\includegraphics[scale=0.8]{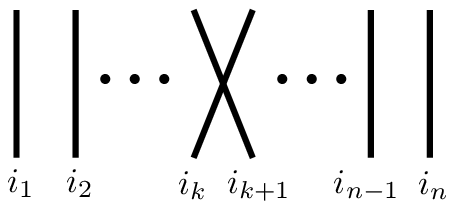}}    
 \end{split}
 \end{equation}
Via this, one can convert the relations \eqref{firstR} -- \eqref{finalrelations}
into a set of  diagrammatic relations for $ \RKLRZ $. 

\medskip

We now have the following Theorem which is an extension of
Theorem 3.2 and Remark 3.7 of \cite{PlazaRyom} to the integral case. 
\begin{theorem}\label{steendavid}
  Let $ n \ge 3$. If $ p > 3 $ then
the homomorphism $ \Phi $ from Lemma \ref{wellknownfunda} induces 
an isomorphism between $ \TLnZp $ and the quotient of $ \RKLRZ  $ given by the relation
\begin{align}\label{rel1}
  e(\ii)  & = 0      &    \mbox{ if } &  i_1  = 0 \mbox{\rm \, mod }  p, \, \, \, 
  i_2  =1 \mbox{\rm \, mod }  p  \mbox{ and }
  i_3  = 2 \mbox{\rm \, mod }  p 
\end{align}  
If $ p=3 $ then $ \Phi $ induces an isomorphism between $ \TLnZp $ and
the quotient of $ \RKLRZ \! \!  $
given by the relation
\begin{align}\label{rel2}
  y_3 e(\ii)  & = 0      &    \mbox{ if } &  i_1  = 0 \mbox{\rm \, mod }  p, \, \, \, 
  i_2  = 1 \mbox{\rm \, mod }  p  \mbox{ and }
  i_3  = 2 \mbox{\rm \, mod }  p 
\end{align}  
\end{theorem}
\begin{dem}
  The proof from \cite{PlazaRyom} carries over.
  It uses properties of Murphy's standard basis 
  that also hold in the present case. These properties lead to a description of 
  $ {\rm ker }\,  \psi $ as the ideal in $  \RKLRZ \! $, given by \eqref{rel1} and
  \eqref{rel2}.
\end{dem}

\medskip
We need the basic ingredients in Hu-Mathas' proof of \ref{humathas}, in the special
case $ \Z_{(p)}  \Si_n$ that we are considering.

\medskip
Let $\{ x_{\s \T}^{\lambda} \, | \, ( \s, \T) \in \std(\lambda)^{\times 2}, \lambda \in \Par \} $
be the specialization $ q= 1 $ of Murphy's standard basis for the Hecke algebra of type $ A_n$,
see \cite{Mat} and \cite{Murphy1}. As already mentioned in the proof of 
Theorem 
\ref{followingkeyresult}, 
it is a cellular basis for $ \Z_{(p)}  \Si_n $ on poset $ (\Par,
\trianglelefteq) $, and the elements $ \{ L_1, L_2, \ldots, L_n \} $ defined in \eqref{defJM}
form a family of $ \JM$-elements for $ \Z_{(p)} \Si_n $ with respect 
to the content function defined in \eqref{contentdef}. 
For $ \QQ \Si_n $, these $ \JM$-elements 
are separating, and so for $ \T \in \std(\lambda) $ we have an idempotent $ E_{\T} \in \QQ \Si_n$,
using the formula in \eqref{IdempotentHecke1}. For $ \s, \T \in \std(\Par) $ we define
\begin{equation} 
f_{\s \T}:= E_{ \s} x_{ \s \T } E_{\T} \in \QQ \Si_n
\end {equation}
Then the elements $ \{ f_{\s \T} \, | \, (\s, \T) \in \std(\lambda)^2, \lambda \in \Par \} $ form
a $ \QQ$-basis for $ \QQ \Si_n$.

\medskip
For $ \lambda \in \ParTwo $ and 
$ \s, \T \in \std(\lambda) $ we define similarly elements $ f_{\s \T} $ in $  \TLnQ$, denoted the
same way, via
\begin{equation}\label{denoted the same way}
f_{\s \T}:= \EE_{ \s} C^{\lambda}_{ \s \T } \EE_{\T} \in \TLnQ
\end {equation}
that form 
a $ \QQ$-basis for $ \TLnQ$. For $ \Phi: \QQ \Si_n \rightarrow \TLnQ $ the homomorphism
from Lemma \ref{wellknownfunda} we have that
$
\Phi(x_{\s \T}^{\lambda} ) = C_{\s \T}^{\lambda} + \mbox{ higher terms}  
$, 
where the higher terms are a linear combination of $ C_{\s_1 \T_1} $ with
$ \s_1 \rhd \s $ and $ \T_1 \rhd \T $, see Theorem 9 of \cite{Harterich}. Using this,
and that $ \Phi(L_i) = \LL_i$ and therefore $ \Phi(E_\T) = \EE_\T $ for $ \T \in \std(\ParTwo) $ we get that 
\begin{equation}\label{compatibility}
  \Phi(f_{\s \T}) = f_{\s \T} \, \, \, \, \, \, \mbox{ for } \s, \T \in \std(\ParTwo)
\end{equation}  

\medskip
For $ p $ a prime and $ \T \in \std(\Par)$ we define the $p$-class $ [ \T] \subseteq \std(\Par) $, 
as in \eqref{firstclass}. There is a well-defined function
from $p$-classes to residue sequences, given by $ [ \T] \mapsto \ii^{\T} := (c_\T(1), c_\T(2), \ldots, c_\T(n) )$.

\medskip
In the proof of the isomorphism {\color{black}{in}} Theorem \ref{humathas}, Hu and Mathas construct 
left and right actions of 
$e(\ii) $, $ y_l $ and $\psi_k  $ 
on $   \Z_{(p)}  \Si_n   $, by defining their actions
on $ \{ f_{\s \T } \}$. Let us explain the
formulas that they used for this.

\medskip
The formulas for $  e(\ii) $ are the simplest. They are given by 
\begin{align}\label{fofei}
&  e(\ii)  f_{\s \T  } :=
\begin{cases} 
    f_{\s \T }  &\mbox{ if } \ii^{\s} = \ii \\
  0 & \mbox{ if } \ii^{\s} \neq \ii
\end{cases}   
&   f_{\s \T  }  e(\ii) :=
\begin{cases} 
   f_{\s \T }  &\mbox{ if } \ii^{\T} = \ii \\
  0 & \mbox{ if } \ii^{\T} \neq \ii
\end{cases}
\end{align}

\medskip
The formulas for $  y_l  $ correspond to taking the nilpotent part
of the $ \JM$-element $ L_i$, just as in the proof of the original isomorphism
Theorem. For $ i \in \Z $  
let $ \hat{i} \in  \Z$ be given via integer division such that 
$ 0 \le  \hat{i}  \le p-1$ and $ \hat{i} \equiv i \, \, {\rm mod }\, \, p   $ and 
consider $ \hat{i} $ as an element of $ \Z_{(p)} $. 
Then 
\begin{align}\label{fofeiJM}
& y_l  f_{\s \T  }:=
     \left(c_{\s}(l) - \widehat{c_{\s}(l)}\right)  f_{\s \T  }   &
  f_{\s \T  } y_l :=
     \left(c_{\T}(l) - \widehat{c_{\T}(l)}\right)  f_{\s \T  }  
\end{align}

The formulas for $  \psi_k $ are a bit more complicated, but also the
most important for us.

\medskip
For $ \s \in \std(\lambda) $ where $ \lambda \in \Par$
and $ k =1,2, \ldots, n-1 $ we set $ \T := \s s_k $ and
$ r= r_\s(k) := c_{\s}(k) - c_{\T}(k)  $. We then 
define $\alpha =  \alpha_\s(k) \in \QQ$ via
\begin{align}\label{alphaF}
  &   \alpha_\s(k) :=
\begin{cases}  
  1 & \text{ if } \T\in \std(\lambda) \text{ and } \T \lhd \s \\
  \dfrac{r^2-1}
       {r^2} & \text{ if } \T\in \std(\lambda) \text{ and } \T \rhd \s \\
0 & \text{ otherwise } 
\end{cases}    
\end{align}  
In the terminology of \cite{hu-mathas2}, 
$ \alpha_\s(k)  $ is a choice of a {\it seminormal coefficient system}.
It is the
{\color{black}{\lq canonical choice\rq\ }}of a seminormal coefficient system, 
since it corresponds to the
{\color{black}{\lq non-diagonal part\rq\ }}of YSF,
see Corollary \ref{YSFsecond}.

\medskip
In order to define the action of $  \psi_k   $ it is enough to define
the left action of 
$ \psi_k e(\ii)    $
and the right action of $  e(\ii) \psi_k    $. 
Suppose that $ \ii^\s=(i_1, i_2,\ldots, i_k , i_{k+1}, \ldots, i_n) $. 
We first define  
$ \beta= \beta_\s(k) \in \QQ $ and $ {\color{black}{ {\widetilde{\beta}}}}=
{\color{black}{ {\widetilde{\beta}}}}_\s(k) \in \QQ $ via
\begin{align}\label{thendefinebeta}
&     \beta_{\s}(k) := 
\begin{cases}  
  \dfrac{\alpha}{1-r}  & \text{ if } i_k \equiv i_{k+1}   \, \,     {\rm mod } \, \, p \\
 \alpha r  & \text{ if } i_k \equiv i_{k+1} +1    \, \,     {\rm mod } \, \,  p \\
 \dfrac{\alpha r}{1-r} & \text{ otherwise } 
\end{cases} 
& {\color{black}{ {\widetilde{\beta}}}}_{\s}(k) := 
\begin{cases}  
  \dfrac{\alpha}{1+r}  & \text{ if } i_k \equiv i_{k+1}   \, \,     {\rm mod } \, \, p \\
- \alpha r  & \text{ if } i_k \equiv i_{k+1} -1    \, \,     {\rm mod } \, \,  p \\
 -\dfrac{\alpha r}{1+r} & \text{ otherwise } 
\end{cases}
\end{align}
Let $ \aaa \in \std(\lambda)$. 
We then have 
\begin{align}
&   \psi_k e(\ii) f_{\s \aaa  } :=
\begin{cases} 
    \beta  f_{\T \aaa  }   - \delta_{i_k , i_{k+1}} \dfrac{1}{r} f_{\s \aaa  } &\mbox{ if } \ii^{\s} = \ii \\
  0 & \mbox{ if } \ii^{\s} \neq \ii
\end{cases} \label{fofei2}    \\
& f_{\aaa \s   }   e(\ii) \psi_k :=
\begin{cases} 
   {\color{black}{ {\widetilde{\beta}}}} f_{\aaa \T   }   - \delta_{i_k , i_{k+1}} \dfrac{1}{r} f_{\aaa \s   } &\mbox{ if } \ii^{\s} = \ii \\
  0 & \mbox{ if } \ii^{\s} \neq \ii
\end{cases} \label{fofei2B}
\end{align}

The formulas in  \eqref{fofei} -- \eqref{fofei2B} are a key
ingredient in Hu and Mathas' proof of Theorem \ref{humathas}, see
Lemma 4.23 in \cite{hu-mathas2}. Note that 
the formulas \eqref{fofei} -- \eqref{fofei2B} in fact over-determine $ F(e(\ii) ),  F(y_l  ) $ and $ F( \psi_k ) $,
since already the left action on the basis $ \{ f_{\s \T} \} $ is enough to determine 
$ F(e(\ii) ),  F(y_l  ) $ and $ F( \psi_k ) $. In other words, the left action determines the right action and
vice versa. 

\medskip

We now return to the homomorphism $ \Phi: \Z_{ (p)} \Si_n \rightarrow \TLnZp$ from Lemma
\ref{wellknownfunda}. We have the following compatibility Theorem. 
\begin{theorem}
 The actions of $ \Phi( e(\ii)) $,  $ \Phi( y_l ) $ and  $ \Phi(  \psi_k ) $ 
are given by the formulas in 
\eqref{fofei} -- \eqref{fofei2},
with the only difference that $ f_{\s \T} $ is now the element of $ \TLnQ$ defined in \eqref{compatibility}.
\end{theorem}
\begin{dem}
  This is an immediate consequence of \eqref{compatibility} and
  the definitions in \eqref{fofei} -- \eqref{fofei2}. 
\end{dem}

\section{Seminormal form for $  \ee   \TLnF   \ee$}\label{Seminormal form for}
We write for simplicity $ \ee := \EE_{ [ \T_{\color{black}{n}}]} \in \TLnZp$, that is 
$ \ee := \Phi( e(\ii) ) $ where $ \ii=(0,-1, -2, \ldots, -n+1 ) $ is the decreasing residue sequence.
This is an idempotent in $ \TLnZp$ and so we obtain an 
{\it idempotent truncated} subalgebra $ \ee   \TLnZp   \ee$ of $ \TLnZp $.
This subalgebra plays an important role for what follows. 
To a certain extent, this runs parallel
to several recent papers, for example \cite{LiPl} and \cite{LPR}, where 
similar idempotent 
truncated algebras have been studied. 
By general principles, 
$ \ee   \TLnZp   \ee$ 
is a subalgebra of $ \TLnZp $, but with 
{\color{black}{unit}}-element $ \ee $.

\medskip
Under the isomorphism from Theorem \ref{steendavid}, the elements
$ \ee   \TLnZp   \ee$ are linear combinations of KLR-diagrams that have top
and bottom residue sequences both equal to
$ \ii = (0,-1, -2, \ldots, -n+1 ) $. 

\medskip
Recall from \eqref{kindofrecursive} that we have fixed natural numbers
$ \NNone $, $ \NNtwo $ and $ \RR $ such that $ \NN= (p-1) +\NNone $ and
$ \NNone= p\NNtwo +\RR $. As in Lemma \ref{kindofrecursive} 
we furthermore have blocks $ B_1, B_2. \ldots, B_{\NNtwo} $ of length $ p$ of consecutive natural numbers. 
The largest number of $ B_i $ is $ I := (i+1)p-1 $ and we define $ S_i \in \Si_n $ as 
\begin{equation}\label{reducedexpSi}
  S_i := s_I (s_{I-1}s_{I+1}) \cdots (s_{I-p+1} s_{I-p+3} \cdots s_{I+p-3} s_{I+p-1}) \cdots
 (s_{I-1}s_{I+1})    s_I 
\end{equation}
$ S_{\color{black}{i}} $ is a reduced expression
for the element of $ \Si_n $ that interchanges the blocks $ B_i $ and $ B_{i+1} $,
respecting the orders of the elements of each block.
We then define $ \UUU_i $ as the element of $ e(\ii ) \RKLRZ \!  e(\ii ) $
that is obtained from $ S_i $ by converting each $ s_j $ to
$ \psi_j $, and finally multiplying on the left and on the right by $ \ee$.
Similar elements have been considered before in \cite{KMR}, \cite{LiPl} and \cite{LPR},
but only for the original KLR-algebra $ \RKLR$ defined over a field. In \cite{LiPl} and \cite{LPR}, 
the $ \UUU_i$'s are called {\it diamonds}. For example, for $ n = 14 $ and $ p=3 $ we have
\begin{align}\label{diamonds} 
& \UUU_1 = \raisebox{-.5\height}{\includegraphics[scale=0.75]{dib76.eps}} & 
& \UUU_2 = \raisebox{-.5\height}{\includegraphics[scale=0.75]{dib77.eps}} 
\end{align} 

\medskip 
Our goal is to describe the left and right actions
on the $ \QQ $-basis $ \{f_{\s \T } \} $ for $ \TLnQ $. 
For this we have the following surprising Theorem, which may be
viewed as 
a generalization of Theorem \ref{YSFfirst}, and then also of
Corollary \ref{YSFsecond},
that is YSF, to the non-semisimple setting. As we shall see, its proof
relies on \eqref{fofei} and \eqref{fofei2}, and so
ultimately on Hu and Mathas' proof of the isomorphism Theorem 
\ref{humathas}. It is valid for $ \ee   \TLnZp   \ee$ and $ \ee   \TLnF   \ee$.

\begin{theorem}\phantomsection\label{YSFthird}
 Suppose that $ \s , \aaa \in [\T_\NN] \cap \std(\ParTwo)$, 
 and that $ i=1,2, \ldots, \NNtwo -1$. Let $ \T := \s \cdot S_i $ and
 suppose that $ \T $ is a standard tableau. 
If  $ \s \rhd \T $
set $ \s_u := \s $, otherwise 
set $ \s_u := \T   $. Let $ \s_d = \s_u \cdot S_i$. 
In the notation of Lemma \ref{kindofrecursive},
define $ f $ as $ f_1 $ if $ \RR= 0 $, otherwise as the first component of $ f_2 $. 
Define $ \greekrho := c_{ f( \s_u )}(i ) - c_{ f( \s_u )}(i+1 )  $ and $ X \in \QQ$ via 
\begin{equation}
  X:= \dfrac{  \Bigl( (\greekrho+1)p -1 \Bigr) \Bigl( (\greekrho+1)p -2 \Bigr)  \cdots \Bigl( \greekrho p +1\Bigr) }
  { \Bigl(\greekrho p -1 \Bigr) \Bigl(\greekrho p -2\Bigr)  \cdots \Bigl( (\greekrho-1)p +1\Bigr) }
\end{equation}
with $ p-1 $ factors in decreasing order in numerator as well as denominator. 
Then the left action of $ \UUU_i  $ is given by 
  \begin{description}
  \item[a)]
$ \UUU_i  f_{ \s_d  \aaa} = 
\dfrac{ \greekrho +1 }{ \greekrho }
f_{\s_d \aaa } + \dfrac{ \greekrho^2-1}{X \greekrho^2} f_{\s_u \aaa }$
  \item[b)]
$ \UUU_i f_{\s_u  \aaa }  = 
\dfrac{ \greekrho -1 }{ \greekrho }
f_{ \s_u \aaa } + X f_{ \s_d \aaa}$
    \end{description}
Suppose next that $ \T  $ is not standard. Then $  \UUU_i  $ acts via 
  \begin{description}
  \item[c)]
$ \UUU_i f_{ \s \aaa} = 
 0  \,\, \, \,\, \, \, \, \, \,  \mbox{ if } i, i+1 \mbox{ are in the same column of } f(\s) $ 
 \item[d)] $ \UUU_i f_{ \s \aaa} = 
 2 f_{ \s \aaa}  \, \mbox{ if } i, i+1 \mbox{ are in the same row of } f(\s) $ 
    \end{description}
\end{theorem}
\begin{dem}
  Let us first prove $ {\bf b) } $. The proof is a book-keeping of the coefficients
  that arise from the applications via \eqref {fofei2} of the $ \psi_i $'s that appear in   
  $ \UUU_i$. By the assumptions, in $\s_u $ the block of numbers $ B_i$ is positioned above the
  block of numbers $ B_{i+1} $
  as indicated below. 
\begin{equation}\label{diamond} 
  \s_u= \, \, \,  \raisebox{-.5\height}{\includegraphics[scale=0.7]{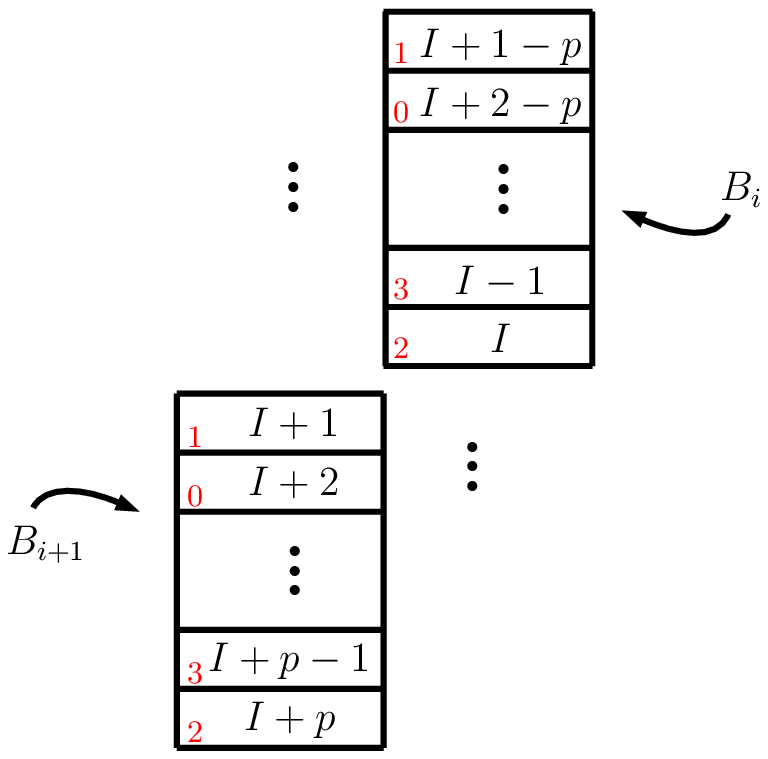}}\, \, \, 
S_i= \,  \raisebox{-.49\height}{\includegraphics[scale=0.7]{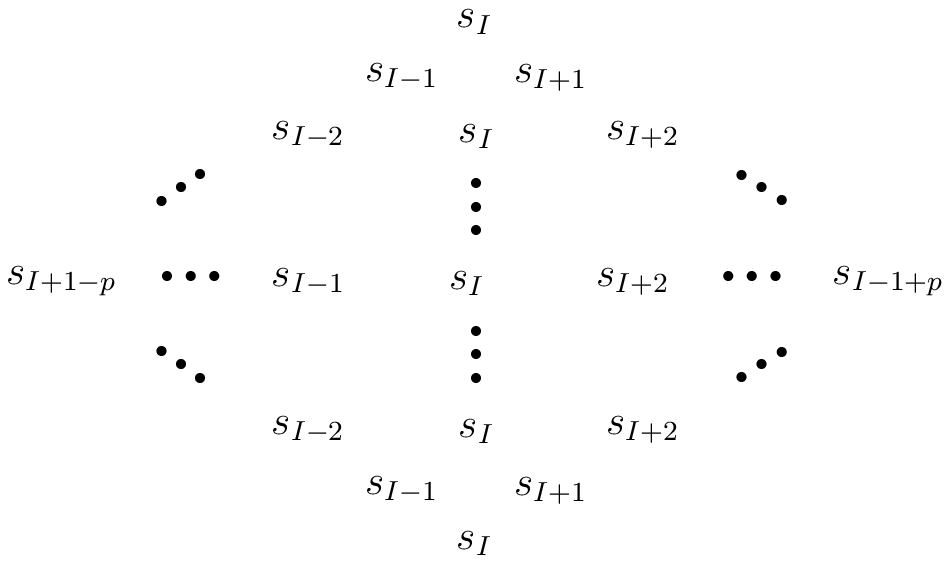}}
\end{equation}

\medskip
For simplicity we write $ f_{\s_u} = f_{ \s_u \aaa} $ and $ f_{\s_d} = f_{ \s_d \aaa} $. 
We first claim that $ \UUU_i  $ maps 
$ f_{\s_u }   $ to a linear combination of $  f_{\s_u} $ and $ f_{\s_d } $, disregarding the coefficients
for the time being. 

\medskip
To show this claim we proceed as follows. When applying $ \psi_I $ to $ f_{\s_u} $,
corresponding to the top row in the
diamond for $ S_i $ in \eqref{diamond}, the residue difference {\color{black}{modulo $p$}}
is $ 1 $, as can be read off
from the red numbers in \eqref{diamond}, 
and so by \eqref{fofei2} the result is a scalar multiple of $  f_{\s_u \cdot s_I} $, that is one term.
Next when applying $ \psi_{I+1} $ and $ \psi_{I-1} $ to $ f_{\s_u \cdot s_I } $,
corresponding to the second row in the 
diamond for $ S_i $ in \eqref{diamond}, the residue difference is $ 2 $
and so by \eqref{fofei2} the result is a multiple of $  f_{\s_u \cdot (s_I s_{I-1} s_{I+1})}  $, that is
one term once again. This pattern repeats itself until 
we reach the middle row of the diamond where the  
residue differences are all $ p $, and so by \eqref{fofei2} these $ \psi_i $'s 
produce two terms each, corresponding to the
two terms in \eqref{fofei2}. The tableau of the first term is given by the action by $ s_i $
whereas the tableau of the second is given by the omission of $ s_i $. 
On the other hand, the $ \psi_i $'s in the lower part of the diamond once again only
produce one term each. This pattern of residue differences can be read off from the KLR-diamonds as well,
see \eqref{diamonds}.

\medskip
We conclude from this that $ \UUU_i $ maps $ f_{\s_u} $ to a linear combination of
$  f_{\s_u \cdot \sigma} $ where $ \sigma $ is a subexpression of $ S_i $ obtained from $ S_i $ by deleting
certain of the $ s_i $'s from the middle row of $ S_i $ and
where $ \s_u \cdot \sigma $ is standard. If $ \sigma $ is the subexpression obtained by 
deleting all the $ s_i $'s of the 
middle row, the resulting term is $   f_{\s_u \cdot \sigma} =  f_{\s_u }  $ and if no $ s_i $ is deleted
the resulting term is $   f_{\s_u \cdot \sigma} =  f_{\s_d }  $, of course.

\medskip
We must however also consider the {\it mixed} cases where some of the $ s_i$'s from
the middle row of $ S_i$ are deleted, but not all. In these cases
we may use Coxeter relations to move a generator
$ s_i \neq s_I $ to the top of $ S_i $ and so we deduce that $ \s_u \cdot \sigma $ is not standard.
In Figure \ref{fig:here is ok} we give an example, using $ p = 5$, and
the indicated tableau $ \s_u$.

\begin{figure}[h]
$  \s_u= \, \, \, $ \raisebox{-.5\height}{\includegraphics[scale=0.6]{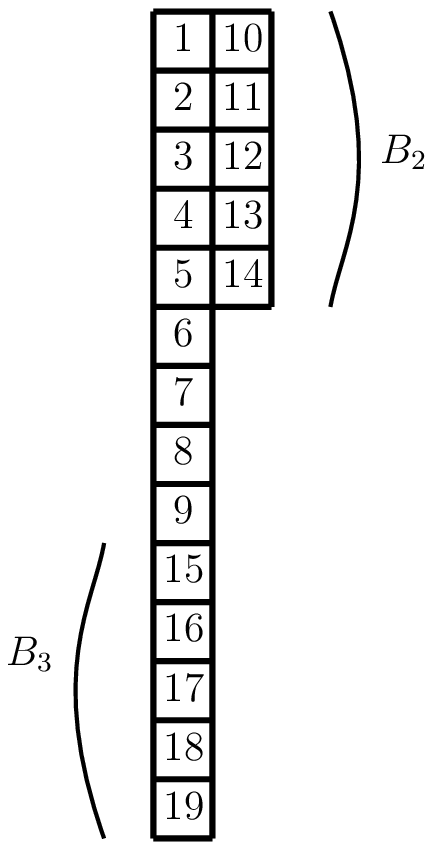}} \, \,
      \raisebox{-.5\height}{\includegraphics[scale=0.6]{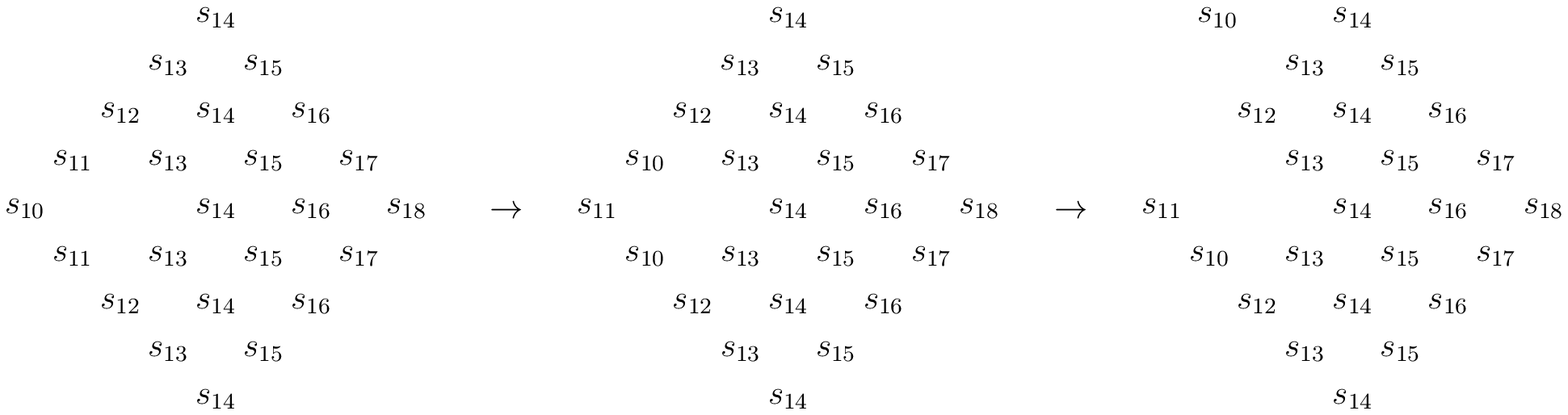}} 
\centering
\caption{Example using $p=5$.}
\label{fig:here is ok}
\end{figure}

It follows from this observation that the part of the action of
$  \psi_{ \sigma}  $ on $ f_{\s_u} $ that gives rise to $ f_{ \s_u \cdot \sigma } $ 
must involve the third case of \eqref{alphaF}, for
at least one of the $ \psi_i$'s, since the other cases produce standard tableaux.
But then the result is zero, proving that the mixed cases do not contribute to the action of $ \UUU_i  $ 
and so the claim is proved.

\begin{figure}[h]
   \raisebox{-.5\height}{\includegraphics[scale=0.8]{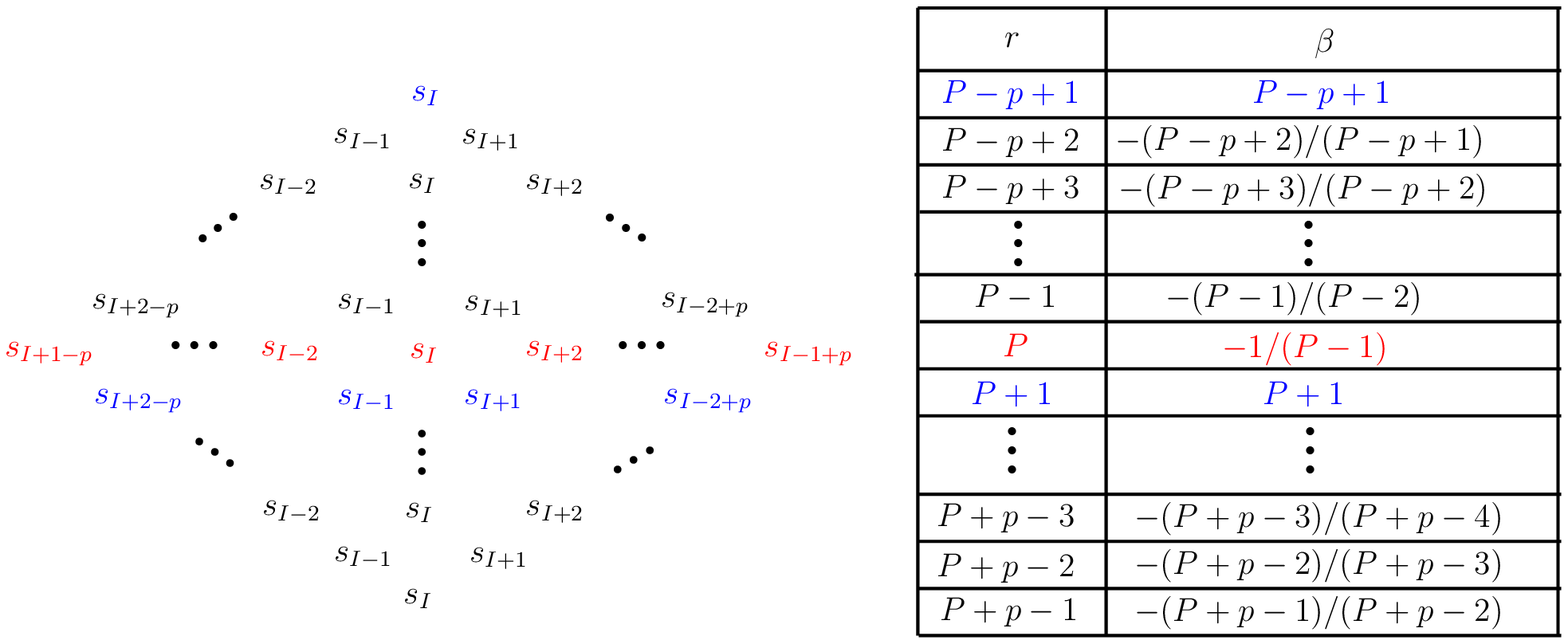}} \, \,  
\centering
\caption{Values of $ r $ and $ \beta $ for each row of the diamond.}
\label{fig:here is}
\end{figure}

\medskip
Let us now calculate the {\color{black}{coefficient}} of $ f_{\s_d}  $
under the action of $ \UUU_i  $ on $ f_{\s_u}$.
The contribution to this coefficient for each $ \psi_i $ of the middle row of the diamond is given by
always choosing the first term of 
\eqref{fofei2}. This implies that the coefficient of $ f_{\s_d}  $ always comes from
{\color{black}{\lq going down\rq\ }}and so
$\alpha = 1 $ for all occurrences of  
\eqref{alphaF} involved in the coefficient of $ f_{\s_d}  $.
The value of $ \beta $, according to \eqref{thendefinebeta}, therefore only depends on $ r $ and 
the relevant residue differences, that are constant along the rows of the diamond. 

\medskip

The table in Figure
\ref{fig:here is}
gives the values of $ r $ and $ \beta $ for each row of the diamond, 
where we write $ P := \greekrho p$, for simplicity. 
The colors in the table correspond to the three cases in the definition of $ \beta $ in
\eqref{thendefinebeta}, with red corresponding to the first case, blue to the second case and
black to the third case. 
To get the coefficient of $ f_{\s_d}  $ we must now multiply all the $ \beta$'s of the table in Figure
\ref{fig:here is}, with multiplicities given by the cardinalities of the rows of the diamond.

\medskip
We first claim 
that the sign of this product is $ +$. To show this
we observe that the number of black or red $ \beta$'s in the table in 
\eqref{thendefinebeta} is
$ p^2$ minus the number of blue $ \beta$'s, that is 
$ p^2 - p = p (p-1) $ which is even, proving the claim. 

\medskip
The product of the $ \beta$'s is therefore
\begin{align}
\begin{split}  
&  {\color{blue}{ \frac{ (P-p+1)} {1}}}  \frac{ (P-p+2)^2}{(P-p+1)^2}  \frac{ (P-p+3)^3}{(P-p+2)^3}
\cdots \frac{ (P-1)^{p-1}}{(P-2)^{p-1}}
       {\color{red}{ \frac{ 1} { (P-1)^p }  }}     {\color{blue}{ \frac{ (P+1)^{p-1} } { 1  }  }}
  \cdots 
 \frac{ (P+p-3)^3} {(P+p-4)^3 }  \frac{ (P+p-2)^2} {(P+p-3)^2 } \frac{ (P+p-1)} {(P+p-2) }  = \\ & \\
 &   {\color{blue}{ \frac{ \cancel{(P-p+1)}} {1}}}
 \frac{ \cancel{ (P-p+2)^2}}{(P-p+1)^{ \cancel{2}}}  \frac{ \cancel{(P-p+3)^3}}{(P-p+2)^{\cancel{3}}}
\cdots \frac{\cancel{ (P-1)^{p-1}}}{(P-2)^{ \cancel{p-1}}}
       {\color{red}{ \frac{ 1} { (P-1)^{\cancel{p }}}  }}     {\color{blue}{ \frac{ (P+1)^{\cancel{p-1}} } { 1  }  }}
  \cdots 
 \frac{ (P+p-3)^{\cancel{3}}} {\cancel{(P+p-4)^3 }}  \frac{ (P+p-2)^{\cancel{2}}} {\cancel{(P+p-3)^2} } \frac{ (P+p-1)} {\cancel{(P+p-2)} } =  
\\ & \\
 & \frac{ (P+1)} {(P-p+1)} \frac{ (P+2)} {(P-p+2)} \cdots
 \frac{ (P+p-2)} {(P-2) }  \frac{ (P+p-1)} {(P-1) } 
\end{split}
\end{align}
Remembering that $ P= \greekrho p $, we conclude from this that the coefficient of
$  f_{\s_d}  $ is 
$ X $ as claimed.

\medskip
In order to determine the coefficient of $  f_{\s_u}  $ we use the same method as for
the coefficient of $  f_{\s_d}  $, with the difference
that this time $ \alpha $ is
{\color{black}{\lq going down\rq\ }}only until reaching the middle row of diamond in which it 
{\color{black}{\lq stands still\rq\ }}and after this point, corresponding to the lower part
of the diamond, $ \alpha $ is
{\color{black}{\lq going up\rq\ }}again. 
Thus the table for
$  f_{\s_u}  $ coincides with the table in Figure \ref{fig:here is} in the upper half of the diamond, but
differs from it in the middle row and below. 
Using the definitions of $ r $, $ \alpha $ and $ \beta$, we then get the following table, where we use
the same color scheme as in Figure \ref{fig:here is}, and once again $ P := \greekrho p $.

\begin{equation}\label{here is once again} 
   \raisebox{-.5\height}{\includegraphics[scale=0.8]{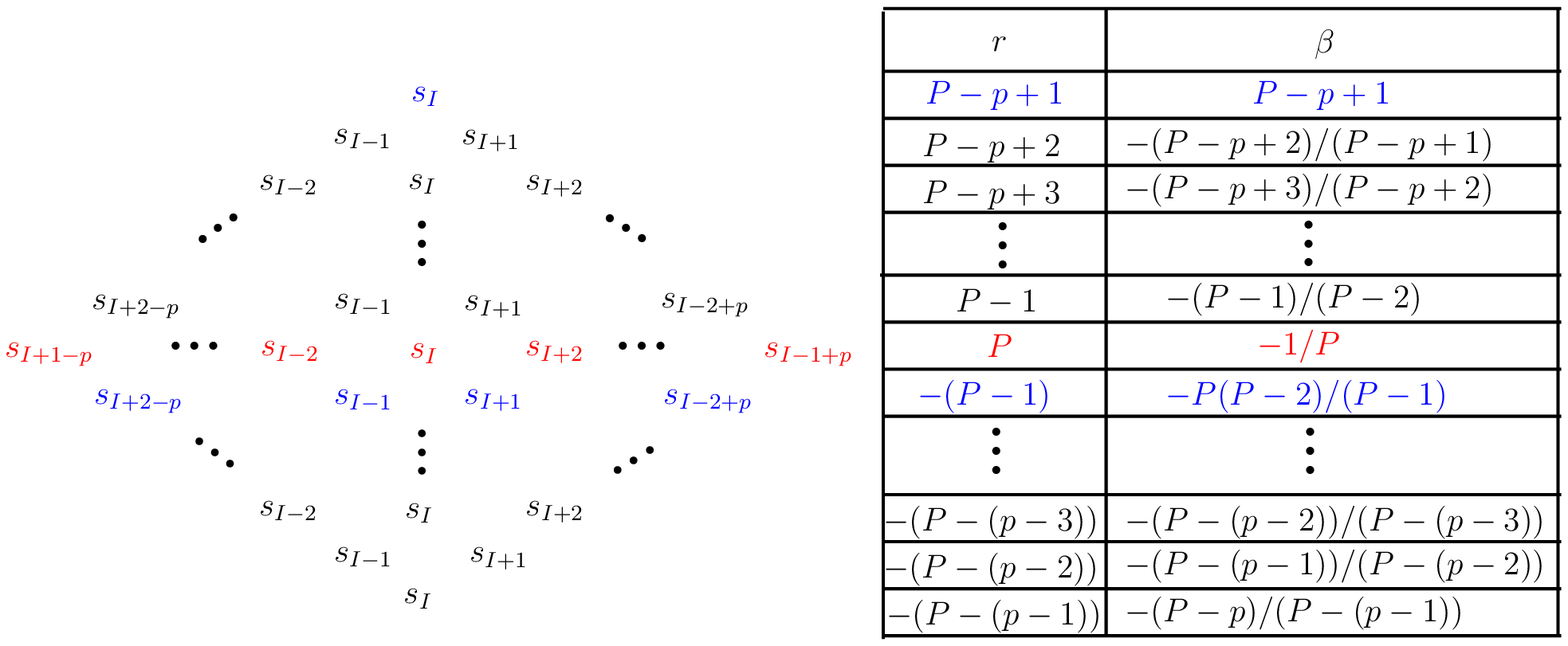}} \, \,
\end{equation}

We must calculate the product of the $\beta$'s that appear in 
\eqref{here is once again}. There is only one $ \beta$
appearing with a positive sign in \eqref{here is once again}, namely the one 
in the first row,  
and so the sign of the product of all the $ \beta$'s is $ (-1)^{ p^2 -1} = 1 $, since $ p $ is an odd
prime. It is now easy to calculate the product of the $ \beta$'s: indeed multiplying the $ \beta $
of the first row with the 
$\beta$ of the last row, the $ \beta$'s of the second row with the $ \beta$'s of the second last row, and so on, 
we find that the product of the $ \beta$'s is  
\begin{equation}
\dfrac{P-p}{P} = \dfrac{\greekrho p -p}{\greekrho p } = \dfrac{\greekrho  -1}{\greekrho} 
\end{equation}  
which proves $ {\bf b) } $

\medskip
The proof of $ {\bf a) } $ is proved with the same methods as the proof of $ {\bf b) } $
and is left to the reader.

{\color{black}{
\medskip
The proof of $ {\bf c) } $ is easy since, 
by the assumption for $ {\bf c) } $, all the numbers of $ B_i$ and $ B_{i+1} $ 
appear in the same column of $ \s $. In particular, $ I $ and $ I+1 $ appear in the same column of
$ \s $ and so indeed $ \UUU_i f_{ \s \aaa} = 0 $ since already $ \psi_I f_{ \s \aaa} = 0 $.

\medskip
The proof of $ {\bf d) } $ is slightly more complicated.
Under the assumption of $ {\bf d) } $, we have that $\s $ and $ \s \cdot S_i $ are as follows
\begin{equation}\label{proof of d} 
  \s= \, \, \,  \raisebox{-.5\height}{\includegraphics[scale=0.7]{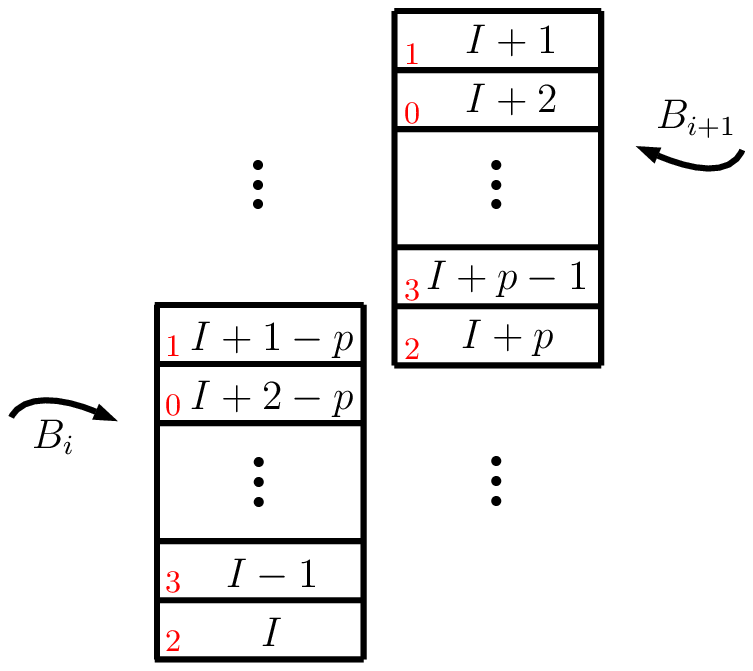}} \, \, \, \, \, \,  \, \, \, \, \, \, 
\s \cdot S_i= \,  \raisebox{-.49\height}{\includegraphics[scale=0.7]{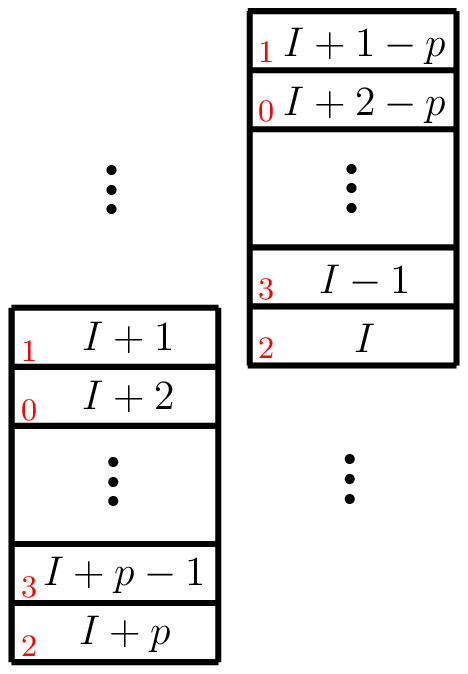}}
\end{equation}
with $ \s \cdot S_i $ non-standard. Using this, and arguing as in the paragraphs following \eqref{diamond},
we get that $ \UUU_i f_{ \s \aaa} =  \lambda f_{ \s \aaa} $ for some $ \lambda $.
With the same notation as before we then get the following table for calculating $ \lambda $.

\begin{equation}\label{here is once again A} 
   \raisebox{-.5\height}{\includegraphics[scale=0.8]{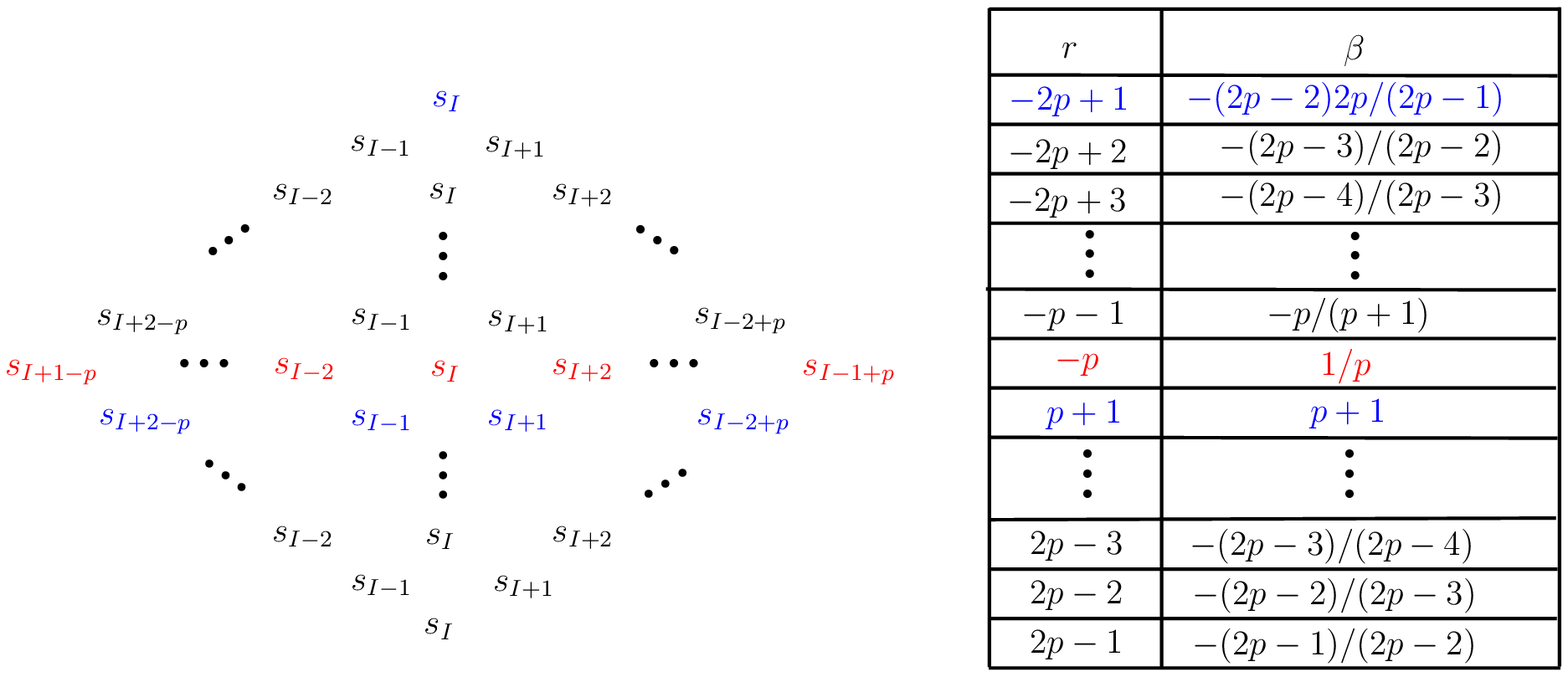}}
\end{equation}

The product of the $ \beta$'s of the table can be calculated by pairing the top $ \beta $ with the
bottom $ \beta$, and so on, and gives $ \lambda =  (-1)^{p-1} 2 = 2 $, 
proving $ {\bf d) } $.
(In fact, this calculation may be viewed as the calculation for the coefficient
of $ f_{\s_d \aaa }$ in $ {\bf a) } $ in the special case  
$  \greekrho =1 $). The proof of the Theorem is finished.
}}
\end{dem}

\medskip
We have the following variant of Theorem \ref{YSFthird} describing the right action
of $ \UUU_i $ on $ \{ f_{\s \T} \} $. Note that the formulas for the right action are the same as the formulas
for the left action, except that $ X $ should be replaced by $ \dfrac{1}{X}$. 

\begin{theorem}\phantomsection\label{YSFthirdB}
Let the notation be the same as in Theorem \ref{YSFthird}. 
Then the right action of $ \UUU_i  $ is given by 
  \begin{description}
  \item[a)]
$  f_{\aaa  \s_d  }  \UUU_i = 
\dfrac{ \greekrho +1 }{ \greekrho }
f_{\aaa \s_d  } + \dfrac{X (\greekrho^2-1)}{ \greekrho^2} f_{\aaa \s_u  }$
  \item[b)]
$  f_{\aaa \s_u   }  \UUU_i = 
\dfrac{ \greekrho -1 }{ \greekrho }
f_{ \aaa \s_u  } + \dfrac{1}{X} f_{\aaa  \s_d }$
    \end{description}
Suppose that $ \T  $ is not standard. Then $ \UUU_i  $ acts via 
  \begin{description}
  \item[c)]
$f_{ \aaa  \s}  \UUU_i = 
 0  \,\, \, \,\, \, \, \, \, \,  \mbox{ if } i, i+1 \mbox{ are in the same column of } f(\s) $ 
 \item[d)] $ f_{\aaa  \s } \UUU_i = 
 2 f_{\aaa  \s }  \, \mbox{ if } i, i+1 \mbox{ are in the same row of } f(\s) $ 
    \end{description}
\end{theorem}

\medskip
Statements similar to the one of the following Corollary, but for
the original KLR-algebra $ \RKLR$ defined over a field, are already present in literature,
see for example \cite{KMR}, \cite{LPR} and \cite{LiPl}, although the proofs in these references
are different from ours, since they rely on KLR-diagrammatics. 
\begin{corollary}\label{corKLRincl}
  Let $ \NNtwo $ be chosen as in \eqref{fixnotation}
{\color{black} and suppose that $ \NNtwo >1 $}. 
  Then there is a {\color{black}{(non-unital)}} injection of Temperley-Lieb algebras given by 
\begin{equation}\label{thereisaninjection}
  \iota_{KLR}:  \TLZpntwo \rightarrow
       {\mathbb {TL}}_{\NN}^{\! {\mathbb Z}_{(p)} }, \, \, \,  \UU_i \mapsto
{\color{black}{\Phi(\UU_i)}}
       \mbox{ for } i =1,2, \ldots, \NNtwo -1
\end{equation}
\end{corollary}

\begin{dem}
  We must show that the left action of the $\UUU_i$'s verify the Temperley-Lieb relations \eqref{eq oneTL},
\eqref{eq twoTL} and \eqref{eq threeTL}. The quadratic relation \eqref{eq oneTL} follows
  immediately from Theorem \ref{YSFthird}, since the $ 2 \times 2 $-matrix $ { \mathbf{M}_{\UUU_i}} $
  expressing the left action of $ \UUU_i $ in terms of $ \{ f_{\aaa \T_d }, f_{\aaa \T_u } \} $ has the form 
\begin{equation}
 \mathbf{M}_{\UUU_i} =  \begin{bmatrix}
\dfrac{\greekrho +1}{\greekrho} & X  \\ 
\dfrac{\greekrho^2 -1}{X\greekrho^{\color{black}{2}}} &   \dfrac{\greekrho -1}{\greekrho}
\end{bmatrix}
\end{equation}
which satisfies $  \mathbf{M}_{\UUU_i}^2 =  {\color{black}{2}}\mathbf{M}_{\UUU_i} $.

\medskip
In order to show relation \eqref{eq twoTL}, we choose $ \s , \T \in [\T_\NN] \cap \std(\ParTwo)$ 
and show that the left action of $ \UUU_i\UUU_{i \pm 1}\UUU_i $ on $f_{ \s \T} $ is equal to the left action
of $ \UUU_i $ on $f_{ \s \T} $. 
Let us focus on  $ \UUU_i\UUU_{i+ 1}\UUU_i $.
We then consider the positions of $ i, i+1 $ and $ i+2 $ in $ f(\s) $ where
$ f $ is as in Theorem \ref{YSFthird}. If $ i, i+1 $ and $ i+2 $ are in different rows of $ f(\s) $,
we have the following possibilities $ \s_1, \s_2, \ldots, \s_6 $ for $ f(\s)$.

\begin{equation}\label{hereA} 
  \raisebox{-.5\height}{\includegraphics[scale=0.65]{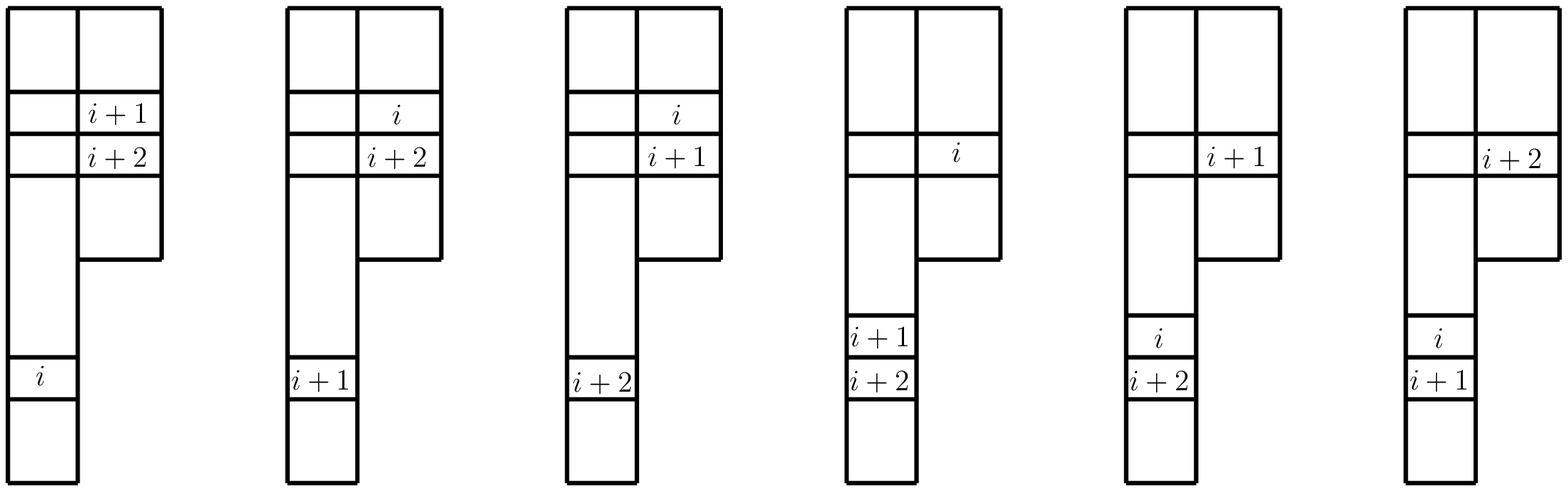}} \, \,
\end{equation}
One now checks for all $ j = 1,2,\ldots, 6 $ that indeed $  \UUU_i\UUU_{i+ 1}\UUU_i f_{\s_j \T}  = 
 \UUU_i f_{\s_j \T}  $. 
For example, using $ \greekrho := c_{\s_2}(i) - c_{\s_1}(i) $ one gets, using Theorem \ref{YSFthird}
repeatedly
\begin{align}
\begin{split}
  &   \UUU_i\UUU_{i+ 1}\UUU_i  f_{\s_1 \T} =
   \UUU_{i+ 1}\UUU_i  \left(\dfrac{\greekrho +1 }{\greekrho} f_{\s_1 \T}+
\dfrac{\greekrho^2 -1 }{X\greekrho^2}
f_{\s_2 \T}\right)   =
\dfrac{\greekrho^2 -1 }{X\greekrho^2}  \UUU_{i+ 1}\UUU_i
f_{\s_2 \T}    \\
& =   
\dfrac{\greekrho^2 -1 }{X\greekrho^2}  \UUU_i \left(  \dfrac{\greekrho  }{\greekrho-1} f_{ \s_2 \T} +
\dfrac{(\greekrho -1)^2 -1 }{X_1(\greekrho -1)^2} 
f_{  \s_3 \T} \right)   = \UUU_i  \left( \dfrac{\greekrho +1 }{X\greekrho}f_{ \s_2 \T} 
\right)   \\
& 
=  \dfrac{\greekrho +1 }{X\greekrho}  \left( \dfrac{\greekrho -1 }{\greekrho} f_{ \s_2  \T} + X f_{ \s_1 \T} \right )
=  \dfrac{\greekrho^2 -1 }{X\greekrho^2}   f_{ \s_2 \T} +  \dfrac{\greekrho +1 }{\greekrho} f_{ \s_1 \T} 
\end{split}
\end{align}
which equals $ \UUU_i  f_{ \s_1 \T}  $. For the other $ \s_j$'s, the verification of
$  \UUU_i\UUU_{i+ 1}\UUU_i f_{\s_j \T}  = 
 \UUU_i  f_{\s_j \T}  $. 
is
done the same way. 

\medskip
If two of the numbers $ i, i+1 $ and $ i+2$ are in the same 
row of $ f(\T)$ we have the following possibilities

\begin{equation}\label{hereAB} 
  \raisebox{-.5\height}{\includegraphics[scale=0.65]{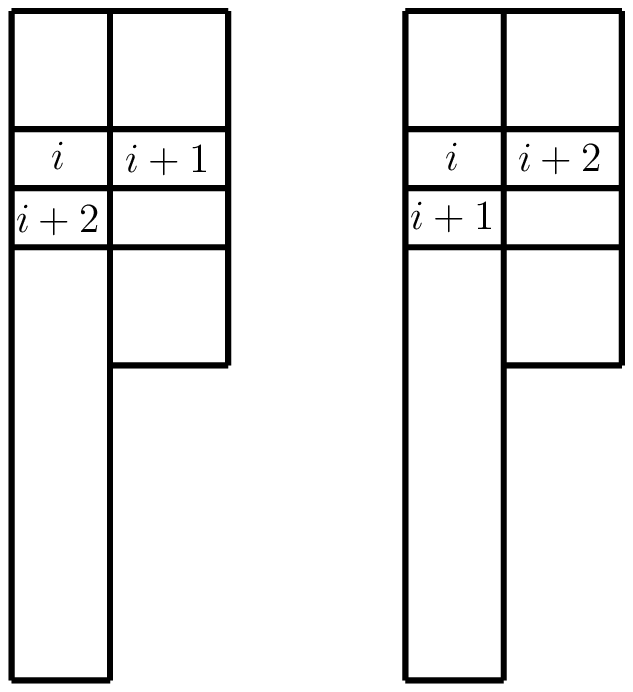}} \, \,
\end{equation}
and in each case one checks that $   \UUU_i\UUU_{i+ 1}\UUU_i  $ and
$   \UUU_i $ act the same way.
The verification of $ \UUU_i\UUU_{i-1 }\UUU_i  = \UUU_i $ is done the same way, 
and finally the verification of relation \eqref{eq threeTL} is trivial.

\medskip
In order to show injectivity of $ \iota_{KLR} $ one first checks 
that throughout the above arguments, one may always replace left actions by right actions.
(This also follows from the theory in 
\cite{hu-mathas2}).

\medskip
Let now $ \{ C_{ \s \T} \, |\,  \s, \T \in \std(\lambda), \lambda \in \ParTwoNuno \} $
be the basis for $\TLZpntwo $, as introduced in the paragraph before 
\eqref{diagrambasisTL}.
From the formulas in 
Theorem \ref{YSFthird} we have that $ \U, \V \in \std(\mu) $ with $ \mu \in \ParTwoNuno $ 
and $ C_{ \s \T}  f_{\U \V } \neq 0 $ implies $ \U \unrhd \T$, and similarly, from the formulas in
Theorem \ref{YSFthirdB}, we have that
$ f_{\U \V } C_{ \s \T}  \neq 0 $ implies $ \V \unrhd \s$.
Moreover, we also {\color{black}{have}} that $ C_{ \T^{\lambda} \U}  f^{}_{\U \V } = \mu^l_{ \U} f_{\T^{\lambda} \V}  $ where $ \mu^l_{  \U} \neq 0$
and that $ f^{}_{\U  \V } C_{ \V \T^{\lambda}} =  \mu^r_{ \V } f_{\U \T^{\lambda} }  $ where $  \mu^r_{ \V } \neq 0$
and where $ \U, \V $ are of shape $ \lambda$. 

\medskip
Suppose now that  
$ 0 \neq C=  \sum_{\s,\T} \lambda_{\s \T} C_{\s \T}
\in  {\rm ker} \, \iota_{KLR}  $. Choose $ (\s_0, \T_0 ) $ such that  
$  \lambda_{\s_0 \T_0} \neq 0 $ and such that
$ (\s_0, \T_0 ) $ is minimal with respect to this property. Then, using 
$ C f_{\s_0 \T_0 } =0  $ we get
$0=  f_{\T_0 \s_0 } C f_{\s_0 \T_0 } = \lambda_{ \s_0 \T_0}c \, \mu^l_{  \s_0} \mu^r_{ \T_0 } f_{\T_0 \T_0}   $, 
where $ c \neq 0 $, 
which implies $ \lambda_{ \s_0 \T_0}  =0 $. This is however 
a contradiction, and so the injectivity of $  \iota_{KLR} $ has been proved. 
\end{dem}

{\color{black}{
\begin{remark}\normalfont
For $ n_2 = 0 $ or $ n_2 = 1 $ the proof of Corollary \ref{corKLRincl} does not make sense.
          In these cases we define $ \iota_{KLR} $ by 
\begin{equation}\label{inthesecases}
  \iota_{KLR}:  \TLZpntwo \rightarrow {\mathbb {TL}}_{\NN}^{\! {\mathbb Z}_{(p)} }, \, \, \,  \one \mapsto
  \EEE_{[\T_{n}]}
\end{equation}
This definition corresponds to the basis case in the induction proof
of Theorem \ref{finalTheo} for all values of $n_2$.
\end{remark}
}}

{\color{black}{
\begin{remark}\normalfont
  In general $ \iota_{KLR}(\TLZpntwo) \subseteq \ee {\mathbb {TL}}_{\NN}^{\! {\mathbb Z}_{(p)} } \ee$,
  but this inclusion is not an equality, since for example $ \ee y_i \ee \in
  \ee {\mathbb {TL}}_{\NN}^{\! {\mathbb Z}_{(p)} } \ee \setminus  \iota_{KLR}(\TLZpntwo) $.
  Over $ \FF$ it is likely that $ \iota_{KLR}(\TLntwoF) $ is the degree zero part of
  $  \ee {\mathbb {TL}}_{\NN}^{\! \FF } \ee $.
\end{remark}
}}

\begin{remark}
  \normalfont
{\color{black}{Let $\TLZpntwo(2^p)$ be the Temperley-Lieb algebra defined over $ {\mathbb Z}_{(p)} $
  with loop parameter $2^p$ and
 let once again $ \NNtwo $ be chosen as in \eqref{fixnotation}. Then there is another
 injection $  \iota_{cab}:  \TLZpntwo(2^p) \rightarrow
       {\mathbb {TL}}_{\NN}^{\! {\mathbb Z}_{(p)} }$ 
       given by replacing each line after the first $p-1 $ lines by $p $ parallel lines.
       For example for $ \NN=14 $ and $ p=3 $ we have
 \begin{align}\label{diamondsref} 
&  \iota_{cab}(\UU_1) = \raisebox{-.5\height}{\includegraphics[scale=0.75]{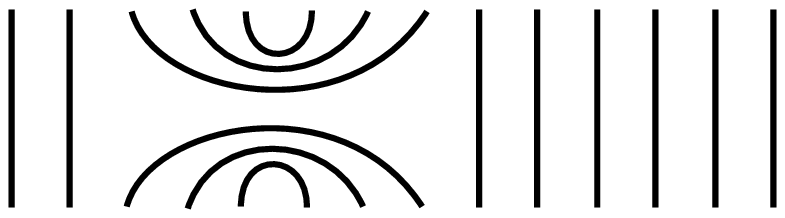}}, \,  & 
& \iota_{cab}(\UU_2) = \raisebox{-.5\height}{\includegraphics[scale=0.75]{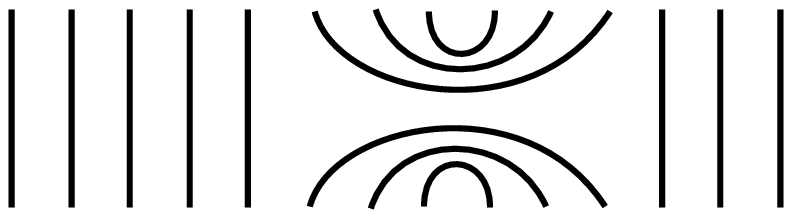}} 
 \end{align}
In view of Fermat's little Theorem, it induces an injection $ \iota_{cab}: \TLFpntwo \rightarrow 
{\mathbb {TL}}_{\NN}^{\! {\mathbb F}_{(p)} }$.
Note that $  \iota_{cab} $ is much simpler to define than $ \iota_{KLR} $ since it does not require
KLR-theory.

\medskip
Let $S_i \in \Si_n $ be as in
\eqref{reducedexpSi}. Then one gets an expression for $ \iota_{cab}(\UU_i) $ by replacing each $ s_j $ in $S_i $
by the generator $ \UU_j$ of $        {\mathbb {TL}}_{\NN}^{\! {\mathbb Z}_{(p)} }$. 
For example for $ \iota_{cab}(\UU_1) $ and $ \iota_{cab}(\UU_2) $ as in \eqref{diamondsref} one gets
 \begin{align}\label{diamondsrefex} 
  \iota_{cab}(\UU_1) = \UU_5  \UU_4  \UU_6 \UU_3  \UU_5  \UU_7  \UU_4 \UU_6 \UU_5,  &  & 
 \iota_{cab}(\UU_2) =\UU_8  \UU_7  \UU_9 \UU_6  \UU_8  \UU_{10}  \UU_7 \UU_9 \UU_8 
\end{align} 
This is parallel to our definition of $ \UUU_i $ in the paragraph following \eqref{reducedexpSi}
where we use $ \psi_i $'s instead of $ \UU_i$'s.
Via these expressions and Theorem \ref{YSFfirst} one may now attempt 
to describe the action of $ \iota_{cab}(\UU_i) $ on $ f_{ \s_d  \aaa}$, 
in the hope of finding formulas similar to the ones of Theorem \ref{YSFthird}, 
but 
already for small values
of $ \NN $ and $ p $ the result is an intractable linear combination of $ f_{\U \T} $'s
for $ \s_d \unlhd \U \unlhd \s_u $ where $ \s_d $ and $ \s_u $ are as in Theorem \ref{YSFthird}.
The reason for this is that YSF, that is Theorem \ref{YSFfirst}, gives rise to two $ f_{\U \T}  $
terms for each $ \UU_i $ 
in $ \iota_{cab}(\UU_i) $, whereas Hu and Mathas' formulas
\eqref{fofei2} and \eqref{fofei2B}
only give rise to one $ f_{\UU \T} $ term for each $ \psi_i $ in $ \UUU_i$, except
for the $ \psi_i $'s in the middle of the diamond. This simpler description of
the action of the $ \psi_i $'s, in comparison with the action of the $ \UU_i $'s, is a key ingredient
in the proofs of Theorem
\ref{YSFthird} and \ref{YSFthirdB} 
and it is a main reason why we need KLR-theory for  
Corollary \ref{corKLRincl} and therefore also, as we shall see, for the main results of this section.

\medskip
Over $ \FF $ it would be interesting to investigate whether the mentioned  linear combination of 
$ f_{\U \T} $'s reduces to the two terms $ f_{\s_d \T} $ and $ f_{\s_u \T} $ since that would imply that 
$ \iota_{cab}  $ and $ \iota_{KLR} $ coincide. In particular, $ \iota_{cab}(\UU_i)   $ would be homogeneous of degree $ 0 $,
although the individual $ \UU_j $-factors of $\iota_{KLR}(\UU_i) $ are not homogeneous.
We thank one of the referees for bringing $ \iota_{cab} $ to our attention.
}}
\end{remark}

{\color{black}{Suppose that $ n_2 > 1 $. }}
Let $ \{  \LLL_1, \LLL_2, \ldots, \LLL_{\NNtwo}  \} $ be the family 
of $ \JM$-elements in $ \TLZpntwo $ 
given by
$ \LLL_i := \Phi( L_i) $ where $ \{  L_1, L_2, \ldots, L_{\NNtwo}  \} 
\subseteq   \Z_{(p)} \Si_{\NNtwo}$ is the original family of $ \JM$-elements in 
\eqref{defJM} and where 
$\Phi: \Z_{(p)} \Si_{\NNtwo} \rightarrow  \TLZpntwo $
is the surjection from Lemma \ref{wellknownfunda}. 
Using the general theory in \cite{Mat-So}, 
we then obtain idempotents $ \EEE_{\T} \in \TLQpntwo$
for $ \T \in \ParTwoNuno$  
that are common eigenvectors for
the $   \LLL_i$'s, via the construction in \eqref{mainconstruction} and
Corollary \ref{finalcorsection2}.
On the other hand,
the inclusion $  \iota_{KLR}:  \TLZpntwo \rightarrow \TLZpn$ from Corollary
\ref{corKLRincl} induces an inclusion $\iota_{KLR}^{\QQ}: \TLQpntwo \subseteq \TLnQ $ and
so we may view the $  \EEE_{\T}$'s as idempotents in $ \TLnQ $ via $ \iota_{KLR}^{\QQ} $.

\medskip
Our next goal is to show, quite surprisingly, that
these new idempotents $ \{  \EEE_{\T} \, | \,  \T \in \std(\ParTwoNuno) $,
viewed as elements
in $ \TLnQ $,  
are closely related to the first idempotents $\{  \EE_{\T} \, | \, \T \in \std(\ParTwo) \}$
in $ \TLnQ $. We start with the following Lemma, which should be compared with
Lemma \ref{2.4}. 

\begin{lemma}\label{inductionbasis}
  Let $ \lambda \in \std(\ParTwo) $ and suppose that 
  $ \T= \T_{\lambda}  \in [\T_\NN] \cap \std(\ParTwo)$ and that
  $ \aaa \in  [\T_\NN] \cap \std(\lambda)$. 
Set $ \s:=  f(\T_{\lambda}) \in \ParTwoNuno $ where
$ f $ is as in Theorem \ref{YSFthird}. Let $f_{\T \aaa} $ and $f_{ \aaa \T} $
be as in \eqref{denoted the same way}. 
Then for for $ i =1,2,\ldots, \NNtwo $ we have that 
\begin{equation}
  \LLL_i f_{\T \aaa}  =  c_{ \s }(i) f_{\T \aaa} \, \, \, \, \,\, \,    \mbox{and} \, \, \, \, \, \, \,
f_{\color{black}{ \aaa \T}} \LLL_i  =  c_{ \s }(i) f_{\T \aaa}
  \end{equation}  
 \end{lemma}
\begin{dem}
Let us show the formula for the left action of $ \LLL_i $. 
Letting $ l_1 $ and $ l_2 $ be the column lengths of $ \s $ we have that
\begin{equation}\label{hereABC} 
\s=  \raisebox{-.5\height}{\includegraphics[scale=0.6]{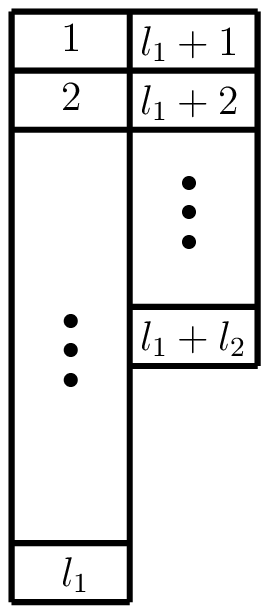}} \, \,
\end{equation}
Once again, we use the recursive formula 
$
  \LLL_{i+1} = (\UUU_i-\one) \LLL_i (\UUU_i-\one) +  \UUU_i-\one. 
$
Together with $ {\bf c) } $ of Theorem \ref{YSFthird}, it reduces the proof to the case $ l_2 = 1 $
and $ i = l_1+1$ where we must show that 
\begin{equation}\label{518}
\LLL_{l_1 +1} f_{\T {\color{black}\aaa}}   =  f_{\T {\color{black}\aaa}} 
\end{equation}
We do so by induction over $ l_1 $. 
The basis of the induction, corresponding to $ l_1 =1 $, is the affirmation that
$ \LLL_{2}  f_{\T \aaa}  =  f_{\T \aaa} \Longleftrightarrow  (\UUU_1- \one) f_{\T\aaa}  =  f_{\T \aaa} $ which is true
by $ {\bf d) } $ of Theorem \ref{YSFthird}.

\medskip
To show the inductive step
$ l_1 -1 \Longrightarrow l_1 $ we write for simplicity $ l:= l_1 $, 
$ \T_d := \T $ and $ \T_u := \T \cdot s_l $ and get via
$ {\bf a) } $ and $ {\bf b) } $ of Theorem \ref{YSFthird} that 
\begin{align}
\begin{split}
&       \LLL_{l +1} f_{\T \aaa} =  \Big( {\color{black}{(}}  \UUU_l-\one) \LLL_l (\UUU_l-\one) +  \UUU_l-\one \Big)  f_{\T_d \aaa} 
  = \LLL_l (\UUU_l-\one)  \Big(  \dfrac{1}{l} f_{\T_d \T} + \dfrac{l^2-1}{X\, l^2}  f_{\T_u \T}\Big) +
  \Big(  \dfrac{1}{l} f_{\T_d \aaa} + \dfrac{l^2-1}{X \, l^2}  f_{\T_u \aaa}\Big)  \\
& = (\UUU_l-\one) \Big(  \dfrac{1-l}{l} f_{\T_d \aaa} + \dfrac{l^2-1}{X\, l^2}  f_{\T_u \aaa}\Big)   +
  \Big(  \dfrac{1}{l} f_{\T_d \aaa} + \dfrac{l^2-1}{X \, l^2}  f_{\T_u \aaa}\Big)  =
 \UUU_l \Big(  \dfrac{1-l}{l} f_{\T_d \aaa} + \dfrac{l^2-1}{X\, l^2}  f_{\T_u \aaa}\Big)  + f_{\T_d \aaa} \\
 & =  \dfrac{1-l}{l} \UUU_l  \Big(   f_{\T_d \aaa} - \dfrac{l+1}{X\, l} f_{\T_u \aaa}\Big) +
 f_{\T_d {\color{black}{\aaa}}   } = f_{\T_d \aaa} =
 f_{\T \aaa}
\end{split}
\end{align}
The proof of the formula for the right action is done the same way. 

\end{dem}  

\medskip
The previous Lemma is the basis step
for the inductive proof of 
the following Theorem which should be compared with
Theorem \ref{mentionedabove}. 
\begin{theorem}\label{indutionseminormalKLR}
Suppose that 
$ \T, \aaa   \in [\T_\NN] \cap \std(\ParTwo)$. 
Set $ \s:=  f(\T) \in \ParTwoNuno $ where
$ f $ is as in Theorem \ref{YSFthird}.
Then for $ i =1,2,\ldots, \NNtwo  $ we have that 
\begin{equation}
  \LLL_i f_{\T \aaa }   = c_{ \s }(i) f_{\T \aaa } \, \, \, \, \,\, \,    \mbox{and} \, \, \, \, \, \, \, 
  f_{ \aaa \T} \LLL_i  = c_{ \s }(i) f_{ \aaa \T }   
\end{equation}
\end{theorem}
\begin{dem}
  As already indicated, the proof is by upwards induction over the dominance
  order in $ \std(\lambda) $, with Lemma \ref{inductionbasis} corresponding to the induction basis.
  The induction step is carried out the same way as the induction step in the proof of
  Theorem \ref{mentionedabove}, with Theorem \ref{YSFthird} replacing Theorem \ref{YSFsecond}. 
  The extra factors $ X$ or $ 1/X $ 
  in the equations corresponding to
 \eqref{YSFinduction}--\eqref{compareB} do not affect the conclusion. 
\end{dem}

\medskip
We have {\color{black}{a series}} of Corollaries to Theorem \ref{indutionseminormalKLR}. 
\begin{corollary}\label{idempotentJMcor}
  Let $ \T $ and $ \s $ be as in Theorem \ref{indutionseminormalKLR}
  and let $ \EE_{\T} \in \TLnQ$ be the idempotent from Corollary \ref{finalcorsection2}. 
  Then we have
\begin{equation}
 \LLL_i \EE_{\T}   = \EE_{\T}  \LLL_i = c_{ \s }(i) \EE_{\T} \, \mbox{ for } i =1,2,\ldots, \NNtwo 
\end{equation}
\end{corollary}
\begin{dem}
  This follows directly from
Theorem \ref{indutionseminormalKLR} together with 
  the construction of $ \EE_{\T}$ in \eqref{mainconstruction}
  and Corollary \ref{finalcorsection2}. 
\end{dem}

\begin{corollary}\label{idempotentJMcordos}
Suppose that 
$ \T  \in [\T_\NN] \cap \std(\ParTwo)$ and that 
$ \s \in \ParTwoNuno $ where $ \NNtwo $ is as in Lemma \ref{313}. 
 Let $\iota_{KLR}^{\QQ}: \TLQpntwo \subseteq \TLnQ $ be the inclusion
given by Corollary
\ref{corKLRincl}. Then we have
\begin{align}\label{622}
 &  \iota_{KLR}^{\QQ}  \big( \EEE_{\s} \big) \cdot \EE_{\T} =\EE_{\T} \,\cdot  \iota_{KLR}^{\QQ} \big( \EEE_{\s} \big) 
= \begin{cases}
 \EE_{\T} & \mbox{ if } f(\T) = \s \\
    0 & \text{ if } f(\T)  \neq \s 
\end{cases}
\end{align}
{\color{black}{In particular 
\begin{equation}\label{623}
  \iota_{KLR}^{\QQ}\big( \EEE_{\s} \big)  = \sum_{ \substack{ \T \in  [\T_\NN] \cap \std(\ParTwo)   \\ f(\T) = \s }}
    \EEE_{ \T}
\end{equation}}}
\end{corollary}
\begin{dem}
{\color{black}{To show \eqref{622}, we first suppose}}
  that $ f(\T) = \s$. Using 
  \eqref{IdempotentHecke1} and Corollary \ref{idempotentJMcor}
  we then get
\begin{equation} \iota_{KLR}^{\QQ}  \big( \EEE_{\s} \big) \cdot \EE_{\T} = 
  \left( \prod_{c \in {\cal C}} \, \, \prod_{\substack{i=1, \ldots, \NNtwo\\ c \neq c_{\s}(i)} } \dfrac{\LLL_i-c}{c_{\s}(i)-c}
  \right)   \EE_{\T} =
\left( \prod_{c \in {\cal C}} \, \, \prod_{\substack{i=1, \ldots, \NNtwo\\ c \neq c_{\s}(i)} } \dfrac{c_{\s}(i)-c}{c_{\s}(i)-c}
  \right)   \EE_{\T} =\EE_{\T} 
\end{equation}
as claimed.
Suppose next that $ f(\T) \neq \s$. Then there is $ i \in \{1,2,\ldots, \NNtwo \} $ such that
$ c_{ f(\T)}(i) \neq c_\s(i)$, since the separability condition \eqref{separationcondition}
is fulfilled, and so $  \iota_{KLR}^{\QQ}  \big( \EEE_{\s} \big)  $ has $ \big(\LLL_i - c_{ f(\T)}(i)  \big) $
as a factor. But by Corollary \ref{idempotentJMcor} we have
$ \big(\LLL_i - c_{ f(\T)}(i)  \big) E_{\T} = 0 $ which implies $ \iota_{KLR}^{\QQ}  \big( \EEE_{\s} \big) E_{\T} =0 $.
The formula for the right action in \eqref{622} is proved the same way.

{\color{black}{
Finally, \eqref{623} is a consequence of \eqref{622} 
since the $ \EEE_{\T} $'s are a complete set of orthogonal idempotents, see 
Corollary \ref{finalcorsection2}, and $ \iota_{KLR}^{\QQ}  \big( \EEE_{\s} \big) \EEE_{\U} =0$ for $ \U \in 
\std(\ParTwo) \setminus [\T_\NN] $.  }}
\end{dem}

\medskip

{\color{black}{

Let 
\begin{equation}
 n+1 = a_k p^k + a_{k-1} p^{k-1} + \ldots +a_1 p + a_0 
\end{equation}
be the expansion of $ n+1 $ in base $ p$ from \eqref{padic}. 
As in \eqref{fixnotation} we have $ n=n_1 + (p-1) $ and $ n_1 = p n_2 +r $
and so  
\begin{equation}\label{theprocess}
 r = a_0\, \,\mbox{ and }    n_2+1 = a_k p^{k-1} + a_{k-1} p^{k-2} + \ldots +a_1
\end{equation}
For our final 
Corollary we allow $n_2 $ to be any natural number or $0$. 
Let $ {\mathcal I}_n $ be the set defined in \eqref{mathcal I}.
\begin{corollary}\label{finalcor}
Choose $ \epsilon_i \in \{ \pm 1\} $ for $ i=1, 2, \ldots, k-1 $ and let 
$m=  (a_k p^{k-1} + \epsilon_{k-1}  a_{k-1} p^{k-2} + \ldots + \epsilon_{1} a_1 )-1  $ be the
corresponding element in $ {\mathcal I}_{n_2} $. Let 
$\iota_{KLR}^{\QQ}: \TLQpntwo \subseteq \TLnQ $ be as above. Then 
\begin{equation}\label{itfollowsfrom}
  \iota_{KLR}^{\QQ}( \EEE_{\T_{m}}) =
  \begin{cases}
\begin{aligned}
&  \EEE_{\T_{ (a_k p^{k} + \epsilon_{k-1}  a_{k-1} p^{k-1} + \ldots + \epsilon_{1} a_1 p +a_0 )-1}} +
 \EEE_{\T_{ (a_k p^{k} + \epsilon_{k-1}  a_{k-1} p^{k-1} + \ldots + \epsilon_{1} a_1 p -a_0 )-1}}
  & \mbox{ if } a_0 \neq 0  \\
& \EEE_{\T_{ (a_k p^{k} + \epsilon_{k-1}  a_{k-1} p^{k-1} + \ldots + \epsilon_{1} a_1 p )-1}}     & \text{ if } a_0 = 0 
\end{aligned}
\end{cases}
\end{equation}  
\end{corollary}
\begin{dem}
  If $ n_2 > 1 $ we get 
\eqref{itfollowsfrom} from \eqref{623} and the definition of $ f$, see Theorem \ref{YSFthird}
and Lemma \ref{kindofrecursive}. If $ n_2 =0 $ or $ n_2 =1 $ we get \eqref{itfollowsfrom}
directly from \eqref{inthesecases}. 
\end{dem}  
\medskip

We now finish the paper by showing how $ ^{p}\!\JWn  $ fits into the picture.
Recall that $ n \ge p $. Repeating the process in \eqref{theprocess} we find that $ n, n_1, n_2  $ and $ r $
belong to sequences of non-negative integers
$ n^i, n_1^i, n_2^i $ and $ r^i$ where 
$ n:=n^0, n_1=n_1^0, n_2=n_2^0 $ and $ r=r^0 $ and where 
\begin{equation}
n^i=n_i^i + (p-1), \,\, \,   n_1^i = p n_2^i +r^i, \, \, \, n^{i+1} = n_2^i \,  \mbox{  for }  i =0,1,\ldots, k-1 
\end{equation}
In fact we have 
\begin{equation}\label{infactwehave}
 r^i = a_i\, \,\mbox{ and }    n_2^i+1 = a_k p^{k-i-1} + a_{k-1} p^{k-i-2} + \ldots +a_{i+1}
\end{equation}
from which we see that $ n_2^i $ is strictly positive, except possibly $ n_2^{k-1} $ which may
be zero. 

\medskip
Using Corollary \ref{corKLRincl} we then get a chain of injections
\begin{equation}\label{thechain}
  {\mathbb {TL}}_{n^{k-1}_2}^{\! {\mathbb Z}_{(p)} } \subseteq
  {\mathbb {TL}}_{n^{k-2}_2}^{\! {\mathbb Z}_{(p)} } \subseteq  \cdots \subseteq  
  {\mathbb {TL}}_{n^{0}_2}^{\! {\mathbb Z}_{(p)} }  \subseteq  
  {\mathbb {TL}}\strut_{n}^{\! {\mathbb Z}_{(p)} } 
\end{equation}
By \eqref{infactwehave} we have $ n^{k-1}_2 = a_k -1 $ and so we have from 
\eqref{thechain} a (non-unital) injection
\begin{equation} \iota_k:  {\mathbb {TL}}_{a_k -1}^{\! {\mathbb Z}_{(p)} }
\subseteq   {\mathbb {TL}}\strut_{n}^{\! {\mathbb Z}_{(p)} }  
\end{equation}
With this we are in position to prove our final Theorem.
It  establishes the promised connection between the 
$p$-Jones-Wenzl idempotents and KLR-theory for the Temperley-Lieb algebra, via the seminormal form
approach to KLR-theory. 
\begin{theorem}\label{finalTheo}
In the above setting we have 
\begin{equation}\label{mustshoweq}
  ^{p}\!\JWn  =   \iota_k( \EEE_{\T_{(a_k -1)}}) 
\end{equation}  
\end{theorem}
\begin{dem}
We proceed by induction on $ k $.
If $ k=1 $ we have $ n+1 = a_1p +a_0 $ and so \eqref{mustshoweq} is the statement
\begin{equation}\label{basemain}
^{p}\!\JWn  = \iota_1( \EEE_{\T_{(a_1 -1)}}) 
\end{equation}  
But by Corollary \ref{finalcor} and the definitions both sides of \eqref{basemain} are equal to $ \EEE_{[\T_{(n])}}  $,
and so the basis of the induction is established.

\medskip
Let us now assume that \eqref{mustshoweq} holds for
$ k-1 $. Since $ n_2 = n_2+1 = (a_k p^{k-1} + a_{k-1} p^{k-2} + \ldots +a_1) -1 $ we then have
\begin{equation}\label{mainstep}
^{p}\!\JWntwo  = \iota_{k-1}( \EEE_{\T_{(a_k -1)}}) 
\end{equation}  
or equivalently 
\begin{equation}\label{mainstepA}
\sum_{ \epsilon_i \in \{\pm 1\}} \EEE_{\T_{(a_k p^{k-1} + \epsilon_{k-1} a_{k-1}p^{k-2} + \ldots + \epsilon_1 a_1) -1}}   = \iota_{k-1}( \EEE_{\T_{(a_k -1)}}) 
\end{equation}     
Applying $ \iota_{KLR}^{\QQ} $ to both sides of \eqref{mainstepA} we arrive via Corollary
\ref{finalcor} at 
\begin{equation}\label{mainstepAB}
  \sum_{ \epsilon_i \in \{\pm 1\}} \EEE_{\T_{(a_k p^{k} + \epsilon_{k-1} a_{k-1}p^{k-1} + \ldots + \epsilon_0 a_0) -1}}
  = \iota_k( \EEE_{\T_{(a_k -1)}}) 
\end{equation}     
that is $ ^{p}\!\JWn  = \iota_k( \EEE_{\T_{(a_k -1)}})  $, as claimed.
The Theorem is proved.

\end{dem}

\medskip
Viewing $ \EEE_{ [  \T_{n_2^i}   ]} $ as an element of
$  {\mathbb {TL}}\strut_{n}^{\! {\mathbb Z}_{(p)} }  $ via \eqref{thechain}, 
we can formulate Theorem \ref{finalTheo} as the statement
\begin{equation}
 ^{p}\!\JWn = \prod_{i=0}^{k-1} \EEE_{ [  \T_{n_2^i}   ]}
\end{equation}
since $ \EEE_{ [  \T_{n_2^i}   ]} \EEE_{ [  \T_{n_2^{k-1}}   ]} = \EEE_{ [  \T_{n_2^{k-1}}   ]} $.
In other words, $ \EEE_\T$ is a summand of $ ^{p}\!\JWn $ if and only
if $f^{(i)}(\T) := ({\overbrace{f \circ \ldots \circ f}^i}) (\T) \in [\T_{n_2^{i-1}}]$ 
for all $ i$.
{\color{black}{For example, for $ n= 12 $ and $ p= 3 $ we get 
using \eqref{tableauxclassA} and \eqref{tableauxclassAB} 
that $  \EEE_\T $ is a summand of $  ^{p}\!\JWdoce $ exactly for $ \T \in \{\T_{12}, \T_{10}, \T_6, \T_4 \} $
in the notation of \eqref{redcolortableaux1}. This is the precise meaning of our statement following \eqref{tableauxclassAB}.}}

\medskip
One could consider this as an incarnation of the fractal structure
of the representation theory of $   \TLnF $ or its Ringel dual $ SL_2(\FF) $,
studied 
for example \cite{E}, \cite{EH}, \cite{steen2}, \cite{TuWe1}.
}}

 \end{document}